\documentclass[10pt]{amsart}
\usepackage{mathrsfs, amsmath,amssymb,amsthm,amsfonts,graphicx}
  \graphicspath{{../../drawings/}}
  \usepackage[dvipsnames]{xcolor}
\usepackage{xcolor}

 \newcommand{\bsp}{\begin{split}}
  \newcommand{\esp}{{\end{split}}}

   \newcommand{\be}{\begin{equation}}

  \newcommand{\ee}{\end{equation}}
      
  \newcommand{\bee}{\begin{equation*}}
  \newcommand{\eee}{\end{equation*}}

 \newcommand{\ric}{\mathrm{Ric}}

    \newcommand{\cO}{\mathcal{O}}

\newtheorem{theorem}{Theorem}[section]
\newtheorem{proposition}[theorem]{Proposition}
\newtheorem{lemma}[theorem]{Lemma}
\newtheorem{corollary}[theorem]{Corollary}

\theoremstyle{definition}
\newtheorem{definition}[theorem]{Definition}
\newtheorem{remark}[theorem]{Remark}

\begin{document}

\title[Uniqueness of compact ancient solutions]{Uniqueness of compact ancient solutions to the higher dimensional Ricci flow}
\author[S.Brendle, P. Daskalopoulos, K. Naff, N. Sesum]{Simon Brendle, Panagiota Daskalopoulos, Keaton Naff  and Natasa Sesum}
\address{Department of Mathematics \\ Columbia University \\ New York NY 10027}
\address{Department of Mathematics \\ Columbia University \\ New York NY 10027}
\address{Department of Mathematics \\ Columbia University \\ New York NY 10027}
\address{Department of Mathematics \\ Rutgers University \\ Piscataway NJ 08854}
\begin{abstract}
In this paper, we study the classification of $\kappa$-noncollapsed ancient solutions to n-dimensional Ricci flow on $S^n$, extending the result in \cite{BDS} to higher dimensions. We prove that such a solution is either isometric to a family of shrinking round spheres, or the Type II ancient solution constructed by Perelman. 
\end{abstract}
\thanks{The first author was supported by the National Science Foundation under grant DMS-1806190 and by the Simons Foundation. The second author was supported by the National Science Foundation under grant DMS-1266172. The fourth author was supported by the National Science Foundation under grants DMS-1056387 and DMS-1811833.} 
\maketitle 

\tableofcontents

\section{Introduction}

In this paper, we consider a solution to the Ricci flow $\frac{\partial}{\partial t} g(t) = -2 \ric_{g(t)}$ on a compact manifold which exists for all times $t\in (-\infty,T)$. We call such a solution an ancient solution. The main focus of this paper is the classification of ancient solutions to Ricci flow in dimensions $n \geq 4$, under natural isotropic curvature conditions that will be discussed below. Ancient solutions play an important role in singularity formation in geometric flows since these solutions occur as limits of sequences of rescalings in regions of high curvature. For example, Perelman's work on the Ricci flow \cite{Perelman1} shows that high curvature regions in a three dimensional Ricci flow are modeled on ancient solutions with nonnegative curvature that are $\kappa$-noncollapsed. In the same paper, Perelman also showed that even in higher dimensions, ancient solutions that occur as blow-up limits around points of high curvature are $\kappa$-noncollapsed. We will focus on these $\kappa$-noncollapsed ancients solutions. Let us begin by briefly reviewing what is known in dimensions two and three. 

In dimension two, ancient solutions to the Ricci flow have been completely classified through a combination of work by Chu, the second author, Hamilton, and the fourth author in three papers \cite{Chu, Daskalopoulos-Sesum, Daskalopoulos-Hamilton-Sesum2}. In particular, there is actually a classification of both collapsed and $\kappa$-noncollapsed ancient solutions. Altogether, there are precisely three (non-flat, non-quotient) ancient solutions: the family of shrinking round spheres, the King solution, and steady cigar soliton. The King solution, independently discovered by King \cite{King} and Rosenau \cite{Rosenau}, resembles two steady cigar solitons which have been cut and glued together to form a compact solution. Of course, the sphere is $\kappa$-noncollapsed, while the cigar and, hence, the King solution are both collapsed. 

In dimension three, there are expected to be many more examples of \textit{collapsed} ancient solutions. For $\kappa$-\textit{noncollapsed} ancient solutions however, Perelman's conjecture \cite{Perelman1}, and its analogue in the compact setting, indicated a simple classification should exist. These conjectures stood for a number of years until several recent breakthroughs made it possible to resolve them in full. In dimension three, noncollapsed ancient solutions have now been completely classified through a combination of results by Angenent, and the first, second, and fourth authors in four papers \cite{Angenent-Brendle-Daskalopoulos-Sesum, Brendle1, Brendle2, BDS} (as well as a pinching result in \cite{Brendle-Huisken-Sinestrari}). Altogether, there are precisely four (non-flat, non-quotient) $\kappa$-noncollapsed ancient solutions: the family of shrinking round spheres, the family of shrinking round cylinders, Perelman's ancient oval solution on $S^3$, and the steady Bryant soliton. Perelman's ancient oval is the higher dimensional analogue of the King solution: it resembles a gluing of two Bryant solitons at very negative times. 

We now turn our attention to dimensions $n \geq 4$. We are interested in classifying ancient solutions which model singularity formation. The question is: what class of singularity models can we understand using the techniques developed in dimension three? The classification of ancient $\kappa$-solutions in dimension three relies on a number of ingredients. As we have mentioned, Perelman's $\kappa$-noncollapsing is crucial and holds in all dimensions. There are two ingredients special to dimension three though, which do not apply in higher dimensions. The first ingredient is the Hamilton-Ivey curvature pinching estimate, which ensures that all blow-up limits have nonnegative curvature. Once one has nonnegative curvature, Hamilton's Harnack inequality holds. The Harnack inequality then implies bounded curvature at bounded distance and, finally, an argument of Perelman upgrades this to bounded curvature. Therefore, in dimension three, ancient $\kappa$-solutions automatically have nonnegative and bounded curvature. No general Hamilton-Ivey-type estimate holds in higher dimensions without an initial curvature positivity assumption. The second ingredient is that the cross section of a noncompact singularity model in dimension three must be compact. Perelman used this property to establish an important ``tube and cap'' structure theorem for ancient $\kappa$-solutions in dimension three. When $n \geq 4$ however, there will be new singularity models, such as generalized cylinders, which will require new arguments to classify. 

The correct assumption on initial data, for the purpose of generalizing the singularity model classification in dimension three to higher dimensions, turns out to be positive isotropic curvature. It was Hamilton who first introduced positive isotropic curvature (PIC) to the Ricci flow in dimension four \cite{Hamilton3}. Importantly, Hamilton established an analogue of the Hamilton-Ivey pinching estimate for PIC initial data. Hamilton's result showed singularities must have nonnegative curvature and compact cross-section, whenever the models are noncompact. In \cite{Brendle3}, the first author generalized Hamilton's results for PIC initial data to dimensions $n \geq 12$. In these dimensions, PIC initial data does not ensure singularity models have nonnegative curvature, but rather they must satisfy a weaker curvature condition known as PIC2. See Section 2 to recall the precise definitions of PIC and PIC2. Importantly, however, this latter curvature condition is still strong enough to ensure Hamilton's Harnack inequality holds \cite{Bre09}. Hamilton and the first author's work justifies the following definition. 
\begin{definition}\label{kappa-solution}
Suppose $n \geq 4$. An ancient $\kappa$-solution is an $n$-dimensional, ancient, complete, nonflat solution of the Ricci flow that is uniformly PIC and weakly PIC2; has bounded curvature; and is $\kappa$-noncollapsed on all scales.
\end{definition}
To summarize, if the initial data of a Ricci flow is PIC and $n = 4$ or $n \geq 12$, then singularity models must be ancient $\kappa$-solutions in the sense above. We expect a similar result to be true for $5 \leq n \leq 11$.  

Our present goal is to complete the classification of ancient $\kappa$-solutions in the sense of Definition \ref{kappa-solution}, extending the classification in dimension three. Having identified the correct curvature assumptions, the program is roughly the same. The first important step was accomplished in \cite{Brendle1}, where the first author showed uniqueness of the Bryant soliton in the class of steady solitons with asymptotic cylindricality. Subsequently, the first author and the second author used \cite{Brendle1} and arguments in \cite{Brendle2} to prove uniqueness of the Bryant soliton among noncompact ancient $\kappa$-solutions in higher dimensions in \cite{BN}. It remains to extend the result of \cite{BDS} for compact ancient solution to higher dimensions, which we complete here. As in \cite{BDS}, the proof is accomplished in two steps. In the first step, we use arguments from \cite{BN} to prove:

\begin{theorem}
\label{rotational.symmetry}
Let $(S^n,g(t))$ be an ancient $\kappa$-solution on $S^n$. Then $(S^n,g(t))$ is rotationally symmetric.
\end{theorem}

Next, we give a complete classification of all ancient $\kappa$-solutions on $S^n$ with rotational symmetry: 

\begin{theorem}
\label{uniqueness.theorem}
Let $(S^n,g_1(t))$ and $(S^n,g_2(t))$ be two ancient $\kappa$-solutions on $S^n$ which are rotationally symmetric. Assume that neither $(S^n,g_1(t))$ nor $(S^n,g_2(t))$ is a family of shrinking round spheres. Then $(S^n,g_1(t))$ and $(S^n,g_2(t))$ coincide up to a reparametrization in space, a translation in time, and a parabolic rescaling.
\end{theorem}

Combining Theorem \ref{rotational.symmetry} and Theorem \ref{uniqueness.theorem}, we can draw the following conclusion:

\begin{theorem} 
Let $(S^n,g(t))$ be an ancient $\kappa$-solution on $S^n$ which is not a family of shrinking round spheres. Then $(S^n,g(t))$ coincides with Perelman's solution up to diffeomorphisms, translations in time, and parabolic rescalings.
\end{theorem}

Let us mention some related work in the mean curvature flow setting. In \cite{Daskalopoulos-Hamilton-Sesum1}, the authors classified compact, convex ancient solutions to the curve shortening flow. In \cite{Brendle-Choi1, Brendle-Choi2}, the authors proved that the bowl soliton is the only ancient solution which is noncompact, noncollapsed, strictly convex, and uniformly two-convex. In \cite{Angenent-Daskalopoulos-Sesum1}, the authors showed that every ancient solution which is compact, noncollapsed, strictly convex, and uniformly two-convex is either the family of shrinking spheres or the ancient oval constructed by White (cf. \cite{White}) and Haslhofer-Hershkovits (cf. \cite{Haslhofer-Hershkovits}). Finally, compact ancient solutions which are collapsed were studied in \cite{Bourni-Langford-Tinaglia}.

\smallskip 

The outline of the paper is as follows:  In Section \ref{structure}, we recall some qualitative properties of ancient $\kappa$-solutions on $S^n$. In particular, an ancient $\kappa$-solution on $S^n$ is either a family of shrinking round spheres, or it has the structure of two caps joined by a tube (in which the solution is nearly cylindrical). In Section \ref{proof.of.rotational.symmetry}, we give the proof of Theorem \ref{rotational.symmetry}. 

In Section \ref{a.priori.estimates}, we derive a-priori estimates for rotationally symmetric solutions. In Section \ref{weights}, we introduce two weight functions $\mu_+(\rho,\tau)$ and $\mu_-(\rho,\tau)$ (one for each cap). These will be used in Section \ref{difference.tip.region} to prove weighted estimates for the linearized equation in each tip region. 

In Section \ref{overview}, we give an overview of the proof of Theorem \ref{uniqueness.theorem}. The proof relies in a crucial way on estimates for the linearized equation in the tip region (Proposition \ref{estimate.for.difference.in.tip.region}) and in the cylindrical region (Proposition \ref{estimate.for.difference.in.cylindrical.region}). These estimates are proved in Section \ref{difference.tip.region} and Section \ref{difference.cylindrical.region}.

\section{Preliminary Results on Structure of Compact Ancient $\kappa$-Solutions}\label{structure}

In this section, we record basic facts about the structure of ancient $\kappa$-solutions as in Section 2 of \cite{BDS}. 
We begin by recalling the definitions of uniformly PIC and weakly PIC2 Riemannian manifolds. 
\begin{definition}\label{pic}
Suppose $n \geq 4$ and $(M, g)$ is a Riemannian manifold of dimension $n \geq 4$. 

\begin{itemize}
\item We say that $(M,g)$ is uniformly PIC if there exists a real number $\alpha > 0$ with the property that $R(\varphi,\bar{\varphi}) \geq \alpha \, |\text{\rm Rm}| \, |\varphi|^2$ for all complex two-vectors of the form $\varphi = (e_1+ie_2) \wedge (e_3+ie_4)$, where $\{e_1,e_2,e_3,e_4\}$ is an orthonormal four-frame. 

\item We say that $(M,g)$ is weakly PIC2 if $R(\varphi,\bar{\varphi}) \geq 0$ for all complex two-vectors of the form $\varphi = (e_1+i\mu e_2) \wedge (e_3+i\lambda e_4)$, where $\{e_1,e_2,e_3,e_4\}$ is an orthonormal four-frame and $\lambda,\mu \in [0,1]$. If the inequality is always strict, we say that $(M,g)$ is strictly PIC2.

\end{itemize}

Equivalently, 
\begin{itemize}
\item  $(M,g)$ is  uniformly PIC 
if there exists $\alpha > 0$ such that 
for every $p \in   M$  and every orthonormal four-frame $e_1, e_2, e_3, e_4 \in T_p M$, 
the curvature tensor $R_{ijkl}$ satisfies 
\[R_{1313} + R_{1414} + R_{2323}  + R_{2424}  -  2R_{1234}  \geq \alpha \, |\mathrm{Rm}|.\]

\item 
$(M,g)$ is  weakly PIC2 
if for every $p \in   M$, every orthonormal four-frame $e_1, e_2, e_3, e_4 \in T_p M$,  and every $\lambda, \mu \in [0,1]$,
the following holds
\[R_{1313} + \lambda^2 R_{1414} + \mu^2 R_{2323}  + \lambda^2 \mu^2 R_{2424}  -  2 \lambda \mu R_{1234}  \geq 0. \]

\end{itemize} 
\end{definition}

Four-dimensional ancient $\kappa$-solutions in the sense of Definition \ref{kappa-solution} automatically satisfy the restricted isotropic curvature pinching condition used to study four-dimensional ancient solutions in \cite{CZ06}. For a proof and to recall the meaning of the pinching condition, see Proposition A.2 in \cite{BN}. 
Note that the restricted isotropic curvature pinching condition implies four-dimensional ancient $\kappa$-solutions have nonnegative curvature operator. The restricted isotropic curvature condition is the assumption under which the authors in \cite{CZ06} developed a theory for ancient $\kappa$-solutions in dimension four, following the work of Hamilton and Perelman. Importantly, for each of the structure results established for ancient $\kappa$-solutions in dimensions $n \geq 5$ by the first author in \cite{Brendle3} under the uniformly PIC and weakly PIC2 assumptions, there is an analogous result for ancient $\kappa$-solutions in dimension $n = 4$ under the restricted isotropic curvature pinching condition, which can be found in \cite{CZ06}. In particular, we note that compactness of ancient $\kappa$-solutions in the sense above is established in \cite{Brendle3} for $n \geq 5$ and follows from work in \cite{CZ06} for $n=4$. 

\smallskip

From now on, we assume $(M, g(t))$ is an ancient $\kappa$-solution which is compact and simply connected. We also assume $(M, g(t))$ {\em is not a family of shrinking round spheres.} Note that because $(M,g(t))$ is compact, the strong maximum principle (Proposition 6.6 in \cite{Brendle3}) implies $(M, g(t))$ is strictly PIC2. By the work of the first author and Schoen, this implies $M$ is diffeomorphic to $S^n$. 

\begin{proposition}
The asymptotic shrinking soliton associated with $(M, g(t))$ is a cylinder.
\end{proposition}
\begin{proof}
The only gradient shrinking Ricci solitons which are uniformly PIC and weakly PIC2 are round the sphere $S^{n}$, the round cylinder $S^{n-1} \times \mathbb{R}$, or a quotient of one of these two by a discrete group of isometries (see Theorem A.1 in \cite{BN}). If the asymptotic soliton has constant curvature, then $(M, g(t))$ would have constant curvature by the pinching result in \cite{Brendle-Huisken-Sinestrari}. This would contradict our assumption that $(M, g(t))$ is not a family of shrinking round spheres. The asymptotic soliton cannot be a compact quotient of the cylinder for a number of reasons. Perhaps the clearest is that these compact quotients of $S^{n-1} \times \mathbb{R}$ do not move self-similarly under the Ricci flow. Alternatively, if the asymptotic soliton is compact, then by smooth Cheeger-Gromov convergence $M$ must be diffeomorphic to a compact quotient of the cylinder, but $M$ is diffeomorphic to $S^n$. Finally, if the asymptotic soliton is a noncompact quotient of the cylinder, then by Theorem A.1 in \cite{Brendle4} the fundamental group of some nontrivial quotient of $S^{n-1}$ would inject into the fundamental group of $M$, which is trivial by assumption. It follows the asymptotic soliton must be isometric to a round cylinder. 
\end{proof}

\begin{proposition}\label{kappa-compactness}
Let $(x_k, t_k)$ be an arbitrary sequence of points in space-time satisfying $\lim_{k \to \infty} t_k = -\infty$. Consider the family of rescaled metrics $g_k(t) := R(x_k , t_k) g\big(t_k +  R(x_k, t_k)^{-1}t\big)$. After passing to a subsequence, the sequence of pointed flows $(M, g_k(t), x_k)$ converges in the Cheeger-Gromov sense to either a family of shrinking cylinders or the Bryant soliton. 
\end{proposition}
\begin{proof}
We have compactness of ancient $\kappa$-solutions for $n \geq 4$ and a classification in the noncompact case by \cite{BN}. Thus, the proof in \cite{BDS} works here. 
\end{proof}

We now fix a large number $L < \infty$ and a small number $\varepsilon_1 > 0$ so that the Neck Improvement Theorem in \cite{BN} (Theorem 4.8) holds. Also, denoting by  $\lambda_1(x,t)$  the smallest eigenvalue of the Ricci tensor at $(x, t)$,  we  fix a small  number  $\theta >0$ with the property that if $(x, t)$ is a spacetime point satisfying $\lambda_1(x, t) \leq \theta R(x,t)$, then the point $(x, t)$ lies at the center of an evolving $\varepsilon_1$-neck.  The existence of $\theta$ is based on a standard contradiction argument which uses compactness of ancient $\kappa$-solutions in higher dimensions. See Lemma A.2 in \cite{BN} for a proof in the noncompact case. The proof in the compact case is nearly identical. 

\begin{proposition}\label{bryant-domains}
Consider a sequence of times $t_k \to -\infty$. If $k$ is sufficiently large, then we can find two disjoint compact domains $\Omega_{1, k}$ and $\Omega_{2,k}$ with the following properties: 
\begin{itemize}
\item $\Omega_{1, k}$ and $\Omega_{2,k}$ are each diffeomorphic to $B^n$.
\item For each $x \in M \setminus (\Omega_{1,k} \cup \Omega_{2,k})$, we have $\lambda_1(x,t_k) < \theta R(x, t_k)$. In particular, the point $(x, t_k)$ lies at the center of an evolving $\varepsilon_1$-neck. 
\item For each $x \in \Omega_{1,k} \cup \Omega_{2,k}$, we have $\lambda_1(x,t_k) > \frac{1}{2} \theta R(x, t_k)$. 
\item $\partial \Omega_{1, k}$ and $\partial \Omega_{2,k}$ are leaves of Hamilton's CMC foliation of $(M, g(t_k))$. 
\item For each $k$, there exists a leaf $\Sigma_k$ of Hamilton's CMC foliation with the property that $\Omega_{1,k}$ and $\Omega_{2,k}$ lie in different connected components of $M \setminus \Sigma_k$, and $\sup_{x \in \Sigma_k} \frac{\lambda_1(x, t_k)}{R(x, t_k)} \to 0$. 
\item The domains $(\Omega_{1,k}, g(t_k))$ and $(\Omega_{2,k}, g(t_k))$ each converge, after rescaling, to a corresponding subset of the Bryant soliton. 
\end{itemize}
\end{proposition}
\begin{proof}
The proof is the same as the proof of Proposition 2.3 in \cite{BDS}.
\end{proof}

\begin{corollary}
If $k$ is sufficiently large, then the set $\{x \in M : \lambda_1(x, t_k) > \frac{1}{2n} R(x,t_k), \, \nabla R(x, t_k) = 0 \}$ consists of exactly two points. One of these points lies in $\Omega_{1, k}$ and the other lies in $\Omega_{2,k}$. In particular, these points are contained in different connected components of $M \setminus \Sigma_k$. 
\end{corollary} 
\begin{proof}
The proof is the same as the proof of Corollary 2.4 in \cite{BDS}. We note that in higher dimensions, at the tip of the $n$-dimensional Bryant soliton we have $\mathrm{Ric} = \frac{1}{n} R\,g$,  $\nabla R = 0$, and $\nabla^2 R < 0$. Moreover, the tip is the unique point on the Bryant soliton where $\nabla R = 0$ and $\lambda_1 > \frac{1}{2n} R$. The claim now follows from the previous proposition. 
\end{proof}

\begin{definition} 
We say that $p$ is a tip of $(M, g(t))$ if $\lambda_1(p, t) > \frac{1}{2n} R(p, t)$ and $\nabla R (p, t) = 0$. 
\end{definition}

By the previous corollary, $(M, g(t))$ has exactly two tips $-t$ is sufficiently large, and the two tips are contained in different components of $\{x \in M : \lambda_1(x, t) > \frac{1}{2} \theta R(x, t)\}$. Let us denote the tips of $(M, g(t))$ by $p_{1, t}$ and $p_{2, t}$. The following proposition is an immediate consequence of Proposition \ref{kappa-compactness}.

\begin{proposition} \label{bryant-at-tips}
Consider a sequence of times $t_k \to -\infty$. Let $p_{1, t_k}$ and $p_{2, t_k}$ denote the tips in  $(M, g(t_k))$. If we rescale the flow around $(p_{1, t_k}, t_k)$ or $(p_{2, t_k}, t_k)$ as in Proposition \ref{kappa-compactness}, then the rescaled flows subsequentially converge to the Bryant soliton in the Cheeger-Gromov sense. 
\end{proposition}

\begin{proposition}
Consider a sequence of times $t_k \to -\infty$. Let $p_{1, t_k}$ and $p_{2, t_k}$ denote the tips in $(M, g(t_k))$. 
Then, we have  $R(p_{1, t_k}, t_k)d_{g(t_k)}(p_{1, t_k}, p_{2, t_k})^2 \to \infty$ and  $R(p_{2, t_k}, t_k)d_{g(t_k)}(p_{1, t_k}, p_{2, t_k})^2 \to \infty$. 
\end{proposition}
\begin{proof} 
We have Perelman's long-range curvature estimate for $n = 4$ by Proposition 3.6 in \cite{CZ06} and for $n \geq 5$ by Theorem 6.13 in \cite{Brendle3}. Thus, the proof of Proposition 2.7 in \cite{BDS} works here.
\end{proof}

\begin{corollary} 
For a given  positive real number $A$, if $-t$ is sufficiently large (depending on $A$), then the balls $B_{g(t)}(p_{1,t}, A R(p_{1, t}, t)^{-\frac{1}{2}})$ and $B_{g(t)}(p_{2,t}, A R(p_{2, t}, t)^{-\frac{1}{2}})$ are disjoint. 
\end{corollary}

\begin{proposition}
For  a sequence of spacetime points $(x_k, t_k)$  with  $t_k \to -\infty$  denote by  $p_{1, t_k}$ and $p_{2, t_k}$ 
 the tips of $(M, g(t_k))$. If $R(p_{1, t_k}, t_k)d_{g(t_k)}(p_{1, t_k},x_k)^2\to \infty$ and $R(p_{2, t_k}, t_k)d_{g(t_k)}(p_{2, t_k},x_k)^2\to \infty$, then $\frac{\lambda_1(x_k, t_k)}{R(x_k, t_k)}\to 0$. 
\end{proposition} 
\begin{proof}
The proof is the same as the proof of Proposition 2.9 in \cite{BDS} 
\end{proof}

\begin{corollary} \label{necks-outside-caps}
Given $\varepsilon > 0$, we can find a time $T \in (-\infty, 0]$ and a large constant $A$ with the following property. If $t \leq T$ and $x \not\in B_{g(t)}(p_{1, t}, A R(p_{1, t}, t)^{-\frac{1}{2}} )\cup  B_{g(t)}(p_{2, t}, A R(p_{2, t}, t)^{-\frac{1}{2}} )$, then $(x, t)$ lies at the center of an evolving $\varepsilon$-neck. 
\end{corollary}

\section{Rotational Symmetry of Compact Ancient $\kappa$-Solutions in Higher Dimensions}\label{proof.of.rotational.symmetry}

In this section, we give a proof of rotational symmetry, extending Section 3 in \cite{BDS}  to higher dimensions. The arguments are essentially the same, except we will use results from \cite{BN}, which is the higher dimensional analogue of part two of \cite{Brendle2}, where the first author first established rotational symmetry of noncompact ancient $\kappa$-solutions in dimension three. Throughout this section, we assume $n \geq 4$ and $(M, g(t))$ is an $n$-dimensional ancient $\kappa$-solution which is compact and simply connected. We also assume that $(M, g(t))$ is not a family of shrinking round spheres. The proof of rotational symmetry is by contradiction. Therefore: \textbf{We will assume throughout this section that $(M, g(t))$ is not rotationally symmetric.}

As in the previous section, let us fix a large number $L < \infty$ and small number $\varepsilon_1 > 0$ so that the Neck Improvement Theorem in \cite{BN} holds. Then let us choose a small number $\theta > 0$ so that if $\lambda_1(x, t) \leq \theta R(x,t)$, the spacetime point $(x, t)$ lies at the center of an evolving $\varepsilon_1$-neck.

We begin with a definition of $\varepsilon$-symmetry of the caps based on the definition used in the noncompact case in \cite{BN}. 
\begin{definition}[Symmetry of Caps]\label{cap-symmetry}
We will say the flow is $\varepsilon$-symmetric at time $\bar t$ if there exists a compact domain $D \subset M$ and a family of time-independent vector fields $\mathcal U = \{U^{(a)} : 1 \leq a \leq {n \choose 2}\}$ which are defined on an open subset containing $D$ such that the following statements hold: 
\begin{itemize}
\item The domain $D$ is a disjoint union of two domains $D_1$ and $D_2$, each of which is diffeomorphic to $B^n$.
\item $\lambda_1(x, \bar t) < \theta R(\bar x, \bar t)$ for all points $x \in M \setminus D$.
\item $\lambda_1(x, \bar t) > \frac{1}{2} \theta R(x, \bar t)$ for all points $x \in D$.
\item $\partial D_1$ and $\partial D_2$ are leaves of Hamilton's CMC foliation of $(M, g(\bar t))$. 
\item For each $x \in M \setminus D$, the point $(x, \bar t)$ is $\varepsilon$-symmetric in the sense of Definition 4.2 in \cite{BN}.
\item The Lie derivative $\mathcal L_{U^{(a)}}(g(t)))$ satisfies the estimate  
\[
\sup_{D_1 \times  [\bar t - \rho_1^2, \bar t]} \sum_{l = 0}^2 \sum_{a = 1}^{{n \choose 2}} \rho_1^{2l} \big|D^l (\mathcal L_{U^{(a)}}(g(t)))|^2 \leq \varepsilon^2,
\]
where $\rho_1^{-2} := \sup_{x \in D_1} R(x, \bar t)$. 
\item The Lie derivative $\mathcal L_{U^{(a)}}(g(t)))$ satisfies the estimate  
\[
\sup_{D_2 \times  [\bar t - \rho_2^2, \bar t]} \sum_{l = 0}^2 \sum_{a = 1}^{{n \choose 2}} \rho_2^{2l} \big|D^l (\mathcal L_{U^{(a)}}(g(t)))|^2 \leq \varepsilon^2,
\]
where $\rho_2^{-2} := \sup_{x \in D_2} R(x, \bar t)$. 
\item If $\Sigma \subset D_1$ is a leaf of Hamilton's CMC foliation of $(M, g(\bar t))$ that has distance at most $50\,r_{\mathrm{neck}}(\partial D_1)$ from $\partial D_1$, then 
\[
\sup_{\Sigma} \sum_{a = 1}^{{n \choose 2}} \rho_1^{-2}|\langle U^{(a)}, \nu \rangle|^2 \leq \varepsilon^2, 
\]
where $\nu$ is the unit normal vector to $\Sigma$ in $(M, g(\bar t))$ and $r_{\mathrm{neck}}(\partial D_1)$ is defined by the identity $\mathrm{area}_{g(\bar t)}(\partial D_1) = \mathrm{area}_{g_{S^{n-1}}}(S^{n-1}) r_{\mathrm{neck}}(\partial D_1)^{n-1}$. 
\item If $\Sigma \subset D_2$ is a leaf of Hamilton's CMC foliation of $(M, g(\bar t))$ that has distance at most $50\,r_{\mathrm{neck}}(\partial D_2)$ from $\partial D_2$, then 
\[
\sup_{\Sigma} \sum_{a = 1}^{{n \choose 2}} \rho_2^{-2}|\langle U^{(a)}, \nu \rangle|^2 \leq \varepsilon^2, 
\]
where $\nu$ is the unit normal vector to $\Sigma$ in $(M, g(\bar t))$ and $r_{\mathrm{neck}}(\partial D_2)$ is defined by the identity $\mathrm{area}_{g(\bar t)}(\partial D_2) = \mathrm{area}_{g_{S^{n-1}}}(S^{n-1}) r_{\mathrm{neck}}(\partial D_2)^{n-1}$. 
\item[$\bullet$] If $\Sigma \subset D_1$ is a leaf of Hamilton's CMC foliation of $(M, g(\bar t))$ that has distance at most $50 \,r_{\mathrm{neck}}(\partial D_1)$ from $\partial D_1$, then 
\[
\sum_{a, b = 1}^{{n \choose 2}} \bigg| \delta_{ab} - \mathrm{area}_{g(\bar t)}(\Sigma)^{-\frac{n+1}{n-1}} \int_\Sigma \langle U^{(a)}, U^{(b)} \rangle_{g(\bar t)} \, d\mu_{g(\bar t)} \bigg|^2 \leq \varepsilon^2. 
\]
\item[$\bullet$] If $\Sigma \subset D_2$ is a leaf of Hamilton's CMC foliation of $(M, g(\bar t))$ that has distance at most $50 \,r_{\mathrm{neck}}(\partial D_2)$ from $\partial D_2$, then 
\[
\sum_{a, b = 1}^{{n \choose 2}} \bigg| \delta_{ab} - \mathrm{area}_{g(\bar t)}(\Sigma)^{-\frac{n+1}{n-1}} \int_\Sigma \langle U^{(a)}, U^{(b)} \rangle_{g(\bar t)} \, d\mu_{g(\bar t)} \bigg|^2 \leq \varepsilon^2. 
\]
\end{itemize}
\end{definition}

By Proposition \ref{bryant-at-tips} and Corollary \ref{necks-outside-caps}, the solution is increasingly symmetric back in time. 

\begin{proposition}
Let $\varepsilon > 0$ be given. If $-t$ is sufficiently large (depending on $\varepsilon$), then the flow is $\varepsilon$-symmetric at time $t$. 
\end{proposition} 

\begin{remark} For $-\bar t$ sufficiently large, the flow is $\varepsilon$-symmetric and the tips of $(M, g(\bar t))$ are contained in different connected components of the domain $D$. In particular, after relabeling $D_1$ and $D_2$ if necessary, we may assume $p_{1, \bar t} \in D_1$ and $p_{2, \bar t} \in D_2$. In view of Corollary \ref{necks-outside-caps} and Definition \ref{cap-symmetry}, we have $\mathrm{diam}_{g(\bar t)}(D_1) \leq C R(p_{1, \bar t}, \bar t)^{-\frac{1}{2}}$ and $\mathrm{diam}_{g(\bar t)}(D_2) \leq C R(p_{2, \bar t}, \bar t)^{-\frac{1}{2}}$ for some constant $C$, which depends only on our choice of $\theta$. By the long-range curvature estimate, this implies $\frac{1}{C} R(p_{1, \bar t}, \bar t) \leq R(x, \bar t) \leq C R(p_{1, \bar t}, \bar t)$ for all $x \in D_1$ and $\frac{1}{C} R(p_{2, \bar t}, \bar t) \leq R(x, \bar t) \leq C R(p_{2, \bar t}, \bar t)$ for all $x \in D_2$.
\end{remark}

Recall by Lemma 9.5 in \cite{Brendle2} (cf. Lemma 5.5 in \cite{BN}): 

\begin{lemma}\label{continuity-of-symmetry}
Suppose that the flow is $\varepsilon$-symmetric at time $\bar t$. If $\tilde t$ is sufficiently close to $\bar t$, then the flow is $2\varepsilon$-symmetric at time $\bar t$. 
\end{lemma}

Now to proceed with the proof by contradiction, consider an arbitrary sequence of positive real numbers $\varepsilon_k \to 0$. For $k$ large, define 
\[
t_k := \inf\{t \in (-\infty, 0] \,:\, \text{The flow is not $\varepsilon_k$-symmetric at time $t$}\}. 
\]

We must have $\limsup_{k\to \infty}  t_k = -\infty$ since otherwise the flow would be symmetric for $-t$ sufficiently large, in contradiction with our assumption.

For sufficiently negative times, we denote by $p_{1, t}$ and $p_{2, t}$ the tips of $(M, g(t))$. Let  and $r_{2,k}^{-2} := R(p_{2,t_k} t_k)$. Since $t_k \to -\infty$, Proposition \ref{bryant-at-tips} implies that if we rescale the flow about either $p_{1, t_k}$ by the factor $r_{1,k}^{-2} := R(p_{1,t_k}, t_k)$ or $p_{2, t_k}$ by the factor $r_{2,k}^{-2} := R(p_{2,t_k}, t_k)$, then the sequence subsequentially converges to the Bryant soliton in the pointed Cheeger-Gromov sense. This gives us the analogue of Lemma 3.5 in \cite{BDS}. 

\begin{lemma}\label{bryant-approximation-at-tips} 
There exists a sequence of real numbers $\delta_k \to 0$ such that the following statements hold when $k$ is sufficiently large: 
\begin{itemize}
\item For each $t \in [t_k - \delta_k^{-1}r_{1,k}^2, t_k]$, we have $d_{g(t)}(p_{1, t_k}, p_t) \leq \delta_k r_{1, k}$ and  $1 -\delta_k \leq r_{1,k}^2 R(p_{1,t}, t) \leq 1 + \delta_k$.
\item For each $t \in [t_k - \delta_k^{-1}r_{2,k}^2, t_k]$, we have $d_{g(t)}(p_{2, t_k}, p_t) \leq \delta_k r_{2, k}$ and  $1 -\delta_k \leq r_{2,k}^2 R(p_{2,t}, t) \leq 1 + \delta_k$.
\item The scalar curvature satisfies $r_{1,k}^{2} R(x, t) \leq 4$ and
\[
\frac{1}{2K}(r_{1,k}^{-1} d_{g(t)}(p_{1,t_k}, x) + 1)^{-1} \leq r_{1,k}^2 R(x, t) \leq 2K( r_{1,k}^{-1} d_{g(t)}(p_{1,t_k}, x) + 1)^{-1}
\]
for all points $(x, t) \in B_{g(t_k)}(p_{1,t_k}, \delta_k^{-1} r_{1,k}) \times  [t_k - \delta_k^{-1}r_{1,k}^2, t_k]$. 
\item The scalar curvature satisfies $r_{2,k}^{2} R(x, t) \leq 4$ and
\[
\frac{1}{2K}(r_{2,k}^{-1} d_{g(t)}(p_{2,t_k}, x) + 1)^{-1} \leq r_{2,k}^2 R(x, t) \leq 2K( r_{2,k}^{-1} d_{g(t)}(p_{2,t_k}, x) + 1)^{-1}
\]
for all points $(x, t) \in B_{g(t_k)}(p_{2,t_k}, \delta_k^{-1} r_{2,k}) \times  [t_k - \delta_k^{-1}r_{2,k}^2, t_k]$. 
\item There exists a nonnegative function $f_1 := f_{1,k} : B_{g(t_k)}(p_{1,t_k}, \delta_k^{-1} r_{1,k}) \times  [t_k - \delta_k^{-1}r_{1,k}^2, t_k] \to \mathbb{R}$ such that $|\mathrm{Ric} - D^2 f_1 | \leq \delta_k r_{1,k}^{-2}$
and 
\[
 |\Delta f_1 + |\nabla f_1|^2 - r_{1,k}^{-2}| \leq \delta_k r_{1,k}^{-2} \quad\text{and} \quad  |\frac{\partial}{\partial t} f_1 + |\nabla f_1|^2 | \leq \delta_k r_{1,k}^{-2}.
\]
Moreover, function $f_1$ satisfies 
\[
\frac{1}{2K}(r_{1,k}^{-1} d_{g(t)}(p_{1,t_k}, x) + 1) \leq f_1(x, t) + 1 \leq 2K (r_{1,k}^{-1} d_{g(t)}(p_{1,t_k}, x) + 1)
\]
 for all points $(x, t) \in B_{g(t_k)}(p_{1,t_k}, \delta_k^{-1} r_{1,k}) \times  [t_k - \delta_k^{-1}r_{1,k}^2, t_k]$.
\item There exists a nonnegative function $f_2 := f_{2,k} : B_{g(t_k)}(p_{2,t_k}, \delta_k^{-1} r_{2,k}) \times  [t_k - \delta_k^{-1}r_{2,k}^2, t_k] \to \mathbb{R}$ such that $|\mathrm{Ric} - D^2 f_2 | \leq \delta_k r_{2,k}^{-2}$
and
\[
 |\Delta f_2 + |\nabla f_2|^2 - r_{2,k}^{-2}| \leq \delta_k r_{2,k}^{-2} \quad\text{and} \quad  |\frac{\partial}{\partial t} f_2 + |\nabla f_2|^2 | \leq \delta_k r_{2,k}^{-2}.
\]
Moreover, function $f_2$ satisfies 
\[
\frac{1}{2K}(r_{2,k}^{-1} d_{g(t)}(p_{2,t_k}, x) + 1) \leq f_2(x, t) + 1 \leq 2K (r_{2,k}^{-1} d_{g(t)}(p_{2,t_k}, x) + 1)
\]
 for all points $(x, t) \in B_{g(t_k)}(p_{2,t_k}, \delta_k^{-1} r_{2,k}) \times  [t_k - \delta_k^{-1}r_{2,k}^2, t_k]$.
\end{itemize}
Here, $K := K(n) \geq 1$ is a universal constant. 
\end{lemma}

The following three results are restatements of Lemma 3.6, Lemma 3.7, and Lemma 3.8 in \cite{BDS}. See also Lemma 5.14 in \cite{BN}.

\begin{lemma}\label{distance-derivative-estimate}
By a suitable choice of $\delta_k$, we can arrange the following holds. If $t \in [t_k - \delta_k^{-1} r_{1,k}^2, t_k]$ and $d_{g(t)}(p_{1, t_k}, x) \leq \delta_k^{-1} r_{1,k}$, then the time derivative of the distance function satisfies the estimate $0 \leq - \frac{d}{dt} d_{g(t)}(p_{1, t_k}, x) \leq 2n \,r_{1,k}^{-1}$. Similarly, if $t \in [t_k - \delta_k^{-1} r_{2,k}^2, t_k]$ and $d_{g(t)}(p_{2, t_k}, x) \leq \delta_k^{-1} r_{2,k}$, then the time derivative of the distance function satisfies the estimate $0 \leq - \frac{d}{dt} d_{g(t)}(p_{2, t_k}, x) \leq 2n\, r_{2,k}^{-1}$.
\end{lemma}

\begin{lemma}
By a suitable choice of $\delta_k$, we can arrange so that the following holds:  the two balls $B_{g(t)}(p_{1, t}, \delta_k^{-2}R(p_{1, t}, t)^{-\frac{1}{2}})$ and $B_{g(t)}(p_{2, t}, \delta_k^{-2} R(p_{2, t}, t)^{-\frac{1}{2}})$ are disjoint for $t \in (-\infty, t_k]$. 
\end{lemma}

\begin{lemma}
If $t \in (-\infty, t_k)$, then the flow is $\varepsilon_k$-symmetric at time $t$. In particular, if $(x, t) \in M \times (-\infty, t_k)$ is a spacetime point satisfying $\lambda_1(x,t) < \frac{1}{2} \theta R(x,t)$, then the point $(x, t)$ is $\varepsilon_k$-symmetric in the sense of Definition 4.2 in \cite{BN}.
\end{lemma}

By Corollary \ref{necks-outside-caps}, we can find a time $T \in (-\infty, 0]$ and a large constant $\Lambda$ with the following properties: 
\begin{itemize}
\item $L \sqrt{\frac{4n^3K}{\Lambda}} \leq 10^{-6}$.
\item If $(\bar x, \bar t) \in M \times (-\infty, T]$ satisfies $d_{g(\bar t)}(p_{1, {\bar t}}, \bar x) \geq \frac{\Lambda}{2} R(p_{1, \bar t}, \bar t)^{-\frac{1}{2}}$ and $d_{g(\bar t)}(p_{2, \bar t}, \bar x) \geq \frac{\Lambda}{2}R(p_{2, \bar t}, \bar t)^{-\frac{1}{2}}$, then $\lambda_1(x, t) < \frac{1}{2} \theta R(x,t)$ for all points $(x, t) \in B_{g(\bar t)}(\bar x, L r_{\mathrm{neck}}(\bar x, \bar t)) \times [\bar t - L r_{\mathrm{neck}}(\bar x, \bar t)^2, \bar t]$, where $r_{\mathrm{neck}}(\bar x, \bar t)^{-2} = \frac{1}{(n-1)(n-2)}R(\bar x, \bar t)$. 
\end{itemize}

The next two results are extensions of Lemma 3.9 and Lemma 3.10 in \cite{BDS} to higher dimensions. The proofs are exactly the same. 

\begin{lemma}
If $(\bar x, \bar t) \in M \times (-\infty, t_k]$ satisfies $d_{g(\bar t)}(p_{1, \bar t}, \bar x) \geq \frac{\Lambda}{2} R(p_{1, \bar t}, \bar t)^{-\frac{1}{2}}$ and $d_{g(\bar t)}(p_{2, \bar t}, \bar x) \geq \frac{\Lambda}{2} R(p_{2, \bar t}, \bar t)^{-\frac{1}{2}}$, then $(\bar x, \bar t)$ is $\frac{\varepsilon_k}{2}$-symmetric. 
\end{lemma}

\begin{lemma}
If $(\bar x, \bar t) \in M \times [t_k - \delta_k^{-1}r_{1, k}^2, t_k]$ satisfies $\Lambda r_{1,k} \leq d_{g(\bar t)}(p_{1, t_k}, \bar x) \leq \delta_k^{-1}r_{1, k}$, then $d_{g(\bar t)}(p_{1, \bar t}, \bar x) \geq \frac{\Lambda}{2} R(p_{1, \bar t}, \bar t)^{-\frac{1}{2}}$ and $d_{g(\bar t)}(p_{2, \bar t}, \bar x) \geq \frac{\Lambda}{2} R(p_{2, \bar t}, \bar t)^{-\frac{1}{2}}$. Similarly, if $(\bar x, \bar t) \in M \times [t_k - \delta_k^{-1}r_{2, k}^2, t_k]$ satisfies $\Lambda r_{2,k} \leq d_{g(\bar t)}(p_{2, t_k}, \bar x) \leq \delta_k^{-1}r_{2, k}$, then $d_{g(\bar t)}(p_{1, \bar t}, \bar x) \geq \frac{\Lambda}{2} R(p_{1, \bar t}, \bar t)^{-\frac{1}{2}}$ and $d_{g(\bar t)}(p_{2, \bar t}, \bar x) \geq \frac{\Lambda}{2} R(p_{2, \bar t}, \bar t)^{-\frac{1}{2}}$. 
\end{lemma}

The proof of the following proposition is the same as the proof of Proposition 3.11 in \cite{BDS}, except for minor difference with how we define the scale of a neck in higher dimensions. Recall that if $(\bar x, \bar t)$ lies at the center of an evolving neck, then we define $r_{\mathrm{neck}}(\bar x, \bar t)^{-2} := \frac{1}{(n-1)(n-2)} R(\bar x, \bar t)$. For the convenience of the reader, we verify this minor difference here. 

\begin{proposition}
If $(\bar x, \bar t) \in M \times [t_k - 2^{-j} \delta_k^{-1} r_{1,k}^2, t_k]$ satisfies $2^{\frac{j}{400}} \Lambda r_{1,k} \leq d_{g(\bar t)}(p_{1, t_k}, \bar x) \leq (400n^3KL)^{-j} \delta_k^{-1} r_{1,k}$, then $(\bar x, \bar t)$ is $2^{-j-1}\varepsilon_k$-symmetric. Similarly, if $(\bar x, \bar t) \in M \times [t_k - 2^{-j} \delta_k^{-1} r_{2,k}^2, t_k]$ satisfies $2^{\frac{j}{400}} \Lambda r_{2,k} \leq d_{g(\bar t)}(p_{2, t_k}, \bar x) \leq (400n^3KL)^{-j} \delta_k^{-1} r_{2,k}$, then $(\bar x, \bar t)$ is $2^{-j-1}\varepsilon_k$-symmetric. 
\end{proposition}
\begin{proof}
The proof is by induction on $j$. The assertion for $j = 0$ follows from the previous two lemmas. 

Assume $j \geq 1$ and that the assertion holds for $j -1$. Suppose $(\bar x, \bar t) \in M \times [t_k - 2^{-j} \delta_k^{-1} r_{1,k}^2, t_k]$ such that $2^{\frac{j}{400}} \Lambda r_{1,k} \leq d_{g(\bar t)}(p_{1,t_k}, \bar x) \leq (400n^3KL)^{-j} \delta_k^{-1} r_{1,k}$. Note that since $\lambda_1(\bar x, \bar t) < \frac{1}{2} \theta R(\bar x, \bar t)$,  $(\bar x, \bar t)$ lies at the center of an $\varepsilon_1$-neck. By Lemma \ref{bryant-approximation-at-tips}, 
\[
r_{\mathrm{neck}}(\bar x, \bar t)^2 = (n-1)(n-2) R(\bar x, \bar t)^{-1} \leq 4Kn^2 r_{1,k} d_{g(\bar t)}(p_{1, t_k}, \bar x). 
\]
Therefore,
\begin{align*}
\bar t - L \, r_{\mathrm{neck}}(\bar x, \bar t)^2 &\geq \bar t - 4KL n^2 r_{1,k} d_{g(\bar t)}(p_{1,t_k}, \bar x) \\
& \geq \bar t - 4KL n^2 (400n^3KL)^{-j} \delta_k^{-1} r_{1,k}^2 \\
& \geq \bar t - 2^{-j} \delta_k^{-1} r_{1,k}^2 \\
& \geq t_k - 2^{-j +1} \delta_k^{-1} r_{1,k}^2. 
\end{align*}
On the other hand, $r_{\mathrm{neck}}(\bar x, \bar t)^2 \leq 4K n^2 r_{1,k} d_{g(\bar t)}(p_{1,t_k}, \bar x) \leq \frac{4K n^2}{\Lambda} d_{g(\bar t)}(p_{1,t_k}, \bar x)^2$. Since $L \sqrt{\frac{4n^3K}{\Lambda}} \leq 10^{-6}$, we obtain 
\[
r_{\mathrm{neck}}(\bar x, \bar t) \leq 10^{-6} L^{-1} d_{g(\bar t)}(p_{1,t_k}, \bar x). 
\]
Consequently, if $x \in B_{g(\bar t)}(\bar x, L \, r_{\mathrm{neck}}(\bar x, \bar t))$, then 
\begin{align*}
d_{g(\bar t)}(p_{1,t_k}, x) &\geq d_{g(\bar t)}(p_{1,t_k}, \bar x) - L r_{\mathrm{neck}}(\bar x, \bar t) \\
& \geq (1 - 10^{-6}) d_{g(\bar t)}(p_{1,t_k}, \bar x)\\
& \geq (1 - 10^{-6}) 2^{\frac{j}{400}} \Lambda r_{1,k} \\
& \geq 2^{\frac{j-1}{400}} \Lambda r_{1,k}. 
\end{align*}
Now on the other hand, $\frac{1}{2}r_{1,k} \leq R(\bar x, \bar t)^{-\frac{1}{2}} \leq r_{\mathrm{neck}}(\bar x, \bar t)$. Putting this together with $r_{\mathrm{neck}}(\bar x, \bar t)^2 \leq 4K n^2 r_{1,k} d_{g(\bar t)}(p_{1, t_k}, \bar x)$, for all $x \in B_{g(\bar t)}(\bar x, L \, r_{\mathrm{neck}}(\bar x, \bar t))$ we obtain
\begin{align*}
d_{g(\bar t)}(p_{1,t_k}, x) + 2nL  \, r_{\mathrm{neck}}(\bar x, \bar t)^2r_{1,k}^{-1} &\leq d_{g(\bar t)}(p_{1,t_k}, \bar x) + L\, r_{\mathrm{neck}}(\bar x, \bar t) + 2n L \,r_{\mathrm{neck}}(\bar x, \bar t)^2 r_{1,k}^{-1} \\
& \leq d_{g(\bar t)}(p_{1,t_k}, \bar x) + (2n + 1) L\, r_{\mathrm{neck}}(\bar x, \bar t)^2 r_{1,k}^{-1}  \\
& \leq 400n^3K L \, d_{g(\bar t)}(p_{1,t_k}, \bar x) \\
& \leq (400n^3K L)^{-j +1} \delta_k^{-1} r_{1,k}. 
\end{align*}
Since by Lemma \ref{distance-derivative-estimate}
\[
d_{g(\bar t)}(p_{1,t_k}, x) \leq d_{g(t)}(p_{1,t_k}, x) \leq d_{g(\bar t)}(p_{1,t_k}, x) + 2n L \, r_{\mathrm{neck}}(\bar x, \bar t)^2 r_{1,k}^{-1}
\]
we conclude 
\[
2^{\frac{j-1}{400}} \Lambda r_{1,k} \leq d_{g(\bar t)}(p_{1,t_k}, x) \leq (400n^3K L)^{-j +1} \delta_k^{-1} r_{1,k}
\]
for all $(x, t) \in B_{g(\bar t)}(\bar x, L \, r_{\mathrm{neck}}(\bar x, \bar t)) \times [\bar t - L \,  r_{\mathrm{neck}}(\bar x, \bar t)^2, \bar t]$. It follows by the induction hypothesis and the Neck Improvement Theorem that the point $(\bar x, \bar t)$ is $2^{-j -1}\varepsilon_k$-symmetric. 
\end{proof}

The remaining arguments in the proof of rotational symmetry in Section 5 of \cite{BN}, go through without change to give us the following final proposition. 

\begin{proposition}
If $k$ is sufficiently large, then the flow is $\frac{\varepsilon_k}{2}$-symmetric at time $t_k$. 
\end{proposition}

The proposition above contradicts the definition of $t_k$ in view of Lemma \ref{continuity-of-symmetry}. This completes the proof of rotational symmetry.

\section{A priori estimates for compact ancient $\kappa$-solutions with rotational symmetry}
\label{a.priori.estimates}

We begin by recalling some basic facts about the Bryant soliton in higher dimensions. For the convenience of the reader we include some further discussion in Appendix \ref{properties.of.Bryant.soliton}.

\begin{proposition}[R.~Bryant \cite{Bryant}]
\label{asymptotics.of.bryant.soliton}
Consider the $n$-dimensional Bryant soliton, normalized so that the scalar curvature at the tip is equal to $1$. Then the metric can be written in the form $\Phi(r)^{-1} \, dr \otimes dr + r^2 \, g_{S^{n-1}}$, where $\Phi(r) = 1-\frac{r^2}{n(n-1)} + O(r^4)$ as $r \to 0$ and $\Phi(r) = (n-2)^2\, r^{-2} + (5-n)(n-2)^3\, r^{-4} + O(r^{-6})$ as $r \to \infty$. 
\end{proposition}

\textbf{Proof.} 
See \cite{Bryant}, Theorem 1 on p.~17.

\begin{proposition}\label{estimate.for.normalized.Bryant.profile}
Let $\eta>0$ be given. If $s$ is sufficiently small (depending on $\eta$), then 
\[
\big|\Phi((1+s)r)^{-1} - \Phi(r)^{-1}\big| \leq \eta \, \big(\Phi(r)^{-1} - 1\big)
\] 
for all $r \geq 0$.

\end{proposition}
\textbf{Proof.}
We define $\chi(r) = r^{-2} \, (\Phi(r)^{-1} - 1)$. The function $\chi(r)$ is a positive smooth function which satisfies $\chi(r) = \frac{1}{n(n-1)} + O(r^2)$ as $r \to 0$ and $\chi(r) = \frac{1}{n-2} + O(r^{-2})$ as $r \to \infty$. Hence, if $s$ is sufficiently small (depending on $\eta$), then 
\[|\chi((1+s)r) - \chi(r)| \leq \frac{\eta}{2} \, \chi(r)\] 
for all $r>0$. Therefore, if $s$ is sufficiently small (depending on $\eta$), then 
\begin{align*} 
|(1+s)^2 \, \chi((1+s)r) - \chi(r)| 
&\leq (1+s)^2 \, |\chi((1+s)r) - \chi(r)| + |(1+s)^2 -1| \, \chi(r) \\ 
&\leq \eta \, \chi(r) 
\end{align*} 
for all $r>0$. This implies 
\[\big | \Phi((1+s)r)^{-1} - \Phi(r)^{-1} \big | \leq \eta \, \big ( \Phi(r)^{-1} - 1 \big )\] 
for all $r>0$. \\

\begin{lemma}
Consider the Bryant soliton, normalized so that the scalar curvature at the tip is equal to $1$. 
Then, 
$$r\, \Phi_r + 2 \Phi = \frac{2 (n-5) (n-2)^3}{r^4} + o(r^{-4}), \qquad \mbox{for}\,\,\, r \gg 1.$$
As a result we have
\[
r\, \Phi_r + 2 \Phi  - \frac{2 (n-5)}{(n-2)} \,  \Phi^2 =  o(r^{-4}). 
\]
\end{lemma} 
\textbf{Proof.} 
According to Proposition \ref{asymptotics.of.bryant.soliton}, we have 
$$\Phi (r) = (n-2)^2\, r^{-2}  + (5-n)(n-2)^3 r^{-4} + o(r^{-4})$$
implying that 
$$r\, \Phi_r  (r) =  -2 (n-2)^2\, r^{-2}  -  4(5-n)(n-2)^3 r^{-4} + o(r^{-4})$$
for $r \gg 1$.
Hence
$$ r \Phi_r + 2\Phi =  \frac{2 (n-5) (n-2)^3}{r^4} + o(r^{-4}).
$$
This proves the first formula. The second one follows from the first and the fact that $\Phi(r)^2 = (n-2)^4 \, r^{-4} + o(r^{-6}).$

\begin{corollary} 
\label{concavity.of.F^2.on.bryant.soliton}
Consider the $n$-dimensional Bryant soliton, normalized so that the scalar curvature at the tip is equal to $1$. Let us write the metric in the form $dz \otimes dz + B(z)^2 \, g_{S^{n-1}}$. Let $\alpha =0$ if $n =4$ and $\alpha >  (n-5)/(n-2)$  if $n \geq 5$. 
Then,  there exists a large constant $L_0$ and such that   $\frac{d^2}{dz^2} \Big (  \frac {B(z)^2}2 \Big )  - \alpha \Big(\frac{d}{dz}B(z)\Big)^4<0$ holds
 if $B(z)^2 \geq \frac{L_0^2}{4}$. 
\end{corollary} 

\textbf{Proof.} 
Since $r \Phi'(r) + 2\Phi(r) = 2(n-5)(n-2)^3\, r^{-4} + O(r^{-6})$ as $r \to \infty$, we conclude that $r \Phi'(r) + 2\Phi(r) < 0$ for $r$ sufficiently large only if $n = 4$. We next observe that $\big ( \frac{d}{dz}  B(z) \big )^2 = \Phi(B(z))$. Differentiating this identity with respect to $z$ gives $2 \, \frac{d^2}{dz^2} B(z) = \Phi'(B(z))$. Thus, we conclude that $(B(z)^2)_{zz}  = 2 B(z) B_{zz} + 2 B_z^2 = B(z) \, \Phi'(B(z)) + 2 \Phi(B(z)) < 0$ if $B(z)$ is sufficiently large and $n = 4$. 
When $n \geq 5$ we need to subtract the  lower order term $2\alpha \,  \Phi^2 > 2(n-5)(n-2)^3\, r^{-4}$.   Since  $\Phi^2 = (n-2)^4 r^{-4} + \cO(r^{-6})$
we obtain a negative sign provided $\alpha > (n-5)/(n-2)$. 
This finishes the proof of the  lemma.  \\ 

\begin{corollary}
\label{asymptotics.of.bryant.soliton.2}
Consider the Bryant soliton, normalized so that the scalar curvature at the tip is equal to $1$. Let us write the metric in the form $dz \otimes dz + B(z)^2 \, g_{S^{n-1}}$. Then $B(z) \, \frac{d}{dz} B(z) \to n-2$ as $z \to \infty$.
\end{corollary}

\textbf{Proof.} 
Note that $r \, \Phi(r)^{\frac{1}{2}} \to n-2$ as $r \to \infty$. Using the identity $\big ( \frac{d}{dz} B(z) \big )^2 = \Phi(B(z))$, we obtain $B(z) \, \frac{d}{dz} B(z) = B(z) \, \Phi(B(z))^{\frac{1}{2}} \to n-2$ as $z \to \infty$. \\

We now assume that $(S^n,g(t))$ is an ancient $\kappa$-solution which is not a family of shrinking round spheres. Let $q\in S^n$ be a reference point chosen as in \cite{Angenent-Brendle-Daskalopoulos-Sesum}. The same proof as the one in \cite{Angenent-Brendle-Daskalopoulos-Sesum} implies that if $t_j \to -\infty$ and if we dilate the flow around the point $(q,t_j)$ by the factor $(-t_j)^{-\frac{1}{2}}$, then the rescaled manifolds converge to a cylinder of radius $\sqrt{2(n-2)}$.  Let $F(z,t)$ denote  the radius of a sphere of symmetry in $(S^n, g(t))$ which has signed distance $z$ from the point $q$. The function $F(z,t)$ satisfies the PDE 
\[F_t(z,t) - F_{zz}(z,t) = -\frac{n-2}{F(z,t)} \, (1-F_z(z,t)^2) - (n-1) \, F_z(z,t) \int_0^z \frac{F_{zz}(z',t)}{F(z',t)} \, dz'.\]  

Furthermore, if
\[G(\xi,\tau) = e^{\frac{\tau}{2}} F(e^{-\frac{\tau}{2}}\xi, -e^{-\tau}) - \sqrt{2(n-2)}\,\]
straightforward computation shows that 
\[
\begin{split}
G_{\tau} &= G_{\xi\xi} - \frac{\xi}{2}\, G_{\xi} + G + (n-2)\frac{G_{\xi}^2}{\sqrt{2(n-2)} + G} - \frac{G^2}{2(\sqrt{2(n-2)}+G)} \\
&- (n-1)G_{\xi}\,\int_0^{\xi} \frac{G_{\xi\xi}^2(\xi',\tau)}{\sqrt{2(n-2)} + G(\xi',\tau)}\, d\xi',
\end{split}
\]
or, equivalently, 
\[
\begin{split}
G_{\tau} &= G_{\xi\xi} - \frac{\xi}{2}\, G_{\xi} + G - \frac{G_{\xi}^2}{\sqrt{2(n-2)} + G} - \frac{G^2}{2(\sqrt{2(n-2)}+G)} \\
&+ (n-1) G_{\xi}\, \left(\frac{G_{\xi}(0, \tau)}{\sqrt{2(n-2)} + G(0,\tau)} - \int_0^{\xi} \frac{G_{\xi}^2(\xi',\tau)}{(\sqrt{2(n-2)} + G(\xi',\tau))^2}\, d\xi'\right).
\end{split}
\] 
We can write this equation as 
\[
G_{\tau} = \mathcal{L} G - \frac{G_{\xi}^2}{\sqrt{2(n-2)}} - \frac{G^2}{2 \sqrt{2(n-2)}} + E(\xi,\tau),
\]
where $E(\xi,\tau)$ is the error term and $\mathcal{L}G =  G_{\xi\xi} - \frac{\xi}{2}\, G_{\xi} + G$. As in \cite{Angenent-Brendle-Daskalopoulos-Sesum}, let $P_+$, $P_0$ and $P_-$ be orthogonal projections associated with the direct sum $\mathcal{H} = \mathcal{H}_+ \oplus \mathcal{H}_0 \oplus \mathcal{H}_-$, where $\mathcal{H}_+$, $\mathcal{H}_0$ and $\mathcal{H}_-$ are the positive, zero and negative eigenspaces with respect to operator $\mathcal{L}$, respectively. Exactly the same reasoning and arguments as in \cite{Angenent-Brendle-Daskalopoulos-Sesum} yield that for $\tau$ sufficiently small the positive mode, i.e. the projection onto $\mathcal{H}_0$ dominates and that
\[
\int_{|\xi| < \delta(\tau)^{-\frac{1}{100}}} e^{-\frac{\xi^2}{4}} E(\xi,\tau)\, (\xi^2-2)\, d\xi = O(A(\tau)^2),
\]
where $A(\tau)$ is the norm of the orthogonal projection of $\hat{G}(\xi,\tau) = G(\xi,\tau) \chi(\delta(\tau)^{\frac{1}{100}}\xi)$ onto $\mathcal{H}_0$. Note that $\chi$ is a cut off function with the support in a parabolic region and $\lim_{\tau\to -\infty} \delta(\tau) = 0$, where both, $\chi$ and $\delta(\tau)$ are defined in the same way as in \cite{Angenent-Brendle-Daskalopoulos-Sesum}. Having the equation for $G$ and the integral estimate above, the same arguments as in \cite{Angenent-Brendle-Daskalopoulos-Sesum} imply the following asymptotics:
\begin{theorem}
\label{asymptotics.compact.case}
Let $(S^n,g(t))$ be a rotationally symmetric ancient $\kappa$-solution which is not isometric to a family of shrinking spheres. Then we can find a reference point $q \in S^n$ such that the following holds. Let $F(z,t)$ denote the radius of the sphere of symmetry in $(S^n,g(t))$ which has signed distance $z$ from the reference point $q$. Then the profile $F(z,t)$ has the following asymptotic expansions:
\begin{itemize}
\item[(i)] Fix a large number $L$. Then, as $t \to -\infty$, we have 
\[\frac 1{2} \, F(z,t)^2   = (n-2) \Big [ (-t)  - \frac{z^2+2t}{4\log(-t)}  \Big ]+ o \Big ( \frac{(-t)}{\log(-t)} \Big ) \] 
for $|z| \leq L\sqrt{-t}$
\item[(ii)] Fix a small number $\theta>0$. Then as $t \to -\infty$, we have 
\[\frac 1{2} \, F(z,t)^2 = (n-2) \big [ (-t)  -  \frac{z^2 +2t }{4\log(-t)} \big ] + o(-t)\] 
for $|z| \leq 2\sqrt{(1 - \theta^2}) \sqrt{(-t) \log(-t)}$.
\item[(iii)] The reference point $q$ has distance $(2 +o(1)) \sqrt{(-t) \log(-t)}$ from each tip. The scalar curvature at each tip is given by $(1+o(1)) \, \frac{\log(-t)}{(-t)}$. Finally, if we rescale the solution around one of the tips, then the rescaled solutions converge to the Bryant soliton as $t \to -\infty$.
\end{itemize}
\end{theorem}

We next  let $H(z,t) :=  \frac 12  F(z,t)^2 + (n-2) \, t$, $K(z, t) := F_z(z,t)^4$ and consider the quantity 
\[Q(z,t) := H_{zz} - \alpha(n) K.\]
where $\alpha(n) = 0$ if $n =4$ and $\alpha(n) = 1 > (n-4)/(n-2)$ if $n \geq 5$. This definition of $\alpha(n)$ ensures both Lemma \ref{boundary.points} and Lemma \ref{Lemma:Q} will hold. 

\begin{lemma}
\label{boundary.points}
Let $L_0$ be the constant in Corollary \ref{concavity.of.F^2.on.bryant.soliton}.  There exists a time $T_0 < 0$ with the following property. If $t \leq T_0$, $F(z,t)^2 = L_0^2 \, \frac{(-t)}{\log (-t)}$, then $Q(z,t) < 0$.
\end{lemma}

\textbf{Proof.} 
By our assumptions on $\alpha(n)$, the proof of this is identical to the proof of Lemma 4.5 in \cite{BDS} once we use Corollary \ref{concavity.of.F^2.on.bryant.soliton}.

\begin{lemma} The function $H(z,t)$ satisfies the equation
 \[H_t(z,t) - H_{zz}(z,t) = (n-3) \frac{H_z(z,t)^2 }{F^2(z,t)}  - (n-1) \, H_z  \int_0^z \frac{F_{zz}(z',t)}{F(z',t)} \, dz'.\]

\end{lemma} 

\begin{proof} We have $H_t= F F_t$, $H_z=F F_z$, $H_{zz} = F F_{zz} + F_z^2$. Hence,
$$H_t - H_{zz} = - F_z^2 - (n-2) (1-F_z^2) + (n-2) - (n-1)  \, F F_z\int_0^z \frac{F_{zz}(z',t)}{F(z',t)} \, dz' $$
which gives
$$H_t - H_{zz} = (n-3) \frac{H_z^2 }{F^2}  - (n-1) \, H_z  \int_0^z \frac{F_{zz}(z',t)}{F(z',t)} \, dz'. $$
\end{proof}

\begin{lemma}
\label{evolution.of.H_zz}
The function $H_{zz}(z,t)$ satisfies the evolution equation 
\[
\begin{split}
&H_{zzt} (z,t) - H_{zzzz}(z,t) \\ 
&= \left( (n-5) \frac{F_z(z,t)}{F(z,t)} - (n-1)\, \int_0^z \frac{F_{zz}(z',t)}{F(z',t)} \, dz' \right)H_{zzz} (z,t)  \\
&  - 4(n-4)\, \frac{F_z(z,t)^2 F_{zz}(z,t)}{F(z,t)} - 4 F_{zz}(z,t)^2.
\end{split}
\]
The function $K(z,t)$ satisfies the evolution equation
\[
\begin{split}
&K_t(z,t) - K_{zz}(z,t) \\
&=\left( (n-5) \frac{F_z(z,t)}{F(z,t)} - (n-1)   \int_0^z \frac{F_{zz}(z',t)}{F(z',t)} \, dz'  \right) K_z(z,t) \\
& + 8 \frac {F_z(z,t)^4F_{zz}(z,t)}{F(z,t)}  - 12 F_z(z,t)^2 F_{zz}(z,t)^2+ 4 (n-2) \frac{1-F_z(z,t)^2}{F(z,t)^2} F_z(z,t)^4.
\end{split}
\]
\end{lemma}

\textbf{Proof.} 
We differentiate twice the equation of $H$ to find 
\bee
\begin{split}
H_{zzt} - H_{zzzz} &= (n-3) \big ( \frac{H_z^2 }{F^2} \big )_{zz}   - (n-1) H_{zzz}   \int_0^z \frac{F_{zz}(z',t)}{F(z',t)} \, dz' \\
&\quad - 2 (n-1)  \frac {H_{zz} F_{zz} }{F} - (n-1) \, H_z \big ( \frac{F_{zz}}{F} \big )_z.
\end{split} 
\eee
Next we use $H_z = F F_z$, $H_{zz} = F  F_{zz} + F_z^2$ and $H_{zzz} = F F_{zzz} + 3 F_z F_{zz}$ to compute 
$$H_z  \big ( \frac{F_{zz}}{F} \big )_z = F F_z  \frac{FF_{zzz} - F_{zz}F_z}{F^2} =  H_{zzz}  \frac{F_z}F - 4 \frac{F_z^2 F_{zz}}{F}$$
and also use 
$$\big ( \frac{H_z^2 }{F^2} \big )_{zz} = (F_z^2)_{zz} = 2 (F_z F_{zz})_z = 2 F_z F_{zzz} + 2 F_{zz}^2 
= 2 \frac{F_z}{F} \, H_{zzz} - 6 \frac{F_z^2 F_{zz}}{F} + 2 F_{zz}^2.$$
Combining the above yields 
\bee
\begin{split}
H_{zzt} - H_{zzzz} &= H_{zzz} \Big ( - (n-1)   \int_0^z \frac{F_{zz}(z',t)}{F(z',t)} \, dz'  + (n-5) \frac{F_z}F \Big ) \\
& \quad  - 6(n-3)  \frac{F_z^2 F_{zz}}{F} + 2(n-3)  F_{zz}^2 - 2 (n-1) F_{zz}^2 \\&\quad - 2(n-1)  \frac{F_z^2 F_{zz} }{F} + 4 (n-1)  \frac{F_z^2 F_{zz}}{F}\\
&= H_{zzz} \Big ( - (n-1)    \int_0^z \frac{F_{zz}(z',t)}{F(z',t)} \, dz'  + (n-5) \frac{F_z}F \Big ) \\
& \quad  - 4(n-4)  \frac{F_z^2 F_{zz}}{F} - 4  F_{zz}^2.\\
\end{split}
\eee
Next, a simple  computation shows that $F_z$ satisfies the equation 
$$F_{zt} = F_{zzz} + \Big ( (n-3) \frac{F_z}{F} - (n-1)   \int_0^z \frac{F_{zz}(z',t)}{F(z',t)} \, dz'  \Big ) F_{zz}  + (n-2) (1-F_z^2) \, \frac{F_z}{F^2}.$$
Set $K=F_z^4$. Using  the previous equation and $K_t= 4F_z^3 F_{zt}$, $K_z=4F_z^3  F_{zz}$, $K_{zz} = 4F_z^3 F_{zzz} + 12 F_z^2 F_{zz}^2$  we find that
$K$ satisfies
\[
\begin{split}
K_t - K_{zz} &= - 12F_z^2 F_{zz}^2  + 4F_z^3\Big ( (n-3) \frac{F_z}{F} - (n-1)    \int_0^z \frac{F_{zz}(z',t)}{F(z',t)} \, dz'  \Big )F_{zz}\\
&\quad+  4 (n-2)(1-F_z^2) \frac{F_z^4}{F^2} \\
&=\Big ( (n-5) \frac{F_z}{F} - (n-1)   \int_0^z \frac{F_{zz}(z',t)}{F(z',t)} \, dz'  \Big ) K_z \\&\quad + 8 F_z^4\frac{F_{zz} }{F}  - 12F_z^2F_{zz}^2 + 4 (n-2)(1-F_z^2) \frac{F_z^4}{F^2} \end{split}
\]
This completes the proof of the lemma. 
\medskip

Combining evolution equations for $H_{zz}$ and $K$ we obtain that $Q= H_{zz} -\alpha K$ satisfies 
\[
\begin{split}
Q_t - Q_{zz} &= \Big ( (n-5) \frac{F_z}{F} - (n-1) \, \int_0^z \frac{F_{zz}(z',t)}{F(z',t)^2} \, dz' \Big ) Q_z - 4(n-4)  \frac{F_z^2 F_{zz}}{F}  \\
&\quad - 4  F_{zz}^2  - 8\alpha  \frac {F_z^4 F_{zz}}{F}  + 12 \alpha \, F_z^2 F_{zz}^2  - 4 \alpha  (n-2) (1-F_z^2) \frac{F_z^4}{F^2}\\
& = \Big ( (n-5) \frac{F_z}{F} - (n-1) \, \int_0^z \frac{F_{zz}(z',t)}{F(z',t)^2} \, dz' \Big ) Q_z  \\
&\quad   -(4(n-4) + 8\alpha F_z^2)\frac {F_z^2 F_{zz}}{F}  - (4-12 \alpha  F_z^2 )F_{zz}^2  - 4 \alpha  (n-2) \frac{ 1-F_z^2}{F^2}F_z^4.
\end{split}
\]

We will use the last equation, Lemma \ref{boundary.points} and the maximum principle to prove the following auxiliary lemma, which will be used to establish Proposition \ref{concavity.of.F2} below.
\medskip 
\begin{lemma}\label{Lemma:Q} Let $L_0$ be chosen as in Corollary \ref{concavity.of.F^2.on.bryant.soliton}. Then, there exists $T_0 \ll -1$ (depending on $L_0$ and dimension $n$) such that   if  for some $t_0 \leq T_0$, we have $Q(z_0,t_0):= \max_{ F(z_0,t_0)^2  \geq \frac{L_0^2 (-t_0)}{\log (-t_0)} } Q(z,t_0)>0$, then  
$$\sup_{ F(z,t)^2  \geq \frac{L_0^2 (-t)}{\log (-t)} } Q(z,t) \geq Q(z_0,t_0)$$
for all $t \leq t_0$.
 \end{lemma} 

\textbf{Proof.} When $n=4$, $\alpha(n)=0$, and $Q(z,t) = H_{zz}(z,t)$. The maximum principle applied to the evolution equation of $H$ yields the claimed inequality, as in the proof of Proposition 4.7 in \cite{BDS}. Hence, we will assume that $n \geq 5$. Let 
\[
q(t) = \sup_{F(z,t)^2 \geq \frac{L_0^2(-t)}{\log(-t)}} Q(z,t).
\]
Let $T_0$ be determined by Lemma \ref{boundary.points} and assume $t_0 \leq T_0$. Assume $Q(z_0, t_0) := q(t_0) > 0$. We must show $q(t) \geq q(t_0)$ for all $t \leq t_0$.

Suppose not. Let $t_1 =\sup\{t \leq t_0 : q(t) < q(t_0)\}$. Then $q(t_1) = q(t_0)$ and there exists a sequence of times $t_{1,j} \nearrow t_1$ such that $q(t_{1,j}) < q(t_0)$. Choose $z_1$ such that $F(z_1, t_1)^2\geq \frac{L_0^2(-t)}{\log(-t_1)}$ and $q(t_1) = Q(z_1, t_1) = Q(z_0, t_0) > 0$. Since $t_1 \leq T_0$, Lemma \ref{boundary.points} implies $Q \leq 0$ when $ F(z,t_1)^2  = L_0^2 \frac{(-t_1)}{\log (-t_1)}$. Hence $z_1$ is an interior spacial maximum. Applying the maximum principle to the evolution equation for $Q$ shows that at such interior maximum we have 
\[
Q_t(z_1, t_1) \leq - \big ( 4(n-4) +8\alpha F_z^2 \big )   \frac{F_z^2 F_{zz}}{F} - (4-12 \alpha F_z^2) \,   F_{zz}^2    - 4 \alpha  (n-2) (1-F_z^2) \frac{F_z^4}{F^2},
\]
where the right hand side is evaluated at $(z_1, t_1)$. 
Using that $|F_z|(z_1, t_1) \ll 1$ (every point in the considered region lies on a very fine neck assuming $-T_0$ is sufficiently large), we have $4-12 \alpha F_z^2 >0$, and hence we conclude the inequality 
\[
Q_t(z_1,t_1) \leq - \big ( 4(n-4) + 8 \alpha F_z^2  \big )   \frac{F_z^2 F_{zz}}{F}     - 4 \alpha  (n-2) (1-F_z^2) \, \frac{F_z^4}{F^2}.
\]
Next express $\displaystyle F_{zz} = \frac{Q - F_z^2 +\alpha F_z^4}{F}$ which gives 
$$  \frac{F_z^2 F_{zz}}{F}  = \frac{F_z^2 Q - F_z^4 +\alpha F_z^6}{F^2}$$   so that 
$$Q_t(z_1,t_1) \leq - \big ( 4(n-4) + 8\alpha F_z^2  \big )   \frac{F_z^2 }{F^2} \, Q + 
4 \big ( (n-4) -  \alpha (n-2) \big ) \frac{F_z^4 }{F^2}   + O(\frac{F_z^6}{F^2}).$$
Using again that $|F_z|(z_1, t_1) \ll 1$, we conclude that if $Q(z_1,t_1) >0$, then $Q_t(z_1,t_1) <0$. But this implies $Q(t_{1,j}, z_1) > Q(z_1, t_1)$ for $t_{1,j}$ sufficiently close to $t_1$ and hence $q(t_{1,j}) > q(t_1) = q(t_0)$, a contradiction. Therefore, we must have $q(t) \geq q(t_0)$ for all $t \leq t_0$, completing the proof.

\medskip

As  consequence we can show the following analogue of Proposition 4.7 in \cite{BDS}. 

\begin{proposition} 
\label{concavity.of.F2}
Let $L_0$ be chosen as in Corollary \ref{concavity.of.F^2.on.bryant.soliton},  $T_0$ be chosen as in Lemma \ref{boundary.points}.
If $t \leq T_0$ and $F(z,t)^2 \geq L_0^2 \, \frac{(-t)}{\log (-t)}$, then $Q(z,t)\leq 0$.
\end{proposition}

\textbf{Proof.} 
Suppose this is false. Then we can find a point $(z_0,t_0)$ such that $t_0 \leq T_0$, $F(z_0,t_0)^2 \geq L_0^2 \, \frac{(-t_0)}{\log (-t_0)}$, and $Q(z_0,t_0) > 0$. In view of Lemma \ref{boundary.points} and  Lemma \ref{Lemma:Q}, we have 
\[\sup_{F(z,t)^2 \geq L_0^2 \, \frac{(-t)}{\log (-t)}} Q(z,t) \geq Q(z_0,t_0) > 0\] 
for each $t \leq t_0$. Let us consider a sequence $t_j \to -\infty$. For $j$ large, we can find a point $z_j$ such that $F(z_j,t_j)^2 \geq L_0^2 \, \frac{(-t_j)}{\log (-t_j)}$ and $Q(z_j,t_j) \geq Q(z_0,t_0) > 0$.   Using the inequality $F_{zz} \leq 0$, we obtain $F_z(z_j,t_j)^2 \geq Q(z_j,t_j) \geq Q(z_0,t_0) > 0$ for $j$ large. Hence, if we rescale around the points $(z_j,t_j)$ and pass to the limit, then the limit cannot be a cylinder. Consequently, the limit of these rescalings must be the Bryant soliton. Hence, after passing to the limit, we obtain a point $z_\infty$ on the Bryant soliton such that $B(z_\infty)^2 \geq L_0^2$ and $$  \frac {d^2}{dz^2 } \Big ( \frac{B(z)^2}2\Big )   - \alpha(n) \Big ( \frac {d}{dz} B(z) \Big )^4 \geq 0, \qquad\mbox{at}\,\, z_\infty.$$   This contradicts Corollary \ref{concavity.of.F^2.on.bryant.soliton}. \\

\medskip

We next recall a crucial estimate from \cite{Angenent-Brendle-Daskalopoulos-Sesum}. \\

\begin{proposition}[cf. \cite{Angenent-Brendle-Daskalopoulos-Sesum}]
\label{precise.estimate.for.F}
Fix a small number $\theta > 0$ and a small number $\eta > 0$. Then 
\[\Big | \frac{1}{2} \, F(z,t)^2 + (n-2)\, t + (n-2)\,\frac{z^2+2t}{4 \log(-t)} \Big | \leq \eta \, \frac{z^2-t}{\log(-t)}\] 
if $F(z,t) \geq \frac{\theta}{400} \sqrt{-t}$ and $-t$ is sufficiently large (depending on $\eta$ and $\theta$). 
\end{proposition}

\textbf{Proof.} 
The proof is analogous to the proof of Proposition 4.8 in \cite{BDS} and is based on the asymptotics of the ancient solution in radially symmetric setting. In \cite{Angenent-Brendle-Daskalopoulos-Sesum} we showed precise asymptotics in the case $n = 3$, but the proof carries over without any changes to higher dimensions due to new results in \cite{BN} and \cite{LZ}. \\

\begin{proposition}
\label{precise.estimate.for.F_z}
Let us fix a small number $\theta > 0$ and a small number $\eta > 0$. Then 
\[\Big | F(z,t) \, F_z(z,t) + \frac{(n-2)\, z}{2 \log(-t)} \Big | \leq \eta \, \frac{|z|+\sqrt{-t}}{\log(-t)}\] 
if $F(z,t) \geq \frac{\theta}{200} \sqrt{-t}$ and $-t$ is sufficiently large (depending on $\eta$ and $\theta$).
\end{proposition}

\textbf{Proof.}
If $n = 4$,  by Proposition \ref{concavity.of.F2} we have that $F^2$ is still concave on the set $F(z,t)^2 \ge L_0^2\, \frac{(-t)}{\log(-t)}$ for appropriately chosen $L_0 \gg 1$ and for all $-t$ sufficiently large, and  the proof is the same to the proof of Proposition 4.9 in \cite{BDS}. 

Now assume $n\geq5$. The beginning of the proof is similar to the proof of Proposition 4.9 in \cite{BDS}. The difference comes  from the fact that $FF_z$ is not monotone decreasing, but instead one needs to use Proposition \ref{concavity.of.F2}. 

Let $\theta \in (0,\frac{1}{2})$ and $\eta \in (0,\frac{1}{2})$ be given. We can find a small positive number $\mu \in (0,\frac{\eta}{n-2})$ and time $T_0$ with the property that $F((1+\mu)z,t) \geq \frac{\theta}{400} \sqrt{-t}$ whenever $F(z,t) \geq \frac{\theta}{200} \sqrt{-t}$ and $t \leq T_0$. Moreover, by Proposition \ref{precise.estimate.for.F}, we can find a time $T \leq T_0$ such that 
\[\Big | \frac{1}{2} \, F(z,t)^2 + (n-2)\, t + (n-2)\,\frac{z^2+2t}{4 \log(-t)} \Big | \leq \eta \, \mu \, \frac{z^2}{16 \log(-t)}\] 
whenever $z \geq 4\sqrt{-t_0}$, $F(z,t) \geq \frac{\theta}{400} \sqrt{-t}$, and $t \leq T$.

Suppose now that $(z_0,t_0)$ is a point in spacetime satisfying $z_0 \geq 4\sqrt{-t_0}$, $F(z_0,t_0) \geq \frac{\theta}{200} \sqrt{-t_0}$, and $t_0 \leq T$.  We assume $-T$ is sufficiently large so that $F_z(z, t_0) < 0$ for all $z \geq 4\sqrt{-t_0}$. By the above, we have $F(z,t_0) \geq \frac{\theta}{400} \sqrt{-t_0}$ for all $z \in [(1-\mu)z_0,(1+\mu)z_0]$.  Consequently, 
\[\Big | \frac{1}{2} \, F(z,t_0)^2 + (n-2)\, t_0 + (n-2)\,\frac{z^2+2t_0}{4 \log(-t_0)} \Big | \leq \eta \, \mu \, \frac{z_0^2}{4 \log(-t_0)}\] 
for all $z \in [(1-\mu)z_0,(1+\mu)z_0]$. This implies 
\begin{equation*}
\inf_{z \in [(1-\mu)z_0,z_0]} \Big ( F(z,t_0) \, F_z(z,t_0) + \frac{(n-2)\,z}{2 \log(-t_0)} \Big ) \leq \eta \, \frac{z_0}{2 \log(-t_0)}
\end{equation*}
and 
\begin{equation*}
\sup_{z \in [z_0,(1+\mu)z_0]} \Big ( F(z,t_0) \, F_z(z,t_0) + \frac{(n-2)\,z}{2 \log(-t_0)} \Big ) \geq -\eta \, \frac{z_0}{2 \log(-t_0)}.
\end{equation*}

Define a function 
\[
S(z,t) := \frac{1}{F(z, t)^2F_z(z,t)^2} - \frac{\alpha}{F(z,t)^2} = \frac{1 - \alpha F_z(z, t)^2}{F(z,t)^2 F_{z}(z,t)^2}. 
\]
Since $F(z, t_0) \geq \frac{\theta}{400} \sqrt{-t_0}$ for all $z \in [(1-\mu)z_0, (1+\mu)z_0]$, we have  $F_z(z,t_0)^2 \ll 1$. In particular, $S(z, t_0) > 0$ if $z \in [(1-\mu)z_0, (1+\mu)z_0]$. Moreover, we may assume $-t_0$ is sufficiently large, depending on $\theta$ and $\eta$, so that 
\[
F(z,t_0)|F_z(z,t_0)| \leq S(z,t_0)^{-\frac{1}{2}} \leq (1+\eta)F(z,t_0)|F_z(z,t_0)|.
\]
Since $F_z < 0$, this means 
\[
 (1+\eta)F(z,t_0)F_z(z,t_0) \leq -S(z,t_0)^{-\frac{1}{2}} \leq F(z,t_0)F_z(z,t_0).
\]
We compute
\begin{align*}
S_z &= -\frac{2}{F^3F_z} -\frac{2F_{zz}}{F^2F_z^3}+ \frac{2\alpha F_z}{F^3}\\
& = -2 \frac{1}{F^3F_z^3}\big(F_z^2 + F F_{zz}  - \alpha F_z^4\big) \\
& = 2 \frac{Q}{F^3(-F_z)^3}.
\end{align*}
Since $F_z < 0$ and by Proposition \ref{concavity.of.F2} $Q\leq 0$, the function $S$ is monotone decreasing in the variable $z$ for $z \in [(1-\mu)z_0, (1+\mu)z_0]$. This means $-S^{-\frac{1}{2}}$ is monotone decreasing as well. So although the function $FF_z$ is not monotone decreasing in higher dimensions, it is very close to the monotone function $-S^{-\frac{1}{2}}$. 

Now it follows from our estimates above that 
\begin{equation*}
\inf_{z \in [(1-\mu)z_0,z_0]} \Big ( -S(z,t_0)^{-\frac{1}{2}}+ \frac{(n-2)\,z}{2 \log(-t_0)} \Big ) \leq \eta \, \frac{z_0}{2 \log(-t_0)}
\end{equation*}
and
\begin{equation*}
\sup_{z \in [z_0,(1+\mu)z_0]} \Big ( -\frac{1}{1+\eta}S(z,t_0)^{-\frac{1}{2}} + \frac{(n-2)\,z}{2 \log(-t_0)} \Big ) \geq -\eta \, \frac{z_0}{2 \log(-t_0)}.
\end{equation*}
Since $-S^{-\frac{1}{2}}$ is monotone decreasing in the relevant region, this implies
\begin{equation*}
-S(z_0,t_0)^{-\frac{1}{2}}+ \frac{(n-2)(1-\mu)z_0}{2 \log(-t_0)} \leq \eta \, \frac{z_0}{2 \log(-t_0)}
\end{equation*}
and
\begin{equation*}
-\frac{1}{1+\eta}S(z_0,t_0)^{-\frac{1}{2}} + \frac{(n-2)(1-\mu)z_0}{2 \log(-t_0)}  \geq -\eta \, \frac{z_0}{2 \log(-t_0)}.
\end{equation*}
Hence
\begin{equation*}
(1+\eta)F(z_0, t_0)F_z(z_0, t_0)+ \frac{(n-2)z_0}{2 \log(-t_0)} \leq ((n-2)\mu +\eta) \, \frac{z_0}{2 \log(-t_0)}
\end{equation*}
and
\begin{equation*}
\frac{1}{1+\eta}F(z_0, t_0) F_z(z_0, t_0)+ \frac{(n-2)z_0}{2 \log(-t_0)}  \geq -((n-2)\mu +\eta) \, \frac{z_0}{2 \log(-t_0)}.
\end{equation*}
Since $\mu \in(0, \frac{\eta}{n-2})$ and evidently $\frac{1}{1 + \eta}  < 1 < 1 + \eta$, it follows that
\begin{equation*}
\Big|F(z_0, t_0)F_z(z_0, t_0)+ \frac{(n-2)z_0}{2 \log(-t_0)}\Big| \leq \eta \frac{z_0}{\log(-t_0)}.
\end{equation*}
To summarize, we have verified the assertion for $z \geq 4\sqrt{-t}$. An analogous argument show the assertion holds for $z \leq-4\sqrt{-t}$. Finally, suppose $|z| \leq 4\sqrt{-t}$. In  dimension three, this case follows from Proposition 5.10 in \cite{Angenent-Brendle-Daskalopoulos-Sesum}. An analogous result holds in higher dimensions; namely $(-\tau) G(\xi, \tau) \to - \frac{1}{4\sqrt{2(n-2)}}(\xi^2 - 2)$ in $C^{\infty}_{\mathrm{loc}}$. The assertion in the region $|z| \leq 4\sqrt{-t}$ (which is equivalent to $|\xi| \leq 4$) follows directly from this result. This completes the proof of Proposition \ref{precise.estimate.for.F_z}.

\begin{corollary}
\label{bound.for.F_z} 
Let us fix a small number $\theta > 0$. Then 
\[|F_z(z,t)| \leq \frac{C(\theta)}{\sqrt{\log(-t)}}\] 
if $F(z,t) \geq \frac{\theta}{200} \sqrt{-t}$ and $-t$ is sufficiently large (depending on $\theta$).
\end{corollary}

\textbf{Proof.} 
Given Theorem \ref{asymptotics.compact.case}, the same reasoning as in \cite{Angenent-Brendle-Daskalopoulos-Sesum} imply that $|z| \leq (2+o(1)) \sqrt{(-t) \log(-t)}$. Hence, the assertion follows from Proposition \ref{precise.estimate.for.F_z}. \\

\begin{proposition} 
\label{higher.derivative.bounds.in.cylindrical.region}
Let us fix a small number $\theta > 0$. Then 
\[F(z,t) \, |F_{zz}(z,t)| + F(z,t)^2 \, |F_{zzz}(z,t)| \leq \frac{C(\theta)}{\sqrt{\log (-t)}}\] 
if $F(z,t) \geq \frac{\theta}{100} \sqrt{-t}$ and $-t$ is sufficiently large (depending on $\theta$).
\end{proposition}

\textbf{Proof.} 
The proof is analogous to the proof of Proposition 4.11 in \cite{BDS}.
\\

\begin{proposition}
\label{derivative.of.F.wrt.t}
Let us fix a small number $\theta > 0$. Then 
\[|n-2 + FF_t| \leq \frac{C(\theta)}{\sqrt{\log(-t)}}\] 
whenever $F \geq \frac{\theta}{100} \sqrt{-t}$, and $-t$ is sufficiently large (depending on $\theta$).
\end{proposition}

\textbf{Proof.} 
Using the evolution equation for $F$, we obtain 
\begin{align*} 
&(n-2) + F(z,t) \, F_t(z,t) \\ 
&= F(z,t) \, F_{zz}(z,t) - F_z(z,t)^2 \\ 
&+ (n-1)  \, F(z,t) \, F_z(z,t) \, \bigg [ F(0,t)^{-1} \, F_z(0,t) - \int_0^z \frac{F_z(z',t)^2}{F(z',t)^2} \, dz' \bigg ]. 
\end{align*} 
The same arguments  in the proof of Proposition 4.12 in \cite{BDS} yield the proof in higher dimensional case as well. \\

\begin{proposition}
\label{neck}
Let $\varepsilon>0$ be given. Then there exists a large number $L$ (depending on $\varepsilon$) and a time $T$ such that the following holds. If $F \geq L \, \sqrt{\frac{(-t)}{\log(-t)}}$ and $t \leq T$ at some point in space-time, then that point lies at the center of an evolving $\varepsilon$-neck. 
\end{proposition}

\textbf{Proof.} This follows from the fact, that follows in an analogous way as in \cite{Angenent-Brendle-Daskalopoulos-Sesum}, that the scalar curvature at each tip is comparable $\frac{\log(-t)}{(-t)}$. \\

\begin{corollary}
\label{higher.derivative.bounds.in.collar.region}
Let $\eta>0$ be given. Then there exists a large number $L$ (depending on $\eta$) and a time $T$ such that $|F_z| + F \, |F_{zz}| + F^2 \, |F_{zzz}| \leq \eta$ whenever $F \geq L \, \sqrt{\frac{(-t)}{\log(-t)}}$ and $t \leq T$.
\end{corollary}

\textbf{Proof.} 
This follows directly from Proposition \ref{neck}. \\

\begin{proposition}
\label{estimate.for.F_z.in.collar.region}
Let $\eta>0$ be given. Then there exist a large number $L \in (\eta^{-1},\infty)$ and a small number $\theta \in (0,\eta)$ (depending on $\eta$), and a time $T$ with the property that 
\[\bigg | (n-2) - \sqrt{\frac{\log(-t)}{(-t)}} \, F \, |F_z| \bigg | \leq \eta\] 
whenever $L \, \sqrt{\frac{(-t)}{\log(-t)}} \leq F \leq 100\theta \sqrt{-2t}$ and $t \leq T$.
\end{proposition}

\textbf{Proof.} 
By Corollary \ref{asymptotics.of.bryant.soliton.2}, we can find a large number $L \in (\eta^{-1},\infty)$ such that $\big | (n-2) - B(z) \, \frac{d}{dz} B(z) \big | \leq \frac{\eta}{2}$ for $z \geq \frac{L}{2}$. Recall that the solution looks like the Bryant soliton near each tip, and the scalar curvature at each tip equals $(1+o(1)) \, \frac{\log(-t)}{(-t)}$. Consequently, 
\[
\bigg | (n-2) - \sqrt{\frac{\log(-t)}{(-t)}} \, F \, |F_z| \bigg | \leq \frac{\eta}{8},
\]
if $F = L \, \sqrt{\frac{(-t)}{\log(-t)}}$ and $-t$ is sufficiently large. Now recall the function $S(z, t)$ defined in Proposition \ref{precise.estimate.for.F_z}:
\[S(z,t) := \frac{1}{F(z,t)^2F_z(z,t)^2} - \frac{\alpha}{F(z,t)^2}.\]
Let $\mu = \min\{\frac{\eta}{4(n-2)}, 1\}$. The function $z \mapsto S(z, t)^{-\frac{1}{2}}$ is monotone decreasing. Moreover, whenever $F \geq L \sqrt{\frac{(-t)}{\log(-t)}}$, $L$ is large, and $-t$ is sufficiently large, $|F_z|$ is small and hence 
\[F|F_z| \leq S^{-\frac{1}{2}} \leq (1 + \mu)F|F_z|.\] 
Consequently,
\begin{align*}
(n-2) - \sqrt{\frac{\log(-t)}{(-t)}} S^{-\frac{1}{2}} &\geq (n-2) - (1+ \mu)\sqrt{\frac{\log(-t)}{(-t)}} \,F \,|F_z|  \\
& \geq -(1+ \mu)\bigg | (n-2) - \sqrt{\frac{\log(-t)}{(-t)}} \, F \, |F_z| \bigg | - \frac{\eta}{4}
\end{align*}
and 
\begin{align*}
(n-2) - \sqrt{\frac{\log(-t)}{(-t)}} S^{-\frac{1}{2}} &\leq (n-2) - \sqrt{\frac{\log(-t)}{(-t)}}\, F\, |F_z|  \\
& \leq (1+\mu)\bigg | (n-2) - \sqrt{\frac{\log(-t)}{(-t)}} \, F \, |F_z| \bigg | + \frac{\eta}{4}
\end{align*}
whenever $F \geq  L \, \sqrt{\frac{(-t)}{\log(-t)}}$ and $-t$ is sufficiently large. In other words, 
\[
\bigg | (n-2) - \sqrt{\frac{\log(-t)}{(-t)}} \,S^{-\frac{1}{2}}\bigg | \leq 2\bigg | (n-2) - \sqrt{\frac{\log(-t)}{(-t)}} \, F \, |F_z| \bigg | + \frac{\eta}{4}
\]
whenever $F \geq  L \, \sqrt{\frac{(-t)}{\log(-t)}}$ and $-t$ is sufficiently large. Similarly, we can show 
\[
\bigg | (n-2) - \sqrt{\frac{\log(-t)}{(-t)}} \,F|F_z|\bigg | \leq \bigg | (n-2) - \sqrt{\frac{\log(-t)}{(-t)}} \,S^{-\frac{1}{2}} \bigg | + \frac{\eta}{2}.
\]
 In particular, from the first estimate in the proof, we obtain
\[
\bigg | (n-2) - \sqrt{\frac{\log(-t)}{(-t)}} \,S^{-\frac{1}{2}}\bigg | \leq \frac{\eta}{2}
\]
if $F = L \, \sqrt{\frac{(-t)}{\log(-t)}}$ and $-t$ is sufficiently large.
 On the other hand, for each $\theta \in (0,\frac{1}{1000})$, Proposition \ref{precise.estimate.for.F} implies
\[z^2 = (4+o(1)) \, \Big(1-\frac{(100\theta)^2}{n-2}\Big) \, (-t) \log(-t)\] 
if $F = 100\theta \sqrt{-2t}$. Using Proposition \ref{precise.estimate.for.F_z}, we obtain 
\[F \, |F_z| = ((n-2)+o(1)) \, \frac{|z|}{2 \log(-t)}\] 
if $F = 100\theta \sqrt{-2t}$. Consequently, 
\[(n-2) - \sqrt{\frac{\log(-t)}{(-t)}} \, F \, |F_z| = (n-2)\, \Big(1 - \sqrt{1-\frac{(100\theta)^2}{n-2}}+o(1)\Big)\]
if $F = 100\theta \sqrt{-2t}$. Therefore, if we choose $\theta$ sufficiently small (depending on $\eta$), then we obtain 
\[
\bigg | (n-2) - \sqrt{\frac{\log(-t)}{(-t)}} \, F \, |F_z| \bigg | \leq \frac{\eta}{8},
\]
if $F = 100\theta \sqrt{-2t}$ and $-t$ is sufficiently large. This implies 
\[
\bigg | (n-2) - \sqrt{\frac{\log(-t)}{(-t)}} \, S^{-\frac{1}{2}}\bigg | \leq \frac{\eta}{2},
\]
if $F = 100\theta \sqrt{-2t}$ and $-t$ is sufficiently large. Because $S^{-\frac{1}{2}}$ is monotone, we conclude that 
\[
\bigg | (n-2) - \sqrt{\frac{\log(-t)}{(-t)}} \, S^{-\frac{1}{2}}\bigg | \leq \frac{\eta}{2},
\]
whenever $L \, \sqrt{\frac{(-t)}{\log(-t)}} \leq F \leq 100\theta \sqrt{-2t}$ and $-t$ is sufficiently large. Finally, this implies 
\[
\bigg | (n-2) - \sqrt{\frac{\log(-t)}{(-t)}} \, F \, |F_z| \bigg | \leq \eta
\]
whenever $L \, \sqrt{\frac{(-t)}{\log(-t)}} \leq F \leq 100\theta \sqrt{-2t}$ and $-t$ is sufficiently large. This completes the proof of Proposition \ref{estimate.for.F_z.in.collar.region}. \\

\begin{proposition}
\label{bound.for.V.in.transition.region}
Let us fix a small number $\theta > 0$. If $-\tau$ is sufficiently large (depending on $\theta$), then $\frac{1}{C(\theta)} \, (-\tau)^{-\frac{1}{2}} \leq V_+(\rho,\tau) \leq C(\theta) \, (-\tau)^{-\frac{1}{2}}$ and $\big | \frac{\partial}{\partial \rho} V_+(\rho,\tau) \big | \leq C(\theta)$ for every $\rho \in [\frac{\theta}{100},100\theta]$.
\end{proposition}

\textbf{Proof.} 
Similarly  to Proposition 4.16 in \cite{BDS}.  \\

In the remainder of this section, we define functions $U_+(r,t)$ and $U_-(r,t)$ so that 
\[U_+(r,t) = \Big ( \frac{\partial}{\partial z} F(z,t) \Big )^2\] 
for $r = F(z,t)$ and $z \geq 2\sqrt{-t}$ and 
\[U_-(r,t) = \Big ( \frac{\partial}{\partial z} F(z,t) \Big )^2\] 
for $r = F(z,t)$ and $z \leq -2\sqrt{-t}$. Let us consider the rescaled functions 
\begin{align*}
V_+(\rho,\tau) &:= \sqrt{U_+(e^{-\frac{\tau}{2}} \rho,-e^{-\tau})}, \\ 
V_-(\rho,\tau) &:= \sqrt{U_-(e^{-\frac{\tau}{2}} \rho,-e^{-\tau})}. 
\end{align*}
For each $\rho \in (0,1)$, we denote by $\xi_+(\rho,\tau)$ the unique positive solution of the equation $F(e^{-\frac{\tau}{2}} \xi,-e^{-\tau}) = e^{-\frac{\tau}{2}} \rho$; moreover, we denote by $\xi_-(\rho,\tau)$ the unique negative solution of the equation $F(e^{-\frac{\tau}{2}} \xi,-e^{-\tau}) = e^{-\frac{\tau}{2}} \rho$. \\

\begin{proposition}
\label{comparison.of.V.with.bryant.soliton.profile}
Fix a small number $\eta>0$. Then we can find a small number $\theta \in (0,\eta)$ (depending on $\eta$) such that, for $-\tau$ sufficiently large, we have 
\[|V_+(\rho,\tau)^{-2} - \Phi((-\tau)^{\frac{1}{2}} \rho)^{-1}| \leq \eta \, (V_+(\rho,\tau)^{-2}-1)\] 
in the region $\{\rho \leq 100\theta\}$. Here, $\Phi$ denotes the profile of the Bryant soliton.
\end{proposition}

\textbf{Proof.} Similarly  to Proposition 4.17 in \cite{BDS}. \\

\begin{proposition}
\label{higher.derivatives.of.V.in.collar.region}
Fix a small number $\eta>0$. Then we can find a large number $L$ (depending on $\eta$) such that, for $-\tau$ sufficiently large, we have 
\[V_+(\rho,\tau) \leq \eta,\] 
\[\Big | \frac{\partial}{\partial \rho} V_+(\rho,\tau) \Big | \leq \eta \, \rho^{-1} \, V_+(\rho,\tau)^{-1}\] 
and 
\[\Big | \frac{\partial^2}{\partial \rho^2} V_+(\rho,\tau) \Big | \leq \eta \, \rho^{-2} \, V_+(\rho,\tau)^{-3}\] 
in the region $\{L \, (-\tau)^{-\frac{1}{2}} \leq \rho \leq \frac{1}{4}\}$.
\end{proposition}

\textbf{Proof.} 
Similarly  to Proposition 4.18 in \cite{BDS}. \\

\begin{corollary}
\label{derivative.of.V.wrt.tau}
Fix a small number $\eta>0$. Then, for $-\tau$ sufficiently large, we have 
\[\Big | \frac{\partial}{\partial \tau} V_+(\rho,\tau) \Big | \leq \eta \, \rho^{-2} \, (V_+(\rho,\tau)^{-1}-1)\] 
in the region $\{\rho \leq \frac{1}{4}\}$. 
\end{corollary}

\textbf{Proof.} 
Similarly  to Proposition 4.19 in \cite{BDS}. \\

\begin{proposition}
\label{derivative.of.xi.wrt.rho}
Fix a small number $\eta>0$. Then we can find a small number $\theta \in (0,\eta)$ (depending on $\eta$) such that, for $-\tau$ sufficiently large, we have 
\[\Big | \frac{\partial}{\partial \rho} \Big ( \frac{\xi_+(\rho,\tau)^2}{4} \Big ) + (n-2) \rho^{-1} \, (V_+(\rho,\tau)^{-2}-1) \Big | \leq \eta \, \rho^{-1} \, (V_+(\rho,\tau)^{-2}-1)\] 
in the region $\{\frac{\theta}{8} \leq \rho \leq 2\theta\}$. 
\end{proposition}

\textbf{Proof.} 
Similarly to Proposition  4.20 in \cite{BDS}. \\

\begin{proposition}
\label{second.derivative.of.xi.wrt.rho}
Fix a small number $\theta>0$. Then, for $-\tau$ large, we have 
\[\Big | \frac{\partial}{\partial \rho} \Big ( \frac{\xi_+(\rho,\tau)^2}{4} \Big ) \Big | \leq C(\theta) \, (-\tau)\] 
and 
\[\Big | \frac{\partial^2}{\partial \rho^2} \Big ( \frac{\xi_+(\rho,\tau)^2}{4} \Big ) \Big | \leq C(\theta) \, (-\tau)^{\frac{3}{2}}\] 
in the region $\{\frac{\theta}{8} \leq \rho \leq 2\theta\}$. 
\end{proposition}

\textbf{Proof.} Similarly  to Proposition 4.21 in \cite{BDS}. \\

\begin{proposition}
\label{derivative.of.xi.wrt.tau}
Fix a small number $\theta>0$. Then, for $-\tau$ large, we have 
\[\Big | \frac{\partial}{\partial \tau} \Big ( \frac{\xi_+(\rho,\tau)^2}{4} \Big ) \Big | \leq o(1) \, (-\tau)\] 
in the region $\{\frac{\theta}{8} \leq \rho \leq 2\theta\}$.
\end{proposition}

\textbf{Proof.}  Similarly  to Proposition 4.22 in \cite{BDS}.  \\

\section{The tip region weights $\mu_+(\rho,\tau)$ and $\mu_-(\rho,\tau)$}

\label{weights}

In this section, we define weights $\mu_+(\rho,\tau)$ and $\mu_-(\rho,\tau)$ which will be needed in the analysis of the linearized equation in the tip region. Let $\theta>0$ be a small positive number, and let $\zeta: \mathbb{R} \to [0,1]$ be a smooth, monotone increasing cutoff function satisfying $\zeta(\rho) = 0$ for $\rho \leq \frac{\theta}{8}$ and $\zeta(\rho) = 1$ for $\rho \geq \frac{\theta}{4}$. We define the weight $\mu_+(\rho,\tau)$ by 
\begin{align*}
\mu_+(\rho,\tau) &= -\zeta(\rho) \, \frac{\xi_+(\rho,\tau)^2}{4} - \int_\rho^\theta \zeta'(\tilde{\rho}) \, \frac{\xi_+(\tilde{\rho},\tau)^2}{4} \, d\tilde{\rho} \\ 
&- (n-2) \int_\rho^\theta (1-\zeta(\tilde{\rho})) \, \tilde{\rho}^{-1} \, \big ( \Phi((-\tau)^{\frac{1}{2}} \, \tilde{\rho})^{-1}-1 \big ) \, d\tilde{\rho},
\end{align*}
where $\Phi$ denotes the profile of the Bryant soliton. We can define a weight $\mu_-(\rho,\tau)$ in analogous fashion. Of course, the cutoff function $\zeta$ and the weights $\mu_+(\rho,\tau)$ and $\mu_-(\rho,\tau)$ depend on the choice of the parameter $\theta$, but we suppress that dependence in our notation. \\

\begin{lemma} 
\label{weight}
The weight $\mu_+(\rho,\tau)$ satisfies $\mu_+(\rho,\tau) = -\frac{\xi_+(\rho,\tau)^2}{4}$ for $\rho \geq \frac{\theta}{4}$. Moreover, $\mu_+(\rho,\tau) \leq 0$ for all $\rho \leq \frac{\theta}{4}$.
\end{lemma}

\textbf{Proof.} 
This follows immediately from the definition of $\mu_+(\rho,\tau)$. \\

\begin{lemma}
\label{derivative.of.mu.wrt.rho}
Fix a small number $\eta>0$. Then we can find a small number $\theta \in (0,\eta)$ (depending on $\eta$) such that, for $-\tau$ sufficiently large, we have 
\[\Big | \frac{\partial \mu_+}{\partial \rho}(\rho,\tau) - (n-2) \rho^{-1} \, (V_+(\rho,\tau)^{-2} - 1) \Big | \leq \eta \, \rho^{-1} \, (V_+(\rho,\tau)^{-2} - 1)\] 
in the tip region $\{\rho \leq 2\theta\}$. 
\end{lemma} 

\textbf{Proof.} 
The  assertion follows  from Proposition  \ref{comparison.of.V.with.bryant.soliton.profile} and Proposition \ref{derivative.of.xi.wrt.rho}
(see  Lemma 5.2 in \cite{BDS} for details).  \\

\begin{lemma}
\label{second.derivative.of.mu.wrt.rho}
If we choose $\theta>0$ sufficiently small, then the following holds. If $-\tau$ is sufficiently large (depending on $\theta$), then 
\[\frac{\partial^2 \mu_+}{\partial \rho^2}(\rho,\tau) \leq \frac{1}{4} \, \Big ( \frac{\partial \mu_+}{\partial \rho}(\rho,\tau) \Big )^2 + \frac{K_*}{4} \, \rho^{-2}\] 
in the tip region $\{\rho \leq 2\theta\}$. Here, $K_*$ is a  constant which depends on the dimension $b$ but is independent of $\theta$.
\end{lemma}

\textbf{Proof.} 
We compute 
\begin{align*} 
\frac{\partial^2 \mu_+}{\partial \rho^2}(\rho,\tau) 
&= -\zeta(\rho) \, \frac{\partial^2}{\partial \rho^2} \Big ( \frac{\xi_+(\rho,\tau)^2}{4} \Big ) - \zeta'(\rho) \, \frac{\partial}{\partial \rho} \Big ( \frac{\xi_+(\rho,\tau)^2}{4} \Big ) \\ 
&- (n-2) [1-\zeta(\rho) + \rho \, \zeta'(\rho)] \, \rho^{-2} \, \big ( \Phi((-\tau)^{\frac{1}{2}} \, \rho)^{-1}-1 \big ) \\ 
&- (n-2) (1-\zeta(\rho)) \, (-\tau)^{\frac{1}{2}} \, \rho^{-1} \, \Phi((-\tau)^{\frac{1}{2}} \, \rho)^{-2} \, \Phi'((-\tau)^{\frac{1}{2}} \, \rho).  
\end{align*}
Recall that $0 \leq \zeta \leq 1$ and $\zeta' \geq 0$. Moreover, we have $\Phi(r)^{-1}-1 \geq \frac{1}{K} \, r^2$ and $|\Phi(r)^{-2} \, \Phi'(r)| \leq Kr$ for all $r \in [0,\infty)$, where $K$ is a universal constant depending only on dimension $n$.  This implies 
\begin{align*} 
\frac{\partial^2 \mu_+}{\partial \rho^2}(\rho,\tau) 
&\leq -\zeta(\rho) \, \frac{\partial^2}{\partial \rho^2} \Big ( \frac{\xi_+(\rho,\tau)^2}{4} \Big ) - \zeta'(\rho) \, \frac{\partial}{\partial \rho} \Big ( \frac{\xi_+(\rho,\tau)^2}{4} \Big ) \\ 
&+ K \, (1-\zeta(\rho)) \, (-\tau), 
\end{align*}
where $K$ is a constant depending on the dimension $n$ but  is independent of $\theta$. Using Proposition \ref{second.derivative.of.xi.wrt.rho}, we obtain 
\[\frac{\partial^2 \mu_+}{\partial \rho^2}(\rho,\tau) \leq o(1) \, (-\tau)^2\] 
in the region $\{\frac{\theta}{8} \leq \rho \leq 2\theta\}$, and 
\[\frac{\partial^2 \mu_+}{\partial \rho^2}(\rho,\tau) \leq K \, (-\tau)\] 
in the region $\{\rho \leq \frac{\theta}{8}\}$. On the other hand, applying Lemma \ref{derivative.of.mu.wrt.rho} with $\eta = \frac{1}{2}$ 
and using that $n-2 \geq 1$ for all $n \geq 4$, we obtain  
\begin{align*} 
\frac{\partial \mu_+}{\partial \rho}(\rho,\tau) 
&\geq \frac{1}{2} \, \rho^{-1} \, (V_+(\rho,\tau)^{-2}-1) \\ 
&\geq \frac{1}{4} \, \rho^{-1} \, (\Phi((-\tau)^{\frac{1}{2}} \rho)^{-1}-1) \\ 
&\geq \frac{1}{4K} \, (-\tau) \, \rho 
\end{align*}
in the region $\{\rho \leq 2\theta\}$, where again $K$ is a constant depending on dimension  but is independent of $\theta$. 
Hence, if $-\tau$ is sufficiently large (depending on $\theta$), then we have 
\[\frac{\partial^2 \mu_+}{\partial \rho^2}(\rho,\tau) \leq \frac{1}{4} \, \Big ( \frac{\partial \mu_+}{\partial \rho}(\rho,\tau) \Big )^2 + 16K^4 \, \rho^{-2}\] 
in the region $\{\rho \leq 2\theta\}$. This completes the proof of Lemma \ref{second.derivative.of.mu.wrt.rho}. \\

\begin{lemma} 
\label{derivative.of.mu.wrt.tau}
Let us fix a small number $\theta>0$. Then, for $-\tau$ large, we have 
\[\Big | \frac{\partial \mu_+}{\partial \tau}(\rho,\tau) \Big | \leq o(1) \, (-\tau)\] 
in the tip region $\{\rho \leq 2\theta\}$.
\end{lemma} 

\textbf{Proof.} 
It follows from Proposition \ref{derivative.of.xi.wrt.tau}, similarly to Proposition 5.3 in \cite{BDS}.  \\

We finish this section  with the following  weighted Poincar\'e inequality.

\begin{proposition}
\label{Poincare.inequality}
If we choose $\theta>0$ sufficiently small, then the following holds. If $-\tau$ is sufficiently large (depending on $\theta$), then 
\[\int_0^{2\theta} \Big ( \frac{\partial \mu_+}{\partial \rho} \Big )^2 \, f^2 \, e^{-\mu_+} \, d\rho \leq 8 \int_0^{2\theta} \Big ( \frac{\partial f}{\partial \rho} \Big )^2 \, e^{-\mu_+} \, d\rho + K_* \int_0^{2\theta} \rho^{-2} \, f^2 \, e^{-\mu_+} \, d\rho\] 
for every smooth function $f$ which is supported in the region $\{\rho \leq 2\theta\}$. Here, $K_*$ is the constant in Lemma \ref{second.derivative.of.mu.wrt.rho}; in particular, $K_*$ depends only on the dimension 
and  is independent of $\theta$. Note that the right hand side is infinite unless $f(0)=0$.
\end{proposition}

\textbf{Proof.} 
We compute 
\[\frac{\partial}{\partial \rho} \Big ( \frac{\partial \mu_+}{\partial \rho} \, f^2 \, e^{-\mu_+} \Big ) = \frac{\partial^2 \mu_+}{\partial \rho^2} \, f^2 \, e^{-\mu_+} + 2 \, \frac{\partial \mu_+}{\partial \rho} \, f \, \frac{\partial f}{\partial \rho} \, e^{-\mu_+} - \Big ( \frac{\partial \mu_+}{\partial \rho} \Big )^2 \, f^2 \, e^{-\mu_+}.\] 
Using Young's inequality, we obtain 
\[\frac{\partial}{\partial \rho} \Big ( \frac{\partial \mu_+}{\partial \rho} \, f^2 \, e^{-\mu_+} \Big ) \leq \frac{\partial^2 \mu_+}{\partial \rho^2} \, f^2 \, e^{-\mu_+} + 2 \, \Big ( \frac{\partial f}{\partial \rho} \Big )^2 \, e^{-\mu_+} - \frac{1}{2} \, \Big ( \frac{\partial \mu_+}{\partial \rho} \Big )^2 \, f^2 \, e^{-\mu_+}.\] 
Hence, Lemma \ref{second.derivative.of.mu.wrt.rho} gives 
\[\frac{\partial}{\partial \rho} \Big ( \frac{\partial \mu_+}{\partial \rho} \, f^2 \, e^{-\mu_+} \Big ) \leq 2 \, \Big ( \frac{\partial f}{\partial \rho} \Big )^2 \, e^{-\mu_+} - \frac{1}{4} \, \Big ( \frac{\partial \mu_+}{\partial \rho} \Big )^2 \, f^2 \, e^{-\mu_+} + \frac{K_*}{4} \, \rho^{-2} \, f^2 \, e^{-\mu_+}.\] 
From this, the assertion follows. \\

\section{Overview of the proof of Theorem \ref{uniqueness.theorem}}

\label{overview}

In this section, we state the four main estimates needed for the proof of Theorem \ref{uniqueness.theorem}, generalizing Section 6 in \cite{BDS}. At the end of this section, we give the proof of Theorem \ref{uniqueness.theorem} assuming these key results. To that end, we consider two ancient $\kappa$-solutions, $(S^n, g_1(t))$ and $(S^n, g_2(t))$ such that neither solution is a family of shrinking round spheres. By the main result of the previous section, we know both solutions are rotationally symmetric. We first choose reference points $q_1, q_2 \in S^n$ such that
\[
\limsup_{t \to -\infty}\, (-t) R_{g_1(t)}(q_1) \leq 50(n-1) \quad \text{and} \quad \limsup_{t \to -\infty}\, (-t) R_{g_2(t)}(q_2) \leq 50(n-1). 
\]
The existence of these points is ensured by the Neck Stability Theorem of Kleiner and Lott. See Proposition 3.1 in \cite{Angenent-Brendle-Daskalopoulos-Sesum} for a proof in dimension three, which also works in higher dimensions. 

Since $(S^n, g_1(t))$ is rotationally symmetric, we can define a profile function $F_1(z,t)$ to be the radius of the sphere of symmetry that has signed distance $z$ from the reference point $q_1$. Similarly we can define $F_2(z,t)$ on $(S^n, g_2(t))$ with respect to $q_2$. These functions, $F_1(z,t)$ and $F_2(z,t)$, satisfy the PDE 
\begin{align*} 
F_t(z,t)  
= F_{zz}(z,t) &- (n-2)F(z,t)^{-1} \, (1-F_z(z,t)^2) \\ 
& - (n-1) \, F_z(z,t) \int_0^z \frac{F_{zz}(z',t)}{F(z',t)} \, dz'. 
\end{align*} 

Our goal is to show that the profile functions $F_1$ and $F_2$ will agree after a reparametrization in space, a translation in time, and a parabolic rescaling. We thus will now define a new function $F_2^{\alpha\beta\gamma}(z,t)$ obtained from $F_2(z,t)$ through a spacial reparametrization, a time translation, and a parabolic rescaling. Here, $(\alpha,\beta,\gamma)$ is a triplet of real numbers satisfying the following admissibility condition previously defined in \cite{BDS}: 

\begin{definition}
\label{admissibility}
Given a real number $\varepsilon \in (0,1)$, we say that the triplet $(\alpha,\beta,\gamma)$ is $\varepsilon$-admissible with respect to time $t_*$ if 
\[|\alpha| \leq \varepsilon \sqrt{-t_*}, \qquad |\beta| \leq \varepsilon \, \frac{(-t_*)}{\log(-t_*)}, \qquad |\gamma| \leq \varepsilon \log (-t_*).\]
\end{definition}

Consider a time $t_\ast < 0$ so that $-t_\ast$ is very large. Suppose $(\alpha, \beta, \gamma)$ is a triplet of real numbers satisfying the criteria of $\varepsilon$-admissibility with respect to the time $t_\ast$, for some $\varepsilon \in (0, 1)$. For each $t \leq t_\ast$, we define a time-translated and parabolically-rescaled metric by 
\[
g_2^{\beta\gamma}(t) := e^{\gamma} g_2(e^{-\gamma}(t - \beta)). 
\]
Of course, $(S^n, g_2^{\beta\gamma}(t))$ is again a rotationally symmetric ancient $\kappa$-solution. We define the time-translated and parabolically-rescaled profile function $F_2^{\beta\gamma}(z,t)$ on the ancient $\kappa$-solution $(S^n, g_2^{\beta\gamma}(t))$ to be the radius of the sphere of symmetry with signed distance $z$ from the reference point $q_2$. Evidently, 
\[
F_2^{\beta\gamma}(z, t) = e^{\frac{\gamma}{2}}F_2\big(e^{-\frac{\gamma}{2}}z, e^{-\gamma}(t - \beta)\big).
\]

Even after a time translation and a parabolic rescaling, it is possible for the profile functions to differ by a translation in space. To account for this, we define a new reference point $q^{\alpha\beta\gamma}_2$ with the property that $q^{\alpha\beta\gamma}_2$ has signed distance $\alpha$ from the original reference point $q_2$ with respect to the metric $g^{\beta\gamma}_2(t_\ast)$. For $t \leq t_\ast$, we define a function $s^{\alpha\beta\gamma}(t)$ to be the signed distance between the sphere of symmetry through $q_2^{\alpha\beta\gamma}$ and the point $q_2$, with respect to $g_2^{\beta\gamma}(t)$. The function $s^{\alpha\beta\gamma}(t)$ is the unique solution of the ODE 
\[
\frac{d}{dt} s^{\alpha\beta\gamma}(t) = (n-1) \int_{0}^{s^{\alpha\beta\gamma}(t)} \frac{F_{2,zz}^{\beta\gamma}(z', t)}{F_2^{\beta\gamma}(z', t)} \, dz', \qquad s^{\alpha\beta\gamma}(t_\ast) = \alpha, 
\]
for $t \leq t_\ast$. This ODE, of course, is just the usual evolution of distance along the Ricci flow and the integrand is the radial component of the Ricci curvature. Now, for $t \leq t_\ast$, we define $F_2^{\alpha\beta\gamma}(z, t)$ to be the radius of the sphere of symmetry in $(S^{n}, g_2^{\beta\gamma}(t))$ which has signed distance $z$ from the point $q_2^{\alpha\beta\gamma}$. The three profile functions are related by the equation
\[
F_2^{\alpha\beta\gamma}(z,t) = F_2^{\beta\gamma}(z + s^{\alpha\beta\gamma}(t), t) = e^{\frac{\gamma}{2}}F_2\big(e^{-\frac{\gamma}{2}}(z+ s^{\alpha\beta\gamma}(t)), e^{-\gamma}(t - \beta)\big). 
\]
In particular, at time $t = t_\ast$, we have 
\[
F_2^{\alpha\beta\gamma}(z,t_\ast) = F_2^{\beta\gamma}(z + \alpha, t) = e^{\frac{\gamma}{2}}F_2\big(e^{-\frac{\gamma}{2}}(z+ \alpha), e^{-\gamma}(t - \beta)\big). 
\]

In the next lemma, we show that for a $\varepsilon$-admissible triplet $(\alpha,\beta,\gamma)$ at $t_\ast$, we expect the new reference point $q^{\alpha\beta\gamma}_2$ to remain suitably close to the original point $q_2$ for all earlier times $t \leq t_\ast$.

\begin{lemma} \label{bound.for.s^alpha}
If $-t_\ast$ is sufficiently large, then the following holds. Suppose the triplet $(\alpha,\beta,\gamma)$ is $\varepsilon$-admissible with respect to time $t_\ast$, where $\varepsilon \in (0, 1)$. Let $s^{\alpha\beta\gamma}(t)$ be the solution of the ODE
\[
\frac{d}{dt} s^{\alpha\beta\gamma}(t) = (n-1) \int_{0}^{s^{\alpha\beta\gamma}(t)} \frac{F_{2,zz}^{\beta\gamma}(z', t)}{F_2^{\beta\gamma}(z', t)} \, dz'
\]
with terminal condition $s^{\alpha\beta\gamma}(t_\ast) = \alpha$. Then $|s^{\alpha\beta\gamma}(t)| \leq \varepsilon \sqrt{-t}$ for all $t \leq t_\ast$. 
\end{lemma}

\textbf{Proof.} The proof is essentially the same as the proof of Lemma 6.2 in \cite{BDS}. Recall that if we rescale the ancient $\kappa$-solution $(S^n, g_2(t))$ around the reference point $q_2$ by the factor $(-t)$, then the solution converges to a round cylinder in pointed Cheeger-Gromov sense at $t\to -\infty$. In particular, in the region $|z| \leq \sqrt{-2t}$, the radial component of the Ricci curvature tends to zero. Consequently, if $-t_\ast$ is sufficiently large, then we will have 
\[
0 \leq - \frac{F_{2,zz}(z,t)}{F_2(z,t)} \leq \frac{1}{(-4(n-1)t)}
\]
whenever $t \leq \frac{1}{2} e^{-\gamma} t_\ast$ and $|z| \leq \sqrt{-2t}$. The first inequality follows from nonnegativity of the Ricci curvature. For the profile function $F^{\beta\gamma}(z,t)$ we replace $t$ by $e^{-\gamma}(t-\beta)$ and $z$ by $e^{-\frac{\gamma}{2}}z$. This implies 
\[
0 \leq - \frac{F^{\beta\gamma}_{2,zz}(z,t)}{F^{\beta\gamma}_2(z,t)} \leq \frac{1}{(-4(n-1)(t-\beta))}. 
\]
whenever $t - \beta \leq \frac{1}{2} t_\ast$ and $|z| \leq \sqrt{-2(t-\beta)}$. By admissibility, $|\beta| \leq \varepsilon \frac{(-t_\ast)}{\log(-t_\ast)} \leq \varepsilon \frac{(-t)}{\log(-t)} \leq \frac{1}{2}(-t)$ whenever $t \leq t_\ast$ and $-t_\ast$ is sufficiently large. This ensures that $2t \leq t -\beta \leq \frac{1}{2} t$ whenever $t \leq t_\ast$. Consequently, 
\[
0 \leq - \frac{F^{\beta\gamma}_{2,zz}(z,t)}{F^{\beta\gamma}_2(z,t)} \leq \frac{1}{(-2(n-1)t)}. 
\]
whenever $t \leq t_\ast$ and $|z| \leq \sqrt{-t}$. Plugging this estimate into the ODE for $s^{\alpha\beta\gamma}(t)$, we obtain
\[
\bigg|\frac{d}{dt}s^{\alpha\beta\gamma}(t)\bigg| \leq \frac{1}{(-2t)}|s^{\alpha\beta\gamma}(t)|,
\]
whenever $t \leq t_\ast$ and $|s^{\alpha\beta\gamma}(t)| \leq \sqrt{-t}$. This gives
\[
\frac{d}{dt} \Big(\frac{|s^{\alpha\beta\gamma}(t)|}{\sqrt{-t}}\Big)\geq 0
\]
whenever $t \leq t_\ast$ and $|s^{\alpha\beta\gamma}(t)| \leq \sqrt{-t}$. At time $t_\ast$, we have $|s^{\alpha\beta\gamma}(t_\ast)| = |\alpha| \leq \varepsilon \sqrt{-t_\ast}$. The the differential inequality above implies $|s^{\alpha\beta\gamma}(t)| \leq \varepsilon \sqrt{-t}$. This completes the proof of Lemma \ref{bound.for.s^alpha}.

Using the admissibility conditions for $(\alpha,\beta,\gamma)$ and the previous lemma, we can estimate the profile function $F^{\alpha\beta\gamma}_2$:

\begin{proposition}
\label{estimate.for.modified.profile}
Fix a small number $\theta > 0$ and a small number $\eta > 0$. Then there exists a small number $\varepsilon>0$ (depending on $\theta$ and $\eta$) with the following property. If the triplet $(\alpha,\beta,\gamma)$ is $\varepsilon$-admissible with respect to time $t_*$ and $-t_*$ is sufficiently large, then  
\[\Big | \frac{1}{2} \, F_2^{\alpha\beta\gamma}(z,t)^2 +(n-2)t + (n-2)\frac{z^2+2t}{4 \log(-t)} \Big | \leq \eta \, \frac{z^2-t}{\log(-t)}\] 
and 
\[\Big | F_2^{\alpha\beta\gamma}(z,t) \, F_{2z}^{\alpha\beta\gamma}(z,t) + \frac{(n-2)z}{2 \log(-t)} \Big | \leq \eta \, \frac{|z|+\sqrt{-t}}{\log(-t)}\] 
whenever $F_2^{\alpha\beta\gamma}(z,t) \geq \frac{\theta}{10} \sqrt{-t}$ and $t \leq t_*$.
\end{proposition}

\textbf{Proof.} The proof is essentially the same as the proof of the Proposition 6.3 in \cite{BDS}. 
Using Proposition \ref{precise.estimate.for.F} and Proposition \ref{precise.estimate.for.F_z}, we obtain 
\[\Big | \frac{1}{2} \, F_2(z,t)^2 + (n-2)t + (n-2)\frac{z^2+2t}{4 \log(-t)} \Big | \leq \frac{\eta}{4} \, \frac{z^2-t}{\log(-t)}\] 
and 
\[\Big | F_2(z,t) \, F_{2z}(z,t) + \frac{(n-2)z}{2 \log(-t)} \Big | \leq \frac{\eta}{4} \, \frac{|z|+\sqrt{-t}}{\log(-t)}\] 
whenever $F_2(z,t) \geq \frac{\theta}{20} \sqrt{-t}$ and $-t$ is sufficiently large. To estimate $F^{\beta\gamma}$, we replace $t$ by $e^{-\gamma} (t-\beta)$ and $z$ by $e^{-\frac{\gamma}{2}} z$. This gives 
\[\Big | \frac{1}{2} \, F_2^{\beta\gamma}(z,t)^2 + (n-2)(t-\beta) + (n-2)\frac{z^2+2(t-\beta)}{4 \log(-(t-\beta)) - 4\gamma} \Big | \leq \frac{\eta}{4} \, \frac{z^2-(t-\beta)}{\log(-(t-\beta))-\gamma}\] 
and 
\[\Big | F_2^{\beta\gamma}(z,t) \, F_{2z}^{\beta\gamma}(z,t) + \frac{(n-2)z}{2 \log(-(t-\beta)) - 2\gamma} \Big | \leq \frac{\eta}{4} \, \frac{|z|+\sqrt{-(t-\beta)}}{\log(-(t-\beta))-\gamma}\] 
whenever $F_2^{\beta\gamma}(z,t) \geq \frac{\theta}{20} \sqrt{-(t-\beta)}$ and $-e^{-\gamma} (t-\beta)$ is sufficiently large. The $\varepsilon$-admissibility assumptions on $(\alpha,\beta,\gamma)$ at time $t_\ast$ ensure that $|\beta| \leq \varepsilon \frac{(-t)}{\log(-t)}$ and $\gamma \leq \varepsilon \log(-t)$ for $t \leq t_\ast$. If $\varepsilon$ is sufficiently small (depending on $\theta$ and $\eta$) and $-t_*$ is sufficiently large (depending on $\theta$ and $\eta$), then we obtain 
\[\Big | \frac{1}{2} \, F_2^{\beta\gamma}(z,t)^2 + (n-2)t + (n-2)\frac{z^2+2t}{4 \log(-t)} \Big | \leq \frac{\eta}{2} \, \frac{z^2-t}{\log(-t)}\] 
and 
\[\Big | F_2^{\beta\gamma}(z,t) \, F_{2z}^{\beta\gamma}(z,t) + \frac{(n-2)z}{2 \log(-t)} \Big | \leq \frac{\eta}{2} \, \frac{|z|+\sqrt{-t}}{\log(-t)}\] 
whenever $F_2^{\beta\gamma}(z,t) \geq \frac{\theta}{10} \sqrt{-t}$ and $t \leq t_*$. By Lemma \ref{bound.for.s^alpha}, $|s^{\alpha\beta\gamma}(t)| \leq \varepsilon \sqrt{-t}$ for $t \leq t_*$. Replacing $z$ by $z + s^{\alpha\beta\gamma}(t)$ and using this estimate, we obtain 
\[\Big | \frac{1}{2} \, F_2^{\alpha\beta\gamma}(z,t)^2 + (n-2)t + (n-2)\frac{z^2+2t}{4 \log(-t)} \Big | \leq \eta \, \frac{z^2-t}{\log(-t)}\] 
and 
\[\Big | F_2^{\alpha\beta\gamma}(z,t) \, F_{2z}^{\alpha\beta\gamma}(z,t) + \frac{(n-2)z}{2 \log(-t)} \Big | \leq \eta \, \frac{|z|+\sqrt{-t}}{\log(-t)}\] 
whenever $F_2^{\alpha\beta\gamma}(z,t) \geq \frac{\theta}{10} \sqrt{-t}$ and $t \leq t_*$. This completes the proof of Proposition \ref{estimate.for.modified.profile}. \\

As in \cite{BDS}, we need to use different functions to describe our solutions and establish estimates in the tip regions. These functions labeled by $U$ are analogous to the profile function $\Phi$ used to describe the Bryant soliton.  As in \cite{BDS}, we define functions $U_{1+}(r,t)$ and $U_{1-}(r,t)$ by 
\[
U_{1+}(r,t) = \Big ( \frac{\partial}{\partial z} F_1(z,t) \Big )^2
\] 
for $r = F_1(z,t)$ and $z \geq 2\sqrt{-t}$ and 
\[
U_{1-}(r,t) = \Big ( \frac{\partial}{\partial z} F_1(z,t) \Big )^2
\] 
for $r = F_1(z,t)$ and $z \leq -2\sqrt{-t}$. Similarly, we define functions $U_{2+}(r,t)$ and $U_{2-}(r,t)$ by 
\[
U_{2+}(r,t) = \Big ( \frac{\partial}{\partial z} F_2(z,t) \Big )^2
\] 
for $r = F_2(z,t)$ and $z \geq 2\sqrt{-t}$ and 
\[
U_{2-}(r,t) = \Big ( \frac{\partial}{\partial z} F_2(z,t) \Big )^2
\] 
for $r = F_2(z,t)$ and $z \leq -2\sqrt{-t}$. 
Then, we define 
\begin{align*} 
U_{2+}^{\beta\gamma}(r,t) &:= U_{2+}(e^{-\frac{\gamma}{2}} r,e^{-\gamma} (t-\beta)), \\ 
U_{2-}^{\beta\gamma}(r,t) &:= U_{2-}(e^{-\frac{\gamma}{2}} r,e^{-\gamma} (t-\beta)). 
\end{align*} 
Recalling that $F^{\alpha\beta\gamma}(z, t) = F^{\beta\gamma}(z + s^{\alpha\beta\gamma}(t), t)$, observe that 
\[
U_{2+}^{\beta\gamma}(r,t) = \Big ( \frac{\partial}{\partial z} F_2^{\alpha\beta\gamma}(z,t) \Big )^2
\] 
for $r = F_2(z,t)$, $z \geq 4\sqrt{-t}$, and $t \leq t_*$, and 
\[
U_{2-}^{\beta\gamma}(r,t) = \Big ( \frac{\partial}{\partial z} F_2^{\alpha\beta\gamma}(z,t) \Big )^2
\] 
for $r = F_2(z,t)$, $z \leq -4\sqrt{-t}$, and $t \leq t_*$. 

For each of the functions above, we will define a function $V$ in the usual rescaled coordinates. For scaling reasons, it is convenient to define the functions labeled by $V$ to be the square-root of the corresponding functions labeled by $U$. As usual, define coordinates $\tau$ and $\rho$ by the identities $t = -e^{-\tau}$ and $r = e^{-\frac{\tau}{2}} \rho$. Then we define:
\begin{align*} 
V_{1+}(\rho,\tau) &:= \sqrt{U_{1+}(e^{-\frac{\tau}{2}} \rho,-e^{-\tau})}, \\ 
V_{1-}(\rho,\tau) &:= \sqrt{U_{1-}(e^{-\frac{\tau}{2}} \rho,-e^{-\tau})}, \\ 
V_{2+}(\rho,\tau) &:= \sqrt{U_{2+}(e^{-\frac{\tau}{2}} \rho,-e^{-\tau})}, \\ 
V_{2-}(\rho,\tau) &:= \sqrt{U_{2-}(e^{-\frac{\tau}{2}} \rho,-e^{-\tau})}, \\ 
V_{2+}^{\beta\gamma}(\rho,\tau) &:= \sqrt{U_{2+}^{\beta\gamma}(e^{-\frac{\tau}{2}} \rho,-e^{-\tau})}, \\ 
V_{2-}^{\beta\gamma}(\rho,\tau) &:= \sqrt{U_{2-}^{\beta\gamma}(e^{-\frac{\tau}{2}} \rho,-e^{-\tau})}. 
\end{align*} 
By solving $-e^{-\tilde \tau} = e^{-\gamma}(-e^{-\tau} - \beta)$ and $e^{-\tilde \tau/2} \tilde \rho = e^{-\gamma/2} e^{-\tau/2} \rho$ for $\tilde \tau$ and $\tilde \rho$, you can confirm that 
\begin{align*} 
V_{2+}^{\beta\gamma}(\rho,\tau) &= V_{2+} \Big ( \frac{\rho}{\sqrt{1+\beta e^\tau}}, \tau + \gamma - \log(1+\beta e^\tau) \Big ), \\ 
V_{2-}^{\beta\gamma}(\rho,\tau) &= V_{2-} \Big ( \frac{\rho}{\sqrt{1+\beta e^\tau}}, \tau + \gamma - \log(1+\beta e^\tau) \Big ). 
\end{align*}

In the following proposition, we recover a version of the estimates established in Proposition \ref{comparison.of.V.with.bryant.soliton.profile} and Corollary \ref{derivative.of.V.wrt.tau} for the modified profile functions $V^{\beta\gamma}_{2+}$ and $V^{\beta\gamma}_{2-}$. 

\begin{proposition}
\label{comparison.of.modified.V.with.bryant.soliton.profile}
Fix a small number $\eta>0$. Then we can find a small number $\theta \in (0,\eta)$ (depending on $\eta$) and a small number $\varepsilon>0$ (depending on $\theta$ and $\eta$) with the following property. If the triplet $(\alpha,\beta,\gamma)$ is $\varepsilon$-admissible with respect to time $t_* = -e^{-\tau_*}$ and $-\tau_*$ is sufficiently large, then 
\[|V_{2+}^{\beta\gamma}(\rho,\tau)^{-2} - \Phi((-\tau)^{\frac{1}{2}} \rho)^{-1}| \leq \eta \, (V_{2+}^{\beta\gamma}(\rho,\tau)^{-2}-1)\] 
for $\rho \leq 10\theta$ and $\tau \leq \tau_*$, and 
\[\Big | \frac{\partial}{\partial \tau} V_{2+}^{\beta\gamma}(\rho,\tau)^{-2} \Big | \leq \eta \, \rho^{-2} \, (V_{2+}^{\beta\gamma}(\rho,\tau)^{-1}-1)\] 
for $\rho \leq \frac{1}{8}$ and $\tau \leq \tau_*$. Here, $\Phi$ denotes the profile of the Bryant soliton.
\end{proposition}

\textbf{Proof.} The proof is identical to the proof of Proposition 6.4 in \cite{BDS}.  \\


We next consider the difference between the two solutions near each of the tips: 
\begin{align*} 
W_+^{\beta\gamma}(\rho,\tau) &:= V_{1+}(\rho,\tau) - V_{2+}^{\beta\gamma}(\rho,\tau), \\ 
W_-^{\beta\gamma}(\rho,\tau) &:= V_{1-}(\rho,\tau) - V_{2-}^{\beta\gamma}(\rho,\tau). 
\end{align*} 
For each $\tau$, we have $W_+^{\beta\gamma}(\rho,\tau) = O(\rho^2)$ and $W_-^{\beta\gamma}(\rho,\tau) = O(\rho^2)$ as $\rho \to 0$. Moreover, let $\mu_+(\rho,\tau)$ and $\mu_-(\rho,\tau)$ denote the weights associated with the solution $(S^n, g_1(t))$.  The following proposition is the first of four key estimates. \\

\begin{proposition}
\label{estimate.for.difference.in.tip.region}
We can choose $\theta>0$ and $\varepsilon>0$ sufficiently small so that the following holds. If $-\tau_*$ is sufficiently large (depending on $\theta$) and the triplet $(\alpha,\beta,\gamma)$ is $\varepsilon$-admissible with respect to time $t_* = -e^{-\tau_*}$, then 
\begin{align*} 
&\sup_{\tau \leq \tau_*} (-\tau)^{-\frac{1}{2}} \int_{\tau-1}^\tau \int_0^\theta V_{1+}^{-2} \, (W_+^{\beta\gamma})^2 \, e^{\mu_+} \\ 
&\leq C(\theta) \, (-\tau_*)^{-1} \sup_{\tau \leq \tau_*} (-\tau)^{-\frac{1}{2}} \int_{\tau-1}^\tau \int_\theta^{2\theta} V_{1+}^{-2} \, (W_+^{\beta\gamma})^2 \, e^{\mu_+}. 
\end{align*} 
An analogous estimate holds for $W_-^{\beta\gamma}$.
\end{proposition}

We will give the proof of Proposition \ref{estimate.for.difference.in.tip.region} in Section \ref{difference.tip.region}. \\

From this point on, we fix $\theta$ small enough so that the conclusion of Proposition \ref{estimate.for.difference.in.tip.region} holds. Let $\chi_{\mathcal{C}}$ denote a smooth, even cutoff function satisfying $\chi_{\mathcal{C}} = 1$ on $[0,\sqrt{4-\frac{\theta^2}{2(n-2)}}]$ and $\chi_{\mathcal{C}} = 0$ on $[\sqrt{4-\frac{\theta^2}{4(n-2)}},\infty)$. Having fixed $\theta$ as in Proposition \ref{estimate.for.difference.in.tip.region}, the factor of $(n-2)$ in higher dimensions ensures there is overlap between the tip region and the collar region where estimates can be played off one another. Moreover, we may assume that $\chi_{\mathcal{C}}$ is monotone decreasing on $[0,\infty)$. 

We define the rescaled profile functions 
\begin{align*} 
G_1(\xi,\tau) &:= e^{\frac{\tau}{2}} \, F_1(e^{-\frac{\tau}{2}} \xi,-e^{-\tau}) - \sqrt{2(n-2)}, \\ 
G_2(\xi,\tau) &:= e^{\frac{\tau}{2}} \, F_2(e^{-\frac{\tau}{2}} \xi,-e^{-\tau}) - \sqrt{2(n-2)}, \\ 
G_2^{\alpha\beta\gamma}(\xi,\tau) &:= e^{\frac{\tau}{2}} \, F_2^{\alpha\beta\gamma}(e^{-\frac{\tau}{2}} \xi,-e^{-\tau}) - \sqrt{2(n-2)}. 
\end{align*} 
Then we consider the difference of the rescaled profile functions in the collar region via
\[
H^{\alpha\beta\gamma}(\xi,\tau) := G_1(\xi,\tau) - G_2^{\alpha\beta\gamma}(\xi,\tau)\] 
and 
\[
H_{\mathcal{C}}^{\alpha\beta\gamma}(\xi,\tau) := \chi_{\mathcal{C}}((-\tau)^{-\frac{1}{2}} \xi) \, H^{\alpha\beta\gamma}(\xi,\tau).
\] 
Using the PDEs for $G_1$ and $G_2^{\alpha\beta\gamma}$, we can derive a PDE for the function $H^{\alpha\beta\gamma}$. As in three dimensions, the leading term in that PDE is given by the operator 
\[
\mathcal{L} f := f_{\xi\xi} - \frac{1}{2} \, \xi \, f_\xi + f.
\]
We will analyze this operator as in \cite{BDS}.  We consider the Hilbert space $\mathcal{H} = L^2(\mathbb{R},e^{-\frac{\xi^2}{4}} \, d\xi)$ and recall that the Hilbert space $\mathcal{H}$ has a natural direct sum decomposition $\mathcal{H} = \mathcal{H}_+ \oplus \mathcal{H}_0 \oplus \mathcal{H}_-$. Furthermore, we recall that $\mathcal{H}_+$ is a two-dimensional subspace spanned by the functions $1$ and $\xi$; $\mathcal{H}_0$ is a one-dimensional subspace spanned by the function $\xi^2-2$; and $\mathcal{H}_-$ is the orthogonal complement of $\mathcal{H}_+ \oplus \mathcal{H}_0$. Finally, let $P_+$, $P_0$, and $P_-$ denote the projection operators associated to the direct sum decomposition $\mathcal{H} = \mathcal{H}_+ \oplus \mathcal{H}_0 \oplus \mathcal{H}_-$.

With these conventions, we write 
\[P_0 H_{\mathcal{C}}^{\alpha\beta\gamma}(\xi,\tau) = \sqrt{2(n-2)} \, a^{\alpha\beta\gamma}(\tau) \, (\xi^2-2),\] 
where 
\[a^{\alpha\beta\gamma}(\tau) := \frac{1}{16\sqrt{2(n-2)\pi}} \int_{\mathbb{R}} e^{-\frac{\xi^2}{4}} \, (\xi^2-2) \, H_{\mathcal{C}}^{\alpha\beta\gamma}(\xi,\tau) \, d\xi.\] 
Moreover, we let $\hat{H}_{\mathcal{C}}^{\alpha\beta\gamma} = P_+ H_{\mathcal{C}}^{\alpha\beta\gamma} + P_- H_{\mathcal{C}}^{\alpha\beta\gamma}$ denote the sum of projections onto the spaces of positive and negative modes.  \\

In the following proposition, we use our freedom of choice in the parameters $(\alpha,\beta,\gamma)$ to ensure the projections our solution $P_0H^{\alpha\beta\gamma}_{\mathcal{C}}$ and $P_+H^{\alpha\beta\gamma}_{\mathcal{C}}$ (i.e. the projections onto the spaces of non-decaying modes of the operator $\mathcal L$) vanish at a particular time $\tau_\ast$. 
\begin{proposition}
\label{choice.of.parameters}
Fix $\theta>0$ and $\varepsilon>0$ small enough so that the conclusion of Proposition \ref{estimate.for.difference.in.tip.region} holds. Let $\delta \in (0,\varepsilon)$ be given. If $-\tau_*$ is sufficiently large (depending on $\delta$), then we can find a triplet $(\alpha,\beta,\gamma)$ (depending on $\tau_*$) such that $P_+ H_{\mathcal{C}}^{\alpha\beta\gamma} = 0$ and $P_0 H_{\mathcal{C}}^{\alpha\beta\gamma} = 0$ at time $\tau_*$. Moreover, if $-\tau_*$ is sufficiently large (depending on $\delta$), then the triplet $(\alpha,\beta,\gamma)$ is $\delta$-admissible with respect to time $t_* = -e^{-\tau_*}$.
\end{proposition}

\textbf{Proof.} 
By definition, $s^{\alpha\beta\gamma}(t_*) = \alpha$. Hence
\[
F_2^{\alpha\beta\gamma}(z,t_*) = e^{\frac{\gamma}{2}} \, F_2(e^{-\frac{\gamma}{2}} (z+\alpha),e^{-\gamma} (t_*-\beta)).
\] 
It follows by a straightforward computation that 
\begin{align*} 
G_2^{\alpha\beta\gamma}(\xi,\tau_*) 
&= \sqrt{1+\beta e^{\tau_*}} \, G_2 \Big ( \frac{\xi + \alpha e^{\frac{\tau_*}{2}}}{\sqrt{1+\beta e^{\tau_*}}},\tau_* + \gamma - \log(1+\beta e^{\tau_*}) \Big ) \\ 
&+ \sqrt{2(n-2)} \, (\sqrt{1+\beta e^{\tau_*}} - 1). 
\end{align*} 
The proof of Proposition \ref{choice.of.parameters} now proceeds as in \cite{Angenent-Daskalopoulos-Sesum2}. This argument relies only on the asymptotics of our solution in the cylindrical region. Since the asymptotics of our ancient solutions to Ricci flow in the cylindrical region are very similar to the cylindrical region asymptotics of ancient solutions to mean curvature flow, the proof of Proposition \ref{choice.of.parameters} is identical to the proof of the corresponding Proposition 4.1 in \cite{Angenent-Daskalopoulos-Sesum2}. \\

From this point on, we assume that the triplet $(\alpha,\beta,\gamma)$ is chosen as in Proposition \ref{choice.of.parameters}, pending our choice of $\tau_\ast$ (which we have not yet fixed). In particular, this will ensure that $a^{\alpha\beta\gamma}(\tau_*) = 0$.  \\

We can now state the remaining three key estimates used in completing the proof of Theorem \ref{uniqueness.theorem}. The first estimate is an estimate for the difference of the solutions in the cylindrical region. 
\begin{proposition}
\label{estimate.for.difference.in.cylindrical.region}
Fix $\theta>0$ small enough so that the conclusion of Proposition \ref{estimate.for.difference.in.tip.region} holds. Suppose that $-\tau_*$ is sufficiently large, and that the triplet $(\alpha,\beta,\gamma)$ is chosen as in Proposition \ref{choice.of.parameters}. Then 
\begin{align*}
&(-\tau_*) \sup_{\tau \leq \tau_*} \int_{\tau-1}^\tau \int_{\mathbb{R}} e^{-\frac{\xi^2}{4}} \, (\hat{H}_{\mathcal{C},\xi}^{\alpha\beta\gamma}(\xi,\tau')^2 + \hat{H}_{\mathcal{C}}^{\alpha\beta\gamma}(\xi,\tau')^2) \, d\xi \, d\tau' \\ 
&\leq C(\theta) \sup_{\tau \leq \tau_*} \int_{\tau-1}^\tau a^{\alpha\beta\gamma}(\tau')^2 \, d\tau' \\ 
&+ C(\theta) \sup_{\tau \leq \tau_*} \int_{\tau-1}^\tau \int_{\{\sqrt{4-\frac{\theta^2}{2(n-2)}} \, (-\tau')^{\frac{1}{2}} \leq |\xi| \leq \sqrt{4-\frac{\theta^2}{4(n-2)}} \, (-\tau')^{\frac{1}{2}}\}} e^{-\frac{\xi^2}{4}} \, H^{\alpha\beta\gamma}(\xi,\tau')^2 \, d\xi \, d\tau'. 
\end{align*}
\end{proposition}

We will give the proof of Proposition \ref{estimate.for.difference.in.cylindrical.region} in Section \ref{difference.cylindrical.region}. \\

Next, by combining Proposition \ref{estimate.for.difference.in.tip.region} and Proposition \ref{estimate.for.difference.in.cylindrical.region}, we can show that in the cylindrical region the norm of $P_0 H_{\mathcal{C}}^{\alpha\beta\gamma}$ dominates over the norm of $\hat{H}_{\mathcal{C}}^{\alpha\beta\gamma}$. More precisely, we have the following result: 

\begin{proposition}
\label{neutral.mode.dominates}
Fix $\theta>0$ small enough so that the conclusion of Proposition \ref{estimate.for.difference.in.tip.region} holds. Suppose that $-\tau_*$ is sufficiently large, and that the triplet $(\alpha,\beta,\gamma)$ is chosen as in Proposition \ref{choice.of.parameters}. Then 
\begin{align*}
&(-\tau_*) \sup_{\tau \leq \tau_*} \int_{\tau-1}^\tau \int_{\mathbb{R}} e^{-\frac{\xi^2}{4}} \, (\hat{H}_{\mathcal{C},\xi}^{\alpha\beta\gamma}(\xi,\tau')^2 + \hat{H}_{\mathcal{C}}^{\alpha\beta\gamma}(\xi,\tau')^2) \, d\xi \, d\tau' \\ 
&\leq C(\theta) \sup_{\tau \leq \tau_*} \int_{\tau-1}^\tau a^{\alpha\beta\gamma}(\tau')^2 \, d\tau'. 
\end{align*}
\end{proposition}

The proof of Proposition \ref{neutral.mode.dominates} will be given in Section \ref{overlap.region}. \\

Using Proposition \ref{neutral.mode.dominates}, we are able to derive an ODE for the function $a^{\alpha\beta\gamma}(\tau)$:

\begin{proposition}
\label{ode.for.a}
Fix $\theta>0$ small enough so that the conclusion of Proposition \ref{estimate.for.difference.in.tip.region} holds. Let $\delta > 0$ be given. Suppose that $-\tau_*$ is sufficiently large (depending on $\delta$), and the triplet $(\alpha,\beta,\gamma)$ is chosen as in Proposition \ref{choice.of.parameters}. Let $Q^{\alpha\beta\gamma}(\tau) := \frac{d}{d\tau} a^{\alpha\beta\gamma}(\tau) - 2 \, (-\tau)^{-1} \, a^{\alpha\beta\gamma}(\tau)$. Then 
\[\sup_{\tau \leq \tau_*} (-\tau) \int_{\tau-1}^\tau |Q^{\alpha\beta\gamma}(\tau')| \, d\tau' \leq \delta \, \sup_{\tau \leq \tau_*} \bigg ( \int_{\tau-1}^\tau a^{\alpha\beta\gamma}(\tau')^2 \, d\tau' \bigg )^{\frac{1}{2}}.\] 
\end{proposition}

The proof of Proposition \ref{ode.for.a} will be given in Section \ref{analysis.of.neutral.mode}. \\

We can now finish the proof of Theorem \ref{uniqueness.theorem}, exactly as in \cite{BDS} in dimension three. For the convenience of the reader, we include a copy of the proof here. 

Using the ODE $\frac{d}{d\tau} a^{\alpha\beta\gamma}(\tau) = 2 \, (-\tau)^{-1} \, a^{\alpha\beta\gamma}(\tau) + Q^{\alpha\beta\gamma}(\tau)$ together with the fact that $a^{\alpha\beta\gamma}(\tau_*) = 0$, we obtain 
\[(-\tau)^2 \, a^{\alpha\beta\gamma}(\tau) = -\int_\tau^{\tau_*} (-\tau')^2 \, Q^{\alpha\beta\gamma}(\tau') \, d\tau'.\] 
This implies 
\begin{align*} 
(-\tau) \, |a^{\alpha\beta\gamma}(\tau)| 
&\leq \int_\tau^{\tau_*} (-\tau') \, |Q^{\alpha\beta\gamma}(\tau')| \, d\tau' \\ 
&\leq \sum_{j=0}^{[\tau_*-\tau]} \int_{\tau_*-j-1}^{\tau_*-j} (-\tau') \, |Q^{\alpha\beta\gamma}(\tau')| \, d\tau' \\ 
&\leq (-\tau) \, \max_{0 \leq j \leq [\tau_*-\tau]} \int_{\tau_*-j-1}^{\tau_*-j} (-\tau') \, |Q^{\alpha\beta\gamma}(\tau')| \, d\tau'.
\end{align*}
We now divide by $-\tau$, and take the supremum over all $\tau \leq \tau_*$. This implies 
\[\sup_{\tau \leq \tau_*} |a^{\alpha\beta\gamma}(\tau)| \leq \sup_{\tau \leq \tau_*} \int_{\tau-1}^\tau (-\tau') \, |Q^{\alpha\beta\gamma}(\tau')| \, d\tau'.\]
On the other hand, Proposition \ref{ode.for.a} gives the following estimate for $Q^{\alpha\beta\gamma}$: 
\[\sup_{\tau \leq \tau_*} (-\tau) \int_{\tau-1}^\tau |Q^{\alpha\beta\gamma}(\tau')| \, d\tau' \leq \delta \, \sup_{\tau \leq \tau_*} |a^{\alpha\beta\gamma}(\tau)|.\] 
Hence, if we choose $\delta$ sufficiently small, and $-\tau_*$ sufficiently large (depending on $\delta$), then $\sup_{\tau \leq \tau_*} |a^{\alpha\beta\gamma}(\tau)| = 0$. Thus, $a^{\alpha\beta\gamma}(\tau) = 0$ for all $\tau \leq \tau_*$. Proposition \ref{neutral.mode.dominates} then implies $\hat{H}_{\mathcal{C}}^{\alpha\beta\gamma}(\xi,\tau) = 0$ for all $\tau \leq \tau_*$. Putting these facts together, we obtain $H_{\mathcal{C}}^{\alpha\beta\gamma}(\xi,\tau) = 0$ for all $\tau \leq \tau_*$. From this, we deduce that $W_+^{\beta\gamma}(\rho,\tau) = 0$ for $\rho \in [\theta,2\theta]$ and $\tau \leq \tau_*$. Proposition \ref{estimate.for.difference.in.tip.region} yields $W_+^{\beta\gamma}(\rho,\tau) = 0$ for $\rho \in [0,2\theta]$ and $\tau \leq \tau_*$. Thus, we conclude that $F_1(z,t) = F_2^{\alpha\beta\gamma}(z,t)$ for all $t \leq t_* = -e^{-\tau_*}$. In other words, the two ancient solutions coincide for $t \leq t_*$. \\

\section{Energy estimates in the tip region and proof of Proposition \ref{estimate.for.difference.in.tip.region}}

\label{difference.tip.region}

In this section, we give the proof of Proposition \ref{estimate.for.difference.in.tip.region}. Let $\omega_T$ denote a nonnegative smooth cutoff function satisfying $\omega_T(\rho)=1$ for $\rho \leq \theta$ and $\omega_T(\rho) = 0$ for $\rho \geq 2\theta$. We define 
\[W_{T+}^{\beta\gamma}(\rho,\tau) := \omega_T(\rho) \, W^{\beta\gamma}_+(\rho,\tau).\]  
To simplify the notation, we will write $W_+$ and $W_{T+}$ instead of $W_+^{\beta\gamma}$ and $W_{T+}^{\beta\gamma}$.   \\

\begin{proposition}
\label{pde.for.W}
The function $W_+(\rho,\tau)$ satisfies the equation 
\begin{align}\label{eqn-W122} 
V_{1+}^{-2} \, \Big ( \frac{\partial W_+}{\partial \tau} + \frac{\rho}{2} \, \frac{\partial W_+}{\partial \rho} \Big ) 
&= \frac{\partial^2 W_+}{\partial \rho^2} + \frac{\partial}{\partial \rho} \Big ( \rho^{-1} \, (n-2) \, (V_{1+}^{-2} - 1) \, W_+ \Big ) \\ 
&+ (n-3) \rho^{-1}  \frac{\partial W_+}{\partial\rho}-2(n-2) \, \rho^{-2} \, W_+ + V_{1+}^{-2} \, \mathcal{B}_+ \, W_+. 
\end{align}
where 
\begin{align*} 
\mathcal{B}_+ 
&:= (n-2)\,\rho^{-2} \, \big ( 1 - V_{1+} \, (V_{2+}^{\beta\gamma})^{-1} \big ) \\ 
& + (n-2)\, \rho^{-1} \, \Big ( 2 \, V_{1+}^{-1} \, \frac{\partial V_{1+}}{\partial \rho} - (V_{2+}^{\beta\gamma})^{-2} \, (V_{1+}+V_{2+}^{\beta\gamma}) \, \frac{\partial V_{2+}^{\beta\gamma}}{\partial \rho} \Big ) \\ 
& + (V_{2+}^{\beta\gamma})^{-2} \, (V_{1+}+V_{2+}^{\beta\gamma}) \, \Big ( \frac{\partial V_{2+}^{\beta\gamma}}{\partial \tau} + \frac{\rho}{2} \, \frac{\partial V_{2+}^{\beta\gamma}}{\partial \rho} \Big ).
\end{align*} 
\end{proposition}

\textbf{Proof.} 
The functions $U_{1+}(r,t)$, $U_{1-}(r,t)$, $U_{2+}^{\beta\gamma}(r,t)$, and $U_{2-}^{\beta\gamma}(r,t)$ all satisfy the same PDE: 
\bee
\begin{split}
U^{-1} \, \frac{\partial U}{\partial t} = \frac{\partial^2 U}{\partial r^2} - \frac{1}{2} U^{-1}  \Big ( \frac{\partial U}{\partial r} \Big )^2 + \frac{(n-2)}{r^2} (U^{-1}-1) \Big (r \frac{\partial U}{\partial r}  + 2U \Big ) + \frac{(n-3)}{r}  \frac{\partial U}{\partial r} .
\end{split}
\eee
Consequently, the functions $V_{1+}(\rho,\tau)$, $V_{1-}(\rho,\tau)$, $V_{2+}^{\beta\gamma}(\rho,\tau)$, and $V_{2-}^{\beta\gamma}(\rho,\tau)$ satisfy the following PDE: 
\[V^{-2} \, \Big ( \frac{\partial V}{\partial \tau} + \frac{\rho}{2} \, \frac{\partial V}{\partial \rho} \Big ) = \frac{\partial^2 V}{\partial \rho^2} + \frac{(n-2)}{\rho^{2}}  \,(V^{-2}-1) \, \Big ( \rho \, \frac{\partial V}{\partial \rho}+V \Big ) + \frac{(n-3)}{\rho}  \frac{\partial V}{\partial \rho}.\] 
The assertion now follows by a  straightforward calculation.

\begin{proposition}
\label{divergence.identity}
The function $W_{T+}(\rho,\tau)$ satisfies 
\begin{align*} 
&\frac{1}{2} \, \frac{\partial}{\partial \tau} \big ( V_{1+}^{-2} W_{T+}^2 e^{\mu_+} \big ) - 
\frac{\partial}{\partial \rho} \Big [ \Big ( \frac{\partial W_{T+}}{\partial \rho} +  (n-2) \rho^{-1} \big ( V_{1+}^{-2} - 1 \big ) \,  W_{T+} \Big ) W_{T+} \, e^{\mu_+} \Big ] \\ 
&+ \frac{\partial}{\partial \rho} \big ( W_+^2 \, \omega_T' \, \omega_T \, e^{\mu_+} \big ) -\frac{1}{2} \frac{\partial}{\partial \rho}\big((n-3) \rho^{-1} W_{T+}^2 e^{\mu_+}\big)\\ 
&\leq -\frac{1}{2} \, \Big ( \frac{\partial W_{T+}}{\partial \rho} + \frac{\partial \mu_+}{\partial \rho} \, W_{T+} \Big )^2 \, e^{\mu_+} -\frac{1}{2}(3n-5) \,\rho^{-2} \, W_{T+}^2 \, e^{\mu_+} \\ 
&+ V_{1+}^{-2} \, \Big ( \frac{1}{2} \, \frac{\partial \mu_+}{\partial \tau} - V_{1+}^{-1} \, \frac{\partial V_{1+}}{\partial \tau} + \frac{\rho}{2} \, \frac{\partial \mu_+}{\partial \rho} + \mathcal{B}_+ \Big ) \, W_{T+}^2 \, e^{\mu_+} \\ 
&+ \frac{1}{2} \, \Big ( \frac{\partial \mu_+}{\partial \rho} - (n-2) \rho^{-1}  \big ( V_{1+}^{-2} - 1 \big )  - \frac{\rho}{2} \, V_{1+}^{-2} \Big )^2 \, W_{T+}^2 \, e^{\mu_+} \\ 
&+ \Big ( \frac{\partial \mu_+}{\partial \rho} -  (n-2) \rho^{-1}  \big ( V_{1+}^{-2} - 1  \big )   + \frac{\rho}{2} \, V_{1+}^{-2}   - (n-3)\, \rho^{-1}\Big )  W_+^2  \omega_T'  \, \omega_T \, e^{\mu_+}  \\&+ (\omega_T')^2  \, W_+^2 e^{\mu_+}  - \frac{1}{2}(n-3)\rho^{-1} W_{T+}^2 \frac{\partial \mu_+}{\partial \rho} e^{\mu_+}
\end{align*} 
\end{proposition}

\textbf{Proof.} Using Proposition \ref{pde.for.W}, we obtain  
\begin{align*} 
&V_{1+}^{-2} \, \Big ( \frac{\partial W_{T+}}{\partial \tau} + \frac{\rho}{2} \, \frac{\partial W_{T+}}{\partial \rho} \Big ) \\ 
&= \frac{\partial^2 W_{T+}}{\partial \rho^2} +  \frac{\partial}{\partial \rho} \Big (  \rho^{-1} (n-2) (V_{1+}^{-2} - 1) \, W_{T+} \Big ) - 2(n-2) \rho^{-2} \, W_{T+} + V_{1+}^{-2} \, \mathcal{B}_+ \, W_{T+} \\ 
&+ \Big ( -2 \, \omega_T' \, \frac{\partial W_+}{\partial \rho} - \omega_T'' \, W_+ - (n-2)  \rho^{-1} \, (V_{1+}^{-2} - 1) \, W_+ \omega_T'
+ \frac{\rho}{2} \, \omega_T' \, V_{1+}^{-2} \, W_+ \Big )\\
&+ (n-3) \rho^{-1}  \Big ( \frac{\partial W_{T+}}{\partial\rho} -  W_+ \, \omega_T' \Big ).
\end{align*}
We next bring in the weight $\mu_+(\rho,\tau)$. A straightforward calculation gives 
\begin{align*} 
&\frac{1}{2} \, \frac{\partial}{\partial \tau} \big ( V_{1+}^{-2} \, W_{T+}^2 \, e^{\mu_+} \big ) - \frac{\partial}{\partial \rho} \Big [ \Big ( \frac{\partial W_{T+}}{\partial \rho} + \rho^{-1} (n-2)  (V_{1+}^{-2} - 1) \, W_{T+} \Big ) \, W_{T+} \, e^{\mu_+} \Big ] \\ 
&+ \frac{\partial}{\partial \rho} \big ( W_+^2 \, \omega_T' \, \omega_T \, e^{\mu_+} \big ) -\frac{1}{2} \frac{\partial}{\partial \rho}\big((n-3) \rho^{-1} W_{T+}^2 e^{\mu_+}\big)\\ 
&= -\Big ( \frac{\partial W_{T+}}{\partial \rho} + \frac{\partial \mu_+}{\partial \rho} \, W_{T+} \Big )^2 \, e^{\mu_+} - \frac{1}{2} (3n-5) \rho^{-2} \, W_{T+}^2 \, e^{\mu_+} \\ 
&+ V_{1+}^{-2} \, \Big ( \frac{1}{2} \, \frac{\partial \mu_+}{\partial \tau} - V_{1+}^{-1} \, \frac{\partial V_{1+}}{\partial \tau} + \frac{\rho}{2} \, \frac{\partial \mu_+}{\partial \rho} + \mathcal{B}_+ \Big ) \, W_{T+}^2 \, e^{\mu_+} \\ 
&+ \Big ( \frac{\partial \mu_+}{\partial \rho} - (n-2)\rho^{-1}  (V_{1+}^{-2} - 1) - \frac{\rho}{2} \, V_{1+}^{-2} \Big ) \Big ( \frac{\partial W_{T+}}{\partial \rho} + \frac{\partial \mu_+}{\partial \rho} \, W_{T+} \Big ) \, W_{T+} \, e^{\mu_+} \\ 
&+ \Big ( \frac{\partial \mu_+}{\partial \rho} -  (n-2)   \rho^{-1}(V_{1+}^{-2} - 1) +
 \frac{\rho}{2} \, V_{1+}^{-2} \Big ) \, W_+^2 \, \omega_T' \, \omega_T \, e^{\mu_+} + (\omega_T')^2 \, W_+^2 \, e^{\mu_+}\\
&-  (n-3)\, \rho^{-1} W_+^2 \, \omega_T   \, \omega_T' e^{\mu_+}- \frac{1}{2}(n-3) \rho^{-1}W_{T+}^2 \frac{\partial \mu_+}{\partial \rho} e^{\mu_+}.
\end{align*} 
The assertion follows now from Young's inequality and combining terms. \\

\begin{corollary}
\label{divergence.identity.2}
Fix a small number $\eta > 0$. Then we can find a small number $\theta \in (0,\eta)$ and a small number $\varepsilon \in (0,\eta)$ (both depending on $\eta$) with the following property. If $-\tau_*$ sufficiently large (depending on $\eta$ and $\theta$) and the triplet $(\alpha,\beta,\gamma)$ is $\varepsilon$-admissible with respect to time $t_* = -e^{-\tau_*}$, then we have 
\begin{align*} 
&\frac{1}{2} \, \frac{\partial}{\partial \tau} \big ( V_{1+}^{-2} \, W_{T+}^2 \, e^{\mu_+} \big ) - \frac{\partial}{\partial \rho} \Big [ \Big ( \frac{\partial W_{T+}}{\partial \rho} + \rho^{-1} (n-2) \big (  V_{1+}^{-2} - 1\big )  W_{T+} \Big ) W_{T+}  e^{\mu_+} \Big ] \\ 
&+ \frac{\partial}{\partial \rho} \big ( W_+^2 \, \omega_T' \, \omega_T \, e^{\mu_+} \big )-\frac{1}{2} \frac{\partial}{\partial \rho}\big((n-3) \rho^{-1} W_{T+}^2 e^{\mu_+}\big)  \\ 
&\leq -\frac{1}{2} \, \Big ( \frac{\partial W_{T+}}{\partial \rho} + \frac{\partial \mu_+}{\partial \rho} \, W_{T+} \Big )^2 \, e^{\mu_+} - \frac{1}{2}(3n-5) \,
\rho^{-2}  W_{T+}^2 \, e^{\mu_+} \\ 
&+ \eta \, \rho^{-2} \, V_{1+}^{-4} \, W_{T+}^2 \, e^{\mu_+} + \eta \, \rho^{-2} \, V_{1+}^{-2} \, W_+^2 \, e^{\mu_+} \, 1_{\{\theta \leq \rho \leq 2\theta\}}.
\end{align*} 
for $\rho \leq 2\theta$ and $\tau \leq \tau_*$.
\end{corollary}

\textbf{Proof.}
By Proposition \ref{comparison.of.V.with.bryant.soliton.profile}, Proposition \ref{higher.derivatives.of.V.in.collar.region}, and Proposition \ref{comparison.of.modified.V.with.bryant.soliton.profile}, we can choose $\theta \in (0,\eta)$ (depending on $\eta$) sufficiently small and $-\tau_*$ sufficiently large (depending on $\eta$ and $\theta$) such that 
\[|\mathcal{B}_+| \leq \eta \, \rho^{-2} \, V_{1+}^{-2}\] 
for $\rho \leq 2\theta$ and $\tau \leq \tau_*$. By Corollary \ref{derivative.of.V.wrt.tau}, Lemma \ref{derivative.of.mu.wrt.rho}, and Lemma \ref{derivative.of.mu.wrt.tau}, we can choose $\theta \in (0,\eta)$ sufficiently small (depending on $\eta$) and $-\tau_*$ sufficiently large (depending on $\eta$ and $\theta$) such that 
\begin{align*}
&\Big | \frac{1}{2} \, \frac{\partial \mu_+}{\partial \tau} - V_{1+}^{-1} \, \frac{\partial V_{1+}}{\partial \tau} + \frac{\rho}{2} \, \frac{\partial \mu_+}{\partial \rho} \Big | \leq \eta \, \rho^{-2} \, V_{1+}^{-2},  \\ 
&\Big | \frac{\partial \mu_+}{\partial \rho} - (n-2)\,\rho^{-1} \, (V_{1+}^{-2} - 1) - \frac{\rho}{2} \, V_{1+}^{-2} \Big | \leq \eta \, \rho^{-1} \, V_{1+}^{-2}, \\ 
&\Big | \frac{\partial \mu_+}{\partial \rho} - (n-2)\,\rho^{-1} \, (V_{1+}^{-2} - 1) + \frac{\rho}{2} \, V_{1+}^{-2} \Big | \leq \eta \, \rho^{-1} \, V_{1+}^{-2} 
\end{align*} 
for $\rho \leq 2\theta$ and $\tau \leq \tau_*$. Note also that for any $\eta > 0$ there exists a $\theta \in (0,\eta)$, 
so that $|(n-3) \rho^{-1} \omega_T' \omega_T | \le \eta\,  \rho^{-2} V_{1+}^{-2} \,  1_{\{\theta\le \rho \le 2\theta\}}$. Finally, recall by Lemma \ref{derivative.of.mu.wrt.rho} $\frac{\partial \mu_+}{\partial \rho} \geq 0$, so 
\[
- \frac{1}{2}(n-3) \rho^{-1}W_{T+}^2 \frac{\partial \mu_+}{\partial \rho} e^{\mu_+} \leq 0.
\]
Hence, the assertion follows from Proposition \ref{divergence.identity}. \\

We now finalize our choice of $\theta$.

\begin{proposition}
\label{integral.estimate}
We can find sufficiently small numbers $\theta>0$, $\lambda>0$, and $\varepsilon>0$ with the following property. If $-\tau_*$ is sufficiently large (depending on $\theta$) and the triplet $(\alpha,\beta,\gamma)$ is $\varepsilon$-admissible with respect to time $t_* = -e^{-\tau_*}$, then 
\begin{align*} 
\frac{1}{2} \, \frac{d}{d\tau} \bigg ( \int_0^{2\theta} V_{1+}^{-2} \, W_{T+}^2 \, e^{\mu_+} \, d\rho \bigg ) 
&\leq -\lambda \, (-\tau) \int_0^{2\theta} V_{1+}^{-2} \, W_{T+}^2 \, e^{\mu_+} \, d\rho \\ 
&+ \int_\theta^{2\theta} \rho^{-2} \, V_{1+}^{-2} \, W_+^2 \, e^{\mu_+} \, d\rho 
\end{align*}
for $\tau \leq \tau_*$.
\end{proposition}

\textbf{Proof.} 
Fix a small number $\eta > 0$. In the following, we choose $\theta$ and $\varepsilon$ sufficiently small (depending on $\eta$), and we choose $-\tau_*$ sufficiently large (depending on $\eta$ and $\theta$). Recall $W_{T+}(2\theta, \tau) = 0$ and, for each $\tau$, $W_+(\rho, \tau) = O(\rho^2)$ as $\rho \to 0$, hence $W_{T+}(\rho, \tau) = O(\rho^2)$ as $\rho \to 0$.  Moreover, we can see from the definition of $\mu_+$ that $\mu_+(\rho, \tau)$ is bounded as $\rho \to 0$. In particular, if we integrate the differential inequality of  Corollary \ref{divergence.identity.2}, the divergence terms vanish. Using Corollary \ref{divergence.identity.2}, we obtain 
\begin{align*} 
&\frac{1}{2} \, \frac{d}{d\tau} \bigg ( \int_0^{2\theta} V_{1+}^{-2} \, W_{T+}^2 \, e^{\mu_+} \, d\rho \bigg ) \\ 
&\leq -\frac{1}{2} \int_0^{2\theta} \Big ( \frac{\partial W_{T+}}{\partial \rho} + \frac{\partial \mu_+}{\partial \rho} \, W_{T+} \Big )^2 \, e^{\mu_+} \, d\rho -\frac{1}{2}(3n-5) \int_0^{2\theta} \rho^{-2} \, W_{T+}^2 \, e^{\mu_+} \, d\rho\\
&+ \eta \int_0^{2\theta} \rho^{-2} \, V_{1+}^{-4}\, W_{T+}^2 \, e^{\mu_+} \, d\rho + \eta \int_\theta^{2\theta} \rho^{-2} \, V_{1+}^{-2} \, W_+^2 \, e^{\mu_+} \, d\rho,
\end{align*} 
for $\tau \leq \tau_*$.

We will next estimate the terms on the right hand side of the above inequality to deduce the statement of the Proposition. 
First, applying  Proposition \ref{Poincare.inequality} to the function $f := e^{\mu_+} W_{T+}$ gives  
\begin{align*} 
0 &\leq 8 \int_0^{2\theta} \Big ( \frac{\partial W_{T+}}{\partial \rho} + \frac{\partial \mu_+}{\partial \rho} \, W_{T+} \Big )^2  \, e^{\mu_+} \, d\rho \\ 
&+ K_* \int_0^{2\theta} \rho^{-2} \, W_{T+}^2 \, e^{\mu_+} \, d\rho - \int_0^{2\theta} \Big ( \frac{\partial \mu_+}{\partial \rho} \Big )^2 \, W_{T+}^2 \, e^{\mu_+} \, d\rho 
\end{align*}
for $\tau \leq \tau_*$. Using Lemma \ref{derivative.of.mu.wrt.rho}, we obtain $(\frac{\partial \mu_+}{\partial \rho})^2 \geq \frac{1}{4} \, \rho^{-2} \, (V_{1+}^{-2}-1)^2$ for $\rho \leq 2\theta$, hence 
\begin{align*} 
0 &\leq 128\eta \int_0^{2\theta} \Big ( \frac{\partial W_{T+}}{\partial \rho} + \frac{\partial \mu_+}{\partial \rho} \, W_{T+} \Big )^2  \, e^{\mu_+} \, d\rho \\ 
&+ 16\eta K_* \int_0^{2\theta} \rho^{-2} \, W_{T+}^2 \, e^{\mu_+} \, d\rho - 4\eta \int_0^{2\theta} \rho^{-2} \, (V_{1+}^{-2}-1)^2 \, W_{T+}^2 \, e^{\mu_+} \, d\rho 
\end{align*}
for $\tau \leq \tau_*$. 

Adding the two inequalities above, we obtain
\begin{align*} 
&\frac{1}{2} \, \frac{d}{d\tau} \bigg ( \int_0^{2\theta} V_{1+}^{-2} \, W_{T+}^2 \, e^{\mu_+} \, d\rho \bigg ) \\ 
&\leq -\Big ( \frac{1}{2}-128\eta \Big ) \int_0^{2\theta} \Big ( \frac{\partial W_{T+}}{\partial \rho} + \frac{\partial \mu_+}{\partial \rho} \, W_{T+} \Big )^2 \, e^{\mu_+} \, d\rho \\ 
&- \Big(\frac{1}{2}(3n-5)-4\eta -16\eta K_*\Big) \int_0^{2\theta} \rho^{-2} \, W_{T+}^2 \, e^{\mu_+} \, d\rho \\ 
&- \eta \int_0^{2\theta} \rho^{-2} \, (4(V_{1+}^{-2} - 1)^2 + 4 -V_{1+}^{-4}) \, W_{T+}^2 \, e^{\mu_+} \, d\rho \\
& + \eta \int_\theta^{2\theta} \rho^{-2} \, V_{1+}^{-2} \, W_+^2 \, e^{\mu_+} \, d\rho 
\end{align*} 
for $\tau \leq \tau_*$. We now choose $\eta>0$ sufficiently small so that $\frac{1}{2}-128\eta > 0$ and $\frac{3n-5}{2} - 4 \eta-16\eta K_* > 0 $. (Here, it is crucial that the constant $K_*$ in the weighted Poincar\'e inequality does not depend on $\theta$.) This ensures that the first two terms on the right hand side of the last estimate have a favorable sign. To estimate the third term on the right hand side, we observe that $\rho^{-2}[4(V_{1+}^{-2} - 1)^2 + 4 - V_{1+}^{-4}] \geq \rho^{-2} V_{1+}^{-4}$. 
 Finally, in view of Proposition \ref{comparison.of.V.with.bryant.soliton.profile}, we can bound $\rho^{-2} V_{1+}^{-4}$ from below by a small positive multiple of $(-\tau) \, V_{1+}^{-2}$. This completes the proof of Proposition \ref{integral.estimate}. \\

We now complete the proof of Proposition \ref{estimate.for.difference.in.tip.region}. Let $\theta$, $\lambda$, and $\varepsilon$ be chosen as in Proposition \ref{integral.estimate}. Let 
\[I(\tau) := \int_{\tau-1}^\tau \int_0^{2\theta} V_{1+}^{-2} \, W_{T+}^2 \, e^{\mu_+}\] 
and 
\[J(\tau) := \int_{\tau-1}^\tau \int_\theta^{2\theta} V_{1+}^{-2} \, W_+^2 \, e^{\mu_+}.\] 
If we choose $-\tau_*$ sufficiently large, then Proposition \ref{integral.estimate} gives 
\[\frac{1}{2} \, I'(\tau) + \lambda \, (-\tau) \, I(\tau) \leq \theta^{-2} \, J(\tau),\] 
hence 
\[\frac{d}{d\tau} (e^{-\lambda \tau^2} \, I(\tau)) \leq 2\theta^{-2} \, e^{-\lambda \tau^2} \, J(\tau)\] 
for $\tau \leq \tau_*$. Clearly, $\lim_{\tau \to -\infty} e^{-\lambda \tau^2} \, I(\tau) = 0$. Consequently, 
\begin{align*} 
e^{-\lambda \tau^2} \, I(\tau) 
&\leq 2\theta^{-2} \int_{-\infty}^\tau e^{-\lambda {\tau'}^2} \, J(\tau') \, d\tau' \\ 
&\leq 2\theta^{-2} \, \Big ( \sup_{\tau' \leq \tau} (-\tau')^{-1} \, J(\tau') \Big ) \, \int_{-\infty}^\tau e^{-\lambda {\tau'}^2} \, (-\tau') \, d\tau' \\ 
&\leq \theta^{-2} \lambda^{-1} \, e^{-\lambda \tau^2} \, \sup_{\tau' \leq \tau} (-\tau')^{-1} \, J(\tau') 
\end{align*} 
for $\tau \leq \tau_*$. This finally gives 
\begin{align*} 
(-\tau)^{-\frac{1}{2}} \, I(\tau) 
&\leq \theta^{-2} \lambda^{-1} \, (-\tau)^{-\frac{1}{2}} \sup_{\tau' \leq \tau} (-\tau')^{-1} \, J(\tau') \\ 
&\leq \theta^{-2} \lambda^{-1} \, (-\tau)^{-1} \sup_{\tau' \leq \tau} (-\tau')^{-\frac{1}{2}} \, J(\tau')  
\end{align*}
for $\tau \leq \tau_*$. Taking the supremum over $\tau \leq \tau_*$ gives 
\[\sup_{\tau \leq \tau_*} (-\tau)^{-\frac{1}{2}} \, I(\tau) \leq \theta^{-2} \lambda^{-1} \, (-\tau_*)^{-1} \, \sup_{\tau \leq \tau_*} (-\tau)^{-\frac{1}{2}} \, J(\tau).\] 
From this, the conclusion of Proposition \ref{estimate.for.difference.in.tip.region} follows immediately. \\

\section{Energy estimates in the cylindrical region and proof of Proposition \ref{estimate.for.difference.in.cylindrical.region}}

\label{difference.cylindrical.region}

In this section, we give the proof of Proposition \ref{estimate.for.difference.in.cylindrical.region}. Throughout this section, we assume that $\theta$ is chosen as in Proposition \ref{estimate.for.difference.in.tip.region}. To simplify the notation, we will write $H$, $H_{\mathcal{C}}$, $\hat{H}_{\mathcal{C}}$, and $a$ instead of $H^{\alpha\beta\gamma}$, $H_{\mathcal{C}}^{\alpha\beta\gamma}$, $\hat{H}_{\mathcal{C}}^{\alpha\beta\gamma}$, and $a^{\alpha\beta\gamma}$. 

Our goal is to study the evolution equation satisfied by the function $H$. The linearized operator 
\[\mathcal{L} f := f_{\xi\xi} - \frac{1}{2} \, \xi \, f_\xi + f\] 
is the same as in \cite{Angenent-Daskalopoulos-Sesum2}, and hence the linear theory from \cite{Angenent-Daskalopoulos-Sesum2} carries over to the Ricci flow case as well. In order for this article to be self-contained, we will state the results from \cite{Angenent-Daskalopoulos-Sesum2} that we will use later, but for the proofs of the same we refer the reader to \cite{Angenent-Daskalopoulos-Sesum2}. 

As in \cite{Angenent-Daskalopoulos-Sesum2}, we consider the Hilbert space $\mathcal{H} = L^2(\mathbb{R},e^{-\frac{\xi^2}{4}} \, d\xi)$. The norm on $\mathcal{H}$ is given by 
\[\|f\|_{\mathcal{H}}^2 := \int_{\mathbb{R}} e^{-\frac{\xi^2}{4}} \, f(\xi)^2 \, d\xi.\] 
Moreover, we denote by $\mathcal{D} \subset \mathcal{H}$ the Hilbert space of all functions $f$ such that $f \in \mathcal{H}$ and $f' \in \mathcal{H}$. The norm on $\mathcal{D}$ is given by 
\[\|f\|_{\mathcal{D}}^2 := \int_{\mathbb{R}} e^{-\frac{\xi^2}{4}} \, (f'(\xi)^2+f(\xi)^2) \, d\xi.\] 
Let $\mathcal{D}^*$ denote the dual space of $\mathcal{D}$. Clearly, the dual space $\mathcal{H}^*$ is a subspace of $\mathcal{D}^*$. After identifying $\mathcal{H}^*$ with $\mathcal{H}$ in the standard way, we can view $\mathcal{H}$ as a subspace of $\mathcal{D}^*$. The restriction of $\|\cdot\|_{\mathcal{D}^*}$ to $\mathcal{H}$ is given by 
\[\|f\|_{\mathcal{D}^*} := \sup \bigg \{ \int_{\mathbb{R}} e^{-\frac{\xi^2}{4}} \, f(\xi) \, g(\xi) \, d\xi: \|g\|_{\mathcal{D}} \leq 1 \bigg \}.\] 
For later reference, we collect some basic facts from \cite{Angenent-Daskalopoulos-Sesum2}. \\

\begin{proposition}
\label{boundedness.of.operators}
The following statements hold:
\begin{itemize}
\item[(i)] The operators $f \mapsto \xi \, f$, $f \mapsto f'$, $f \mapsto -f' + \frac{1}{2} \, \xi \, f$ are bounded from $\mathcal{D}$ to $\mathcal{H}$. 
\item[(ii)] The operators $f \mapsto \xi \, f$, $f \mapsto f'$, $f \mapsto -f' + \frac{1}{2} \, \xi \, f$ are bounded from $\mathcal{H}$ to $\mathcal{D}^*$. 
\item[(iii)] The operators $f \mapsto \xi^2 \, f$, $f \mapsto \xi \, f'$, $f \mapsto f''$ are bounded from $\mathcal{D}$ to $\mathcal{D}^*$.
\item[(iv)] The operator $f \mapsto \int_0^\xi f$ is bounded from $\mathcal{H}$ to $\mathcal{D}$.
\end{itemize}
\end{proposition}

\textbf{Proof.} 
Statements (i), (ii), and (iii) were proved in \cite{Angenent-Daskalopoulos-Sesum2}. To prove statement (iv), let us consider a function $f \in \mathcal{H}$, and let $g(\xi) := \int_0^\xi f(\xi') \, d\xi'$. Then $g(\xi)^2 \leq \xi \int_0^\xi f(\xi')^2 \, d\xi'$ for $\xi \geq 0$. Using Fubini's theorem, we obtain 
\begin{align*} 
\int_0^\infty e^{-\frac{\xi^2}{4}} \, g(\xi)^2 \, d\xi 
&\leq \int_0^\infty e^{-\frac{\xi^2}{4}} \, \xi \, \bigg ( \int_0^\xi f(\xi')^2 \, d\xi' \bigg ) \, d\xi \\ 
&= \int_0^\infty \bigg ( \int_{\xi'}^\infty e^{-\frac{\xi^2}{4}} \, \xi \, d\xi \bigg ) \, f(\xi')^2 \, d\xi' \\ 
&= 2 \int_0^\infty e^{-\frac{\xi'^2}{4}} \, f(\xi')^2 \, d\xi'. 
\end{align*}
An analogous argument gives $\int_{-\infty}^0 e^{-\frac{\xi^2}{4}} \, g(\xi)^2 \, d\xi \leq 2 \int_{-\infty}^0 e^{-\frac{\xi'^2}{4}} \, f(\xi')^2 \, d\xi'$. Therefore, $\|g\|_{\mathcal{H}} \leq C \, \|f\|_{\mathcal{H}}$. Since $g'=f$, it follows that $\|g\|_{\mathcal{D}} \leq C \, \|f\|_{\mathcal{H}}$, as claimed. \\

For a time-dependent function $f$, we introduce the following norms:
\begin{align*} 
&\|f\|_{\mathcal{H},\infty,\tau_*}^2 := \sup_{\tau \leq \tau_*} \int_{\tau-1}^\tau \|f(\cdot,\tau')\|_{\mathcal{H}}^2 \, d\tau', \\ 
&\|f\|_{\mathcal{D},\infty,\tau_*}^2 := \sup_{\tau \leq \tau_*} \int_{\tau-1}^\tau \|f(\cdot,\tau')\|_{\mathcal{D}}^2 \, d\tau', \\ 
&\|f\|_{\mathcal{D}^*,\infty,\tau_*}^2 := \sup_{\tau \leq \tau_*} \int_{\tau-1}^\tau \|f(\cdot,\tau')\|_{\mathcal{D}^*}^2 \, d\tau'.
\end{align*}
The following energy estimate was proved in \cite{Angenent-Daskalopoulos-Sesum2}: \\

\begin{proposition}
\label{linear.estimate}
Let $g: (-\infty,\tau_*] \to \mathcal{D}^*$ be a bounded function. Let $f: (-\infty,\tau_*] \to \mathcal{D}$ be a bounded function which satisfies the linear equation 
\[\frac{\partial}{\partial \tau} f(\tau) - \mathcal{L} f(\tau) = g(\tau).\] 
Then the function $\hat{f} := P_+ f + P_- f$ satisfies the estimate 
\[\sup_{\tau \leq \tau_*} \|\hat{f}(\tau)\|_{\mathcal{H}} + \Lambda^{-1} \, \|\hat{f}\|_{\mathcal{D},\infty,\tau_*} \leq \|P_+ f(\tau_*)\|_{\mathcal{H}} + \Lambda \, \|g\|_{\mathcal{D}^*,\infty,\tau_*},\] 
where $\Lambda$ is a universal constant.
\end{proposition}

\textbf{Proof.} 
See \cite{Angenent-Daskalopoulos-Sesum2}, Lemma 6.6. \\

We now continue with the proof of Proposition \ref{estimate.for.difference.in.cylindrical.region}. The functions $G_1(\xi,\tau)$ and $G_2^{\alpha\beta\gamma}(\xi,\tau)$ satisfy the equation 
\begin{align*} 
G_\tau(\xi,\tau) 
&= G_{\xi\xi}(\xi,\tau) - \frac{1}{2} \, \xi \, G_\xi(\xi,\tau) \\
&+ \frac{1}{2} \, (\sqrt{2(n-2)}+G(\xi,\tau)) - (n-2)\,(\sqrt{2(n-2)}+G(\xi,\tau))^{-1} \\ 
&-(\sqrt{2(n-2)}+G(\xi,\tau))^{-1} \, G_\xi(\xi,\tau)^2 \\ 
&+ (n-1) \, G_\xi(\xi,\tau) \, \bigg [ \frac{G_\xi(0,\tau)}{\sqrt{2(n-2)}+G(0,\tau)} - \int_0^\xi \frac{G_\xi(\xi',\tau)^2}{(\sqrt{2(n-2)}+G(\xi',\tau))^2} \, d\xi' \bigg ].
\end{align*} 
Note that the two terms on the second line above can be written
\[
\frac{1}{2} \, (\sqrt{2(n-2)}+G) - (n-2)\,(\sqrt{2(n-2)}+G)^{-1} = G - \frac{1}{2}(\sqrt{2(n-2)} + G)^{-1}G^2. 
\]
Consequently, the difference $H(\xi,\tau) = G_1(\xi,\tau) - G_2^{\alpha\beta\gamma}(\xi,\tau)$ satisfies 
\[H_\tau(\xi,\tau) = H_{\xi\xi}(\xi,\tau) - \frac{1}{2} \, \xi \, H_\xi(\xi,\tau) + H(\xi,\tau) + \sum_{k=1}^6 E_k(\xi,\tau),\] 
where 
\begin{align*} 
E_1(\xi,\tau) &= \Big [ (n-2)\,(\sqrt{2(n-2)}+G_1(\xi,\tau))^{-1} (\sqrt{2(n-2)}+G_2^{\alpha\beta\gamma}(\xi,\tau))^{-1} - \frac{1}{2} \Big ] \, H(\xi,\tau) \\ 
E_2(\xi,\tau) &= (\sqrt{2(n-2)}+G_1(\xi,\tau))^{-1} (\sqrt{2(n-2)}+G_2^{\alpha\beta\gamma}(\xi,\tau))^{-1} \, G_{1\xi}(\xi,\tau)^2 \, H(\xi,\tau), \\ 
E_3(\xi,\tau) &= -(\sqrt{2(n-2)}+G_2^{\alpha\beta\gamma}(\xi,\tau))^{-1} \, (G_{1\xi}(\xi,\tau) + G_{2\xi}^{\alpha\beta\gamma}(\xi,\tau)) \, H_\xi(\xi,\tau) \\ 
E_4(\xi,\tau) &= (n-1) \, \bigg [ \frac{G_{1\xi}(0,\tau)}{\sqrt{2(n-2)}+G_1(0,\tau)} - \int_0^\xi \frac{G_{1\xi}(\xi',\tau)^2}{(\sqrt{2(n-2)}+G_1(\xi',\tau))^2} \, d\xi' \bigg ] \, H_\xi(\xi,\tau), \\ 
E_5(\xi,\tau) &= (n-1) \, G_{2\xi}^{\alpha\beta\gamma}(\xi,\tau) \, \frac{H_\xi(0,\tau)}{\sqrt{2(n-2)}+G_1(0,\tau)} \\ &- (n-1) \, G_{2\xi}^{\alpha\beta\gamma}(\xi,\tau) \, \frac{G_{2\xi}^{\alpha\beta\gamma}(0,\tau) \, H(0,\tau)}{(\sqrt{2(n-2)}+G_1(0)) (\sqrt{2(n-2)}+G_2^{\alpha\beta\gamma}(0,\tau))}, \\ 
E_6(\xi,\tau) &= (n-1) \, G_{2\xi}^{\alpha\beta\gamma}(\xi,\tau) \, \bigg [ -\int_0^\xi \frac{(G_{1\xi}(\xi',\tau)+G_{2\xi}^{\alpha\beta\gamma}(\xi',\tau)) \, H_\xi(\xi',\tau)}{(\sqrt{2(n-2)}+G_2^{\alpha\beta\gamma}(\xi',\tau))^2} \, d\xi' \\ 
&+ \int_0^\xi \frac{(2\sqrt{2(n-2)}+G_1(\xi',\tau)+G_2^{\alpha\beta\gamma}(\xi',\tau)) \, H(\xi',\tau) \, G_{1\xi}(\xi',\tau)^2}{(\sqrt{2(n-2)}+G_1(\xi',\tau))^2 (\sqrt{2(n-2)}+G_2^{\alpha\beta\gamma}(\xi',\tau))^2} \, d\xi' \bigg ]. 
\end{align*}

Consequently, the function $H_{\mathcal{C}}(\xi,\tau) = \chi_{\mathcal{C}}((-\tau)^{-\frac{1}{2}} \xi) \, H(\xi,\tau)$ satisfies 
\[H_{\mathcal{C},\tau}(\xi,\tau) = H_{\mathcal{C},\xi\xi}(\xi,\tau) - \frac{1}{2} \, \xi \, H_{\mathcal{C},\xi}(\xi,\tau) + H_{\mathcal{C}}(\xi,\tau) + \sum_{k=1}^{10} E_{\mathcal{C},k}(\xi,\tau),\] 
where 
\begin{align*} 
E_{\mathcal{C},1}(\xi,\tau) &= \Big [ (n-2)(\sqrt{2(n-2)}+G_1(\xi,\tau))^{-1} (\sqrt{2(n-2)}+G_2^{\alpha\beta\gamma}(\xi,\tau))^{-1} - \frac{1}{2} \Big ] \, H_{\mathcal{C}}(\xi,\tau), \\ 
E_{\mathcal{C},2}(\xi,\tau) &= (\sqrt{2(n-2)}+G_1(\xi,\tau))^{-1} (\sqrt{2(n-2)}+G_2^{\alpha\beta\gamma}(\xi,\tau))^{-1} \, G_{1\xi}(\xi,\tau)^2 \, H_{\mathcal{C}}(\xi,\tau), \\ 
E_{\mathcal{C},3}(\xi,\tau) &= -(\sqrt{2(n-2)}+G_2^{\alpha\beta\gamma}(\xi,\tau))^{-1} \, (G_{1\xi}(\xi,\tau) + G_{2\xi}^{\alpha\beta\gamma}(\xi,\tau)) \, H_{\mathcal{C},\xi}(\xi,\tau) \\ 
E_{\mathcal{C},4}(\xi,\tau) &= (n-1) \, \bigg [ \frac{G_{1\xi}(0,\tau)}{\sqrt{2(n-2)}+G_1(0,\tau)} - \int_0^\xi \frac{G_{1\xi}(\xi',\tau)^2}{(\sqrt{2(n-2)}+G_1(\xi',\tau))^2} \, d\xi' \bigg ] \, H_{\mathcal{C},\xi}(\xi,\tau), \\ 
E_{\mathcal{C},5}(\xi,\tau) &= (n-1) \, \chi_{\mathcal{C}}((-\tau)^{-\frac{1}{2}} \xi) \, G_{2\xi}^{\alpha\beta\gamma}(\xi,\tau) \, \frac{H_\xi(0,\tau)}{\sqrt{2(n-2)}+G_1(0,\tau)} \\ &- (n-1) \, \chi_{\mathcal{C}}((-\tau)^{-\frac{1}{2}} \xi) \, G_{2\xi}^{\alpha\beta\gamma}(\xi,\tau) \, \frac{G_{2\xi}^{\alpha\beta\gamma}(0,\tau) \, H(0,\tau)}{(\sqrt{2(n-2)}+G_1(0)) (\sqrt{2(n-2)}+G_2^{\alpha\beta\gamma}(0,\tau))}, \\ 
E_{\mathcal{C},6}(\xi,\tau) &= (n-1) \, \chi_{\mathcal{C}}((-\tau)^{-\frac{1}{2}} \xi) \, G_{2\xi}^{\alpha\beta\gamma}(\xi,\tau) \\ 
&\cdot \bigg [ -\int_0^\xi \frac{(G_{1\xi}(\xi',\tau)+G_{2\xi}^{\alpha\beta\gamma}(\xi',\tau)) \, H_\xi(\xi',\tau)}{(\sqrt{2(n-2)}+G_2^{\alpha\beta\gamma}(\xi',\tau))^2} \, d\xi' \\ 
&+ \int_0^\xi \frac{(2\sqrt{2(n-2)}+G_1(\xi',\tau)+G_2^{\alpha\beta\gamma}(\xi',\tau)) \, H(\xi',\tau) \, G_{1\xi}(\xi',\tau)^2}{(\sqrt{2(n-2)}+G_1(\xi',\tau))^2 (\sqrt{2(n-2)}+G_2^{\alpha\beta\gamma}(\xi',\tau))^2} \, d\xi' \bigg ], \\ 
E_{\mathcal{C},7}(\xi,\tau) &= (\sqrt{2(n-2)}+G_2^{\alpha\beta\gamma}(\xi,\tau))^{-1} \, (G_{1\xi}(\xi,\tau) + G_{2\xi}^{\alpha\beta\gamma}(\xi,\tau)) \\ 
&\cdot (-\tau)^{-\frac{1}{2}} \, \chi_{\mathcal{C}}'((-\tau)^{-\frac{1}{2}} \xi) \, H(\xi,\tau), \\ 
E_{\mathcal{C},8}(\xi,\tau) &= -(n-1) \, \bigg [ \frac{G_{1\xi}(0,\tau)}{\sqrt{2(n-2)}+G_1(0,\tau)} - \int_0^\xi \frac{G_{1\xi}(\xi',\tau)^2}{(\sqrt{2(n-2)}+G_1(\xi',\tau))^2} \, d\xi' \bigg ] \\ 
&\cdot (-\tau)^{-\frac{1}{2}} \, \chi_{\mathcal{C}}'((-\tau)^{-\frac{1}{2}} \xi) \, H(\xi,\tau), \\ 
E_{\mathcal{C},9}(\xi,\tau) &= (-\tau)^{-1} \, \chi_{\mathcal{C}}''((-\tau)^{-\frac{1}{2}} \xi) \, H(\xi,\tau) \\ 
&+ \frac{1}{2} \, (-\tau)^{-\frac{3}{2}} \xi \, \chi_{\mathcal{C}}'((-\tau)^{-\frac{1}{2}} \xi) \, H(\xi,\tau), \\ 
E_{\mathcal{C},10}(\xi,\tau) &= -2 \, (-\tau)^{-\frac{1}{2}} \, \frac{\partial}{\partial \xi} \big [ \chi_{\mathcal{C}}'((-\tau)^{-\frac{1}{2}} \xi) \, H(\xi,\tau) \big ] \\ 
&+ \frac{1}{2} \, (-\tau)^{-\frac{1}{2}} \xi \, \chi_{\mathcal{C}}'((-\tau)^{-\frac{1}{2}} \xi) \, H(\xi,\tau).
\end{align*}

In the following, we will estimate the terms $\sum_{k=1}^6 \|E_{\mathcal{C},k}\|_{\mathcal{H},\infty,\tau_*}$ and $\sum_{k=7}^{10} \|E_{\mathcal{C},k}\|_{\mathcal{D}^*,\infty,\tau_*}$. To that end, we need the following estimates for the functions $G_1(\xi,\tau)$ and $G_2^{\alpha\beta\gamma}(\xi,\tau)$:

\begin{proposition} 
\label{estimates.for.G_1.and.G_2}
Fix a small number $\theta > 0$ and a small number $\eta > 0$. Then there exists a small number $\varepsilon>0$ (depending on $\theta$ and $\eta$) with the following property. If the triplet $(\alpha,\beta,\gamma)$ is $\varepsilon$-admissible with respect to time $t_* = -e^{-\tau_*}$ and $-\tau_*$ is sufficiently large, then  
\begin{align*} 
&\Big | (\sqrt{2(n-2)}+G_1(\xi,\tau))^2 - 2(n-2) + (n-2)\,\frac{\xi^2-2}{2 (-\tau)} \Big | \leq \eta \, \frac{\xi^2+1}{(-\tau)}, \\ 
&\Big | (\sqrt{2(n-2)}+G_2^{\alpha\beta\gamma}(\xi,\tau))^2 - 2(n-2) + (n-2)\,\frac{\xi^2-2}{2 (-\tau)} \Big | \leq \eta \, \frac{\xi^2+1}{(-\tau)} 
\end{align*} 
and 
\begin{align*} 
&\Big | (\sqrt{2(n-2)}+G_1(\xi,\tau)) \, G_{1\xi}(\xi,\tau) + \frac{(n-2)\,\xi}{2 (-\tau)} \Big | \leq \eta \, \frac{|\xi|+1}{(-\tau)}, \\ 
&\Big | (\sqrt{2(n-2)}+G_2^{\alpha\beta\gamma}(\xi,\tau)) \, G_{2\xi}^{\alpha\beta\gamma}(\xi,\tau) + \frac{(n-2)\,\xi}{2 (-\tau)} \Big | \leq \eta \, \frac{|\xi|+1}{(-\tau)} 
\end{align*} 

for $|\xi| \leq \sqrt{4-\frac{\theta^2}{8(n-2)}} \, (-\tau)^{\frac{1}{2}}$ and $\tau \leq \tau_*$. 
\end{proposition}

\textbf{Proof.} 
This follows directly from Proposition \ref{precise.estimate.for.F}, Proposition \ref{precise.estimate.for.F_z}, and Proposition \ref{estimate.for.modified.profile}. \\

In order to estimate the term $\|E_{\mathcal{C},6}\|_{\mathcal{H},\infty,\tau_*}$, we need the following pointwise estimate: 

\begin{lemma}
\label{pointwise.estimate.for.E6}
We have 
\begin{align*} 
|E_{\mathcal{C},6}(\xi,\tau)| 
&\leq C(\theta) \, (-\tau)^{-\frac{1}{2}} \, |G_{2\xi}^{\alpha\beta\gamma}(\xi,\tau)| \, \bigg | \int_0^\xi |H_{\mathcal{C}}(\xi',\tau)| \, d\xi' \bigg | \\ 
&+ C(\theta) \, (-\tau)^{-\frac{1}{2}} \, |G_{2\xi}^{\alpha\beta\gamma}(\xi,\tau)| \, (|H_{\mathcal{C}}(\xi,\tau)|+|H(0,\tau)|).
\end{align*}
\end{lemma}

\textbf{Proof.} 
The proof is analogous to the proof of Lemma 8.4 in \cite{BDS}. \\

In order to estimate the term $\|E_{\mathcal{C},5}\|_{\mathcal{H},\infty,\tau_*}$, we need the following estimate for $H_\xi(0,\tau)$: 

\begin{lemma}
\label{derivative.of.H.at.0}
We have 
\[\sup_{\tau \leq \tau_*} \bigg ( \int_{\tau-1}^\tau H_\xi(0,\tau')^2 \, d\tau' \bigg )^{\frac{1}{2}} \leq C \, \|H_{\mathcal{C}}\|_{\mathcal{H},\infty,\tau_*} + C \sum_{k=1}^6 \|E_{\mathcal{C},k}\|_{\mathcal{H},\infty,\tau_*}.\]
\end{lemma}

\textbf{Proof.} 
In the region $\{|\xi| \leq 1\}$, we have $\frac{\partial}{\partial \tau} H_{\mathcal{C}} = \mathcal{L} H_{\mathcal{C}} + \sum_{k=1}^6 E_{\mathcal{C},k}$. Using standard interior estimates for linear parabolic equations and the embedding of the Sobolev space $H^2([-1,1])$ into $C^1([-1,1])$, we obtain 
\[\sup_{\tau \leq \tau_*} \bigg ( \int_{\tau-1}^\tau H_{\mathcal{C},\xi}(0,\tau')^2 \, d\tau' \bigg )^{\frac{1}{2}} \leq C \, \|H_{\mathcal{C}}\|_{\mathcal{H},\infty,\tau_*} + C \sum_{k=1}^6 \|E_{\mathcal{C},k}\|_{\mathcal{H},\infty,\tau_*}.\]
Since $H_{\mathcal{C},\xi}(0,\tau) = H_\xi(0,\tau)$, the assertion follows. \\

\begin{lemma}
\label{E1.E6}
We have 
\[\sum_{k=1}^6 \|E_{\mathcal{C},k}\|_{\mathcal{H},\infty,\tau_*} \leq C(\theta) \, (-\tau_*)^{-\frac{1}{2}} \, \|H_{\mathcal{C}}\|_{\mathcal{D},\infty,\tau_*}.\]
\end{lemma}

\textbf{Proof.} The proof is analogous to the proof of Lemma 8.6 in \cite{BDS}.

\begin{lemma}
\label{E7.E10}
We have 
\begin{align*} 
&\sum_{k=7}^9 \|E_{\mathcal{C},k}\|_{\mathcal{H},\infty,\tau_*} + \|E_{\mathcal{C},10}\|_{\mathcal{D}^*,\infty,\tau_*} \\ 
&\leq C(\theta) \, (-\tau_*)^{-\frac{1}{2}} \, \Big \| H \, 1_{\{\sqrt{4-\frac{\theta^2}{2(n-2)}} \, (-\tau)^{\frac{1}{2}} \leq |\xi| \leq \sqrt{4-\frac{\theta^2}{4(n-2)}} \, (-\tau)^{\frac{1}{2}}\}} \Big \|_{\mathcal{H},\infty,\tau_*}. 
\end{align*}
\end{lemma}

\textbf{Proof.} 
Using Proposition \ref{boundedness.of.operators}, we obtain 
\[\|E_{\mathcal{C},10}\|_{\mathcal{D}^*,\infty,\tau_*} \leq C \, (-\tau_*)^{-\frac{1}{2}} \, \|\chi_{\mathcal{C}}'((-\tau)^{-\frac{1}{2}} \xi) \, H\|_{\mathcal{H},\infty,\tau_*}.\] 
This gives the desired estimate for $E_{\mathcal{C},10}$. The estimates for $E_{\mathcal{C},7}$, $E_{\mathcal{C},8}$, and $E_{\mathcal{C},9}$ follow directly from the respective definitions. This completes the proof of Lemma \ref{E7.E10}. \\

We now complete the proof of Proposition \ref{estimate.for.difference.in.cylindrical.region}. To that end, we apply Proposition \ref{linear.estimate} to the function $H_{\mathcal{C}}$. Since $P_+ H_{\mathcal{C}}(\tau_*) = 0$, we obtain 
\[\sup_{\tau \leq \tau_*} \|\hat{H}_{\mathcal{C}}(\tau)\|_{\mathcal{H}} + \Lambda^{-1} \, \|\hat{H}_{\mathcal{C}}\|_{\mathcal{D},\infty,\tau_*} \leq \Lambda \sum_{k=1}^{10} \|E_{\mathcal{C},k}(\xi,\tau)\|_{\mathcal{D}^*,\infty,\tau_*}\] 
by Proposition \ref{estimate.for.difference.in.cylindrical.region}. We use Lemma \ref{E1.E6} and Lemma \ref{E7.E10} to estimate the terms on the right hand side. This gives 
\begin{align*} 
&\sup_{\tau \leq \tau_*} \|\hat{H}_{\mathcal{C}}(\tau)\|_{\mathcal{H}} + \Lambda^{-1} \, \|\hat{H}_{\mathcal{C}}\|_{\mathcal{D},\infty,\tau_*} \\ 
&\leq C(\theta) \, (-\tau_*)^{-\frac{1}{2}} \, \|\hat{H}_{\mathcal{C}}\|_{\mathcal{D},\infty,\tau_*} + C(\theta) \, (-\tau_*)^{-\frac{1}{2}} \, \|P_0 H_{\mathcal{C}}\|_{\mathcal{D},\infty,\tau_*} \\ 
&+ C(\theta) \, (-\tau_*)^{-\frac{1}{2}} \, \Big \| H \, 1_{\{\sqrt{4-\frac{\theta^2}{2(n-2)}} \, (-\tau)^{\frac{1}{2}} \leq |\xi| \leq \sqrt{4-\frac{\theta^2}{4(n-2)}} \, (-\tau)^{\frac{1}{2}}\}} \Big \|_{\mathcal{H},\infty,\tau_*}. 
\end{align*}
If $-\tau_*$ is sufficiently large, the first term on the right hand side can be absorbed into the left hand side. This completes the proof of Proposition \ref{estimate.for.difference.in.cylindrical.region}. \\

\section{Analysis of the overlap region and proof of Proposition \ref{neutral.mode.dominates}}

\label{overlap.region}

In this section, we give the proof of Proposition \ref{neutral.mode.dominates}. We remind the reader that $\theta$ is chosen as in Proposition \ref{estimate.for.difference.in.tip.region}. We also recall that $\chi_{\mathcal{C}}$ is a smooth cutoff, which satisfies $\chi_{\mathcal{C}} = 1$ on  $[0, \sqrt{4 - \frac{\theta^2}{2(n-2)}}]$ and $\chi_{\mathcal{C}} = 0$ on $[\sqrt{4 - \frac{\theta^2}{4(n-2)}}, \infty)$.  We also assume $\chi_{\mathcal{C}}$ is monotone decreasing on $[0, \infty)$. As before, we write $H$, $H_{\mathcal{C}}$, $\hat{H}_{\mathcal{C}}$, and $a$ instead of $H^{\alpha\beta\gamma}$, $H_{\mathcal{C}}^{\alpha\beta\gamma}$, $\hat{H}_{\mathcal{C}}^{\alpha\beta\gamma}$, and $a^{\alpha\beta\gamma}$. We begin by recalling the following elementary lemma from \cite{BDS}:

\begin{lemma}[Lemma 9.1 in \cite{BDS}]
\label{weighted.Poincare.inequality.cylindrical.region}
Assume that $4 \leq L_1 < L_2 < L_3$. Then 
\begin{align*} 
L_2^2 \int_{\{L_2 \leq \xi \leq L_3\}} e^{-\frac{\xi^2}{4}} \, f(\xi)^2 \, d\xi 
&\leq C \int_{\{L_1 \leq \xi \leq L_3\}} e^{-\frac{\xi^2}{4}} \, f'(\xi)^2 \, d\xi \\ 
&+ C \, (L_2-L_1)^{-2} \int_{\{L_1 \leq \xi \leq L_2\}} e^{-\frac{\xi^2}{4}} \, f(\xi)^2 \, d\xi, 
\end{align*}
where $C$ is a numerical constant that is independent of $L_1$, $L_2$, $L_3$, and $f$.
\end{lemma} 

The following lemma relates the function $H(\xi,\tau)$ to the function $W_+(\rho,\tau)$:

\begin{lemma}
\label{pointwise.estimate.in.overlap.region}
If we choose $-\tau_*$ sufficiently large (depending on $\theta$), then 
\[\big | H_\xi(\xi,\tau) + W_+(\sqrt{2(n-2)}+G_1(\xi,\tau),\tau) \big | \leq C(\theta) \, |H(\xi,\tau)|\] 
provided that $\sqrt{4-400 \, \frac{\theta^2}{n-2}} \, (-\tau)^{\frac{1}{2}} \leq \xi \leq \sqrt{4-\frac{\theta^2}{100(n-2)}} \, (-\tau)^{\frac{1}{2}}$ and $\tau \leq \tau_*$. 
\end{lemma}

\textbf{Proof.} 
The proof of this lemma  is analogous to the proof of Lemma 9.2 in \cite{BDS}. Replacing $\theta$ by $\theta/\sqrt{n-2}$ in higher dimensions accounts for the minor changes in Proposition \ref{precise.estimate.for.F} and Proposition \ref{estimate.for.modified.profile} in higher dimensions. The details are left to the reader.

\begin{lemma} 
\label{overlap.1}
We have 
\begin{align*} 
&(-\tau) \int_{\{\sqrt{4-\frac{\theta^2}{2(n-2)}} \, (-\tau)^{\frac{1}{2}} \leq \xi \leq \sqrt{4-\frac{\theta^2}{4(n-2)}} \, (-\tau)^{\frac{1}{2}}\}} e^{-\frac{\xi^2}{4}} \, H(\xi,\tau)^2 \, d\xi \\ 
&\leq C(\theta) \, (-\tau)^{-\frac{1}{2}} \int_{\frac{\theta}{4}}^\theta V_{1+}(\rho,\tau)^{-2} \, W_+(\rho,\tau)^2 \, e^{\mu_+(\rho,\tau)} \, d\rho \\ 
&+ C(\theta) \int_{\{\sqrt{4-\frac{\theta^2}{n-2}} \, (-\tau)^{\frac{1}{2}} \leq \xi \leq \sqrt{4-\frac{\theta^2}{2(n-2)}} \, (-\tau)^{\frac{1}{2}}\}} e^{-\frac{\xi^2}{4}} \, H(\xi,\tau)^2 \, d\xi,
\end{align*} 
provided that $\tau \leq \tau_*$ and $-\tau_*$ is sufficiently large.
\end{lemma}

\textbf{Proof.} 
The proof of this lemma is analogous to the proof of Lemma 9.3 in \cite{BDS}. The first step, as in \cite{BDS} is to apply Lemma \ref{weighted.Poincare.inequality.cylindrical.region} with $L_1 = \sqrt{4-\frac{\theta^2}{n-2}} \, (-\tau)^{\frac{1}{2}}$, $L_2 = \sqrt{4-\frac{\theta^2}{2(n-2)}} \, (-\tau)^{\frac{1}{2}}$, $L_3 = \sqrt{4-\frac{\theta^2}{4(n-2)}} \, (-\tau)^{\frac{1}{2}}$, and $f(\xi) = H(\xi,\tau)$. Besides replacing $\theta$ by $\theta/\sqrt{n-2}$, the remainder of the proof goes through unchanged. 


\begin{lemma} 
\label{overlap.2}
We have 
\begin{align*} 
&(-\tau)^{-\frac{1}{2}} \int_\theta^{2\theta} V_{1+}(\rho,\tau)^{-2} \, W_+(\rho,\tau)^2 \, e^{\mu_+(\rho,\tau)} \, d\rho \\ 
&\leq C(\theta) \int_{\{\sqrt{4-\frac{16\theta^2}{n-2}} \, (-\tau)^{\frac{1}{2}} \leq \xi \leq \sqrt{4-\frac{\theta^2}{n-2}} \, (-\tau)^{\frac{1}{2}}\}} e^{-\frac{\xi^2}{4}} \, (H_\xi(\xi,\tau)^2 + H(\xi,\tau)^2) \, d\xi 
\end{align*} 
provided that $\tau \leq \tau_*$ and $-\tau_*$ is sufficiently large.
\end{lemma}

\textbf{Proof.} 
The proof is analogous to the proof of Lemma 9.4 in \cite{BDS}. \\

\begin{proposition}
\label{control.overlap.region}
We have 
\begin{align*} 
&\sup_{\tau \leq \tau_*} (-\tau) \int_{\tau-1}^\tau \int_{\{\sqrt{4-\frac{\theta^2}{2(n-2)}} \, (-\tau')^{\frac{1}{2}} \leq |\xi| \leq \sqrt{4-\frac{\theta^2}{4(n-2)}} \, (-\tau')^{\frac{1}{2}}\}} e^{-\frac{\xi^2}{4}} \, H(\xi,\tau')^2 \, d\xi \, d\tau' \\ 
&\leq C(\theta) \sup_{\tau \leq \tau_*} \int_{\tau-1}^\tau \int_{\mathbb{R}} e^{-\frac{\xi^2}{4}} \, (H_{\mathcal{C},\xi}(\xi,\tau')^2 + H_{\mathcal{C}}(\xi,\tau')^2) \, d\xi \, d\tau'. 
\end{align*} 
\end{proposition}

\textbf{Proof.} 
The proof is analogous to the proof of Proposition 9.5 in \cite{BDS}. \\


After these preparations, we now finish the proof of Proposition \ref{neutral.mode.dominates}. Using Proposition \ref{control.overlap.region}, we obtain
\begin{align*} 
&\sup_{\tau \leq \tau_*} (-\tau) \int_{\tau-1}^\tau \int_{\{\sqrt{4-\frac{\theta^2}{2(n-2)}} \, (-\tau')^{\frac{1}{2}} \leq |\xi| \leq \sqrt{4-\frac{\theta^2}{4(n-2)}} \, (-\tau')^{\frac{1}{2}}\}} e^{-\frac{\xi^2}{4}} \, H(\xi,\tau')^2 \, d\xi \, d\tau' \\ 
&\leq C(\theta) \sup_{\tau \leq \tau_*} \int_{\tau-1}^\tau a(\tau')^2 \, d\tau' \\ 
&+ C(\theta) \sup_{\tau \leq \tau_*} \int_{\tau-1}^\tau \int_{\mathbb{R}} e^{-\frac{\xi^2}{4}} \, (\hat{H}_{\mathcal{C},\xi}(\xi,\tau')^2 + \hat{H}_{\mathcal{C}}(\xi,\tau')^2) \, d\xi \, d\tau'. 
\end{align*} 
Combining this estimate with Proposition \ref{estimate.for.difference.in.cylindrical.region} gives 
\begin{align*}
&(-\tau_*) \sup_{\tau \leq \tau_*} \int_{\tau-1}^\tau \int_{\mathbb{R}} e^{-\frac{\xi^2}{4}} \, (\hat{H}_{\mathcal{C},\xi}(\xi,\tau')^2 + \hat{H}_{\mathcal{C}}(\xi,\tau')^2) \, d\xi \, d\tau' \\ 
&\leq C(\theta) \sup_{\tau \leq \tau_*} \int_{\tau-1}^\tau a(\tau')^2 \, d\tau'
\end{align*} 
if $-\tau_*$ is chosen sufficiently large. This completes the proof of Proposition \ref{neutral.mode.dominates}. \\

\section{Analysis of the neutral mode and proof of Proposition \ref{ode.for.a}}

\label{analysis.of.neutral.mode}

In this final section, we give the proof of Proposition \ref{ode.for.a}. 

\begin{lemma}
\label{control.overlap.region.2}
We have 
\begin{align*} 
&\sup_{\tau \leq \tau_*} (-\tau) \int_{\tau-1}^\tau \int_{\{\sqrt{4-\frac{\theta^2}{2(n-2)}} \, (-\tau')^{\frac{1}{2}} \leq |\xi| \leq \sqrt{4-\frac{\theta^2}{4(n-2)}} \, (-\tau')^{\frac{1}{2}}\}} e^{-\frac{\xi^2}{4}} \, H(\xi,\tau')^2 \, d\xi \, d\tau' \\ 
&\leq C(\theta) \sup_{\tau \leq \tau_*} \int_{\tau-1}^\tau a(\tau')^2 \, d\tau'. 
\end{align*} 
\end{lemma}

\textbf{Proof.} 
This follows by combining Proposition \ref{neutral.mode.dominates} and Proposition \ref{control.overlap.region}. \\

We next establish an improved version of Lemma \ref{derivative.of.H.at.0}: 

\begin{lemma}
\label{derivative.of.H.at.0.improved.version}
We have 
\[(-\tau_*) \sup_{\tau \leq \tau_*} \int_{\tau-1}^\tau H_\xi(0,\tau')^2 \, d\tau' \leq C(\theta) \sup_{\tau \leq \tau_*} \int_{\tau-1}^\tau a(\tau')^2 \, d\tau'.\] 
\end{lemma}

\textbf{Proof.} The proof of this lemma is analogous to the proof of Lemma 10.2 in \cite{BDS}. \\


After these preparations, we now study the evolution of the function $a(\tau)$. Using the evolution equation $\frac{\partial}{\partial \tau} H_{\mathcal{C}} = \mathcal{L} H_{\mathcal{C}} + \sum_{k=1}^{10} E_{\mathcal{C},k}$, we obtain 
\[\frac{d}{d\tau} a(\tau) = \sum_{k=1}^{10} I_k(\tau),\] 
where 
\[I_k(\tau) = \frac{1}{16\sqrt{2(n-2)\pi}} \int_{\mathbb{R}} e^{-\frac{\xi^2}{4}} \, (\xi^2-2) \, E_{\mathcal{C},k}(\xi,\tau) \, d\xi.\] 
In the remainder of this section, we estimate the terms $I_k(\tau)$. \\

\begin{lemma} 
\label{I1}
Let $\delta > 0$ be given. If $-\tau_*$ is sufficiently large (depending on $\delta$), then  
\[\sup_{\tau \leq \tau_*} (-\tau) \int_{\tau-1}^\tau |I_1(\tau') - (-\tau')^{-1} \, a(\tau')| \, d\tau' \leq \delta \, \sup_{\tau \leq \tau_*} \bigg ( \int_{\tau-1}^\tau a(\tau')^2 \, d\tau' \bigg )^{\frac{1}{2}}.\] 
\end{lemma}

\textbf{Proof.} The proof of this lemma is analogous to the proof of Lemma 10.3 in \cite{BDS}. \\


\begin{lemma} 
\label{I2}
Let $\delta > 0$ be given. If $-\tau_*$ is sufficiently large (depending on $\delta$), then  
\[\sup_{\tau \leq \tau_*} (-\tau) \int_{\tau-1}^\tau |I_2(\tau')| \, d\tau' \leq \delta \, \sup_{\tau \leq \tau_*} \bigg ( \int_{\tau-1}^\tau a(\tau')^2 \, d\tau' \bigg )^{\frac{1}{2}}.\] 
\end{lemma}

\textbf{Proof.} The proof of this lemma is analogous to the proof of Lemma 10.4 in \cite{BDS}.\\

\begin{lemma} 
\label{I3}
Let $\delta > 0$ be given. If $-\tau_*$ is sufficiently large (depending on $\delta$), then  
\[\sup_{\tau \leq \tau_*} (-\tau) \int_{\tau-1}^\tau |I_3(\tau') - (-\tau')^{-1} \, a(\tau')| \, d\tau' \leq \delta \, \sup_{\tau \leq \tau_*} \bigg ( \int_{\tau-1}^\tau a(\tau')^2 \, d\tau' \bigg )^{\frac{1}{2}}.\] 
\end{lemma}

\textbf{Proof.} The proof of this lemma is analogous to the proof of Lemma 10.5 in \cite{BDS}.\\


\begin{lemma} 
\label{I4}
Let $\delta > 0$ be given. If $-\tau_*$ is sufficiently large (depending on $\delta$), then 
\[\sup_{\tau \leq \tau_*} (-\tau) \int_{\tau-1}^\tau |I_4(\tau')| \, d\tau' \leq \delta \, \sup_{\tau \leq \tau_*} \bigg ( \int_{\tau-1}^\tau a(\tau')^2 \, d\tau' \bigg )^{\frac{1}{2}}.\] 
\end{lemma}

\textbf{Proof.} The proof of this lemma is analogous to the proof of Lemma 10.6 in \cite{BDS}. \\

\begin{lemma} 
\label{I5}
Let $\delta > 0$ be given. If $-\tau_*$ is sufficiently large (depending on $\delta$), then 
\[\sup_{\tau \leq \tau_*} (-\tau) \int_{\tau-1}^\tau |I_5(\tau')| \, d\tau' \leq \delta \, \sup_{\tau \leq \tau_*} \bigg ( \int_{\tau-1}^\tau a(\tau')^2 \, d\tau' \bigg )^{\frac{1}{2}}.\] 
\end{lemma}

\textbf{Proof.} The proof of this lemma is analogous to the proof of Lemma 10.7 in \cite{BDS}. \\

\begin{lemma} 
\label{I6}
Let $\delta > 0$ be given. If $-\tau_*$ is sufficiently large (depending on $\delta$), then 
\[\sup_{\tau \leq \tau_*} (-\tau) \int_{\tau-1}^\tau |I_6(\tau')| \, d\tau' \leq \delta \, \sup_{\tau \leq \tau_*} \bigg ( \int_{\tau-1}^\tau a(\tau')^2 \, d\tau' \bigg )^{\frac{1}{2}}.\] 
\end{lemma}

\textbf{Proof.} The proof of this lemma is analogous to the proof of Lemma 10.8 in \cite{BDS}. \\


\begin{lemma} 
\label{I7.I10}
Let $\delta > 0$ be given. If $-\tau_*$ is sufficiently large (depending on $\delta$), then 
\[\sup_{\tau \leq \tau_*} (-\tau) \int_{\tau-1}^\tau \sum_{k=7}^{10} |I_k(\tau')| \, d\tau' \leq \delta \, \sup_{\tau \leq \tau_*} \bigg ( \int_{\tau-1}^\tau a(\tau')^2 \, d\tau' \bigg )^{\frac{1}{2}}.\] 
\end{lemma}

\textbf{Proof.} The proof of this lemma is analogous to the proof of Lemma 10.9 in \cite{BDS}.\\

Proposition \ref{ode.for.a} follows immediately from Lemma \ref{I1} -- Lemma \ref{I7.I10} together with the identity $\frac{d}{d\tau} a(\tau) = \sum_{k=1}^{10} I_k(\tau)$. 

\appendix
\section{The Bryant soliton}

\label{properties.of.Bryant.soliton}

In \cite{Bryant} Bryant showed that up to constant multiples, there is only one complete, steady, rotationally symmetric soliton in dimension three that is not flat. It has positive sectional curvature. The maximum scalar curvature is equal to $1$, and is attained at the center of rotation. The complete soliton can be written in the form $g = dz \otimes dz + B(z)^2 \, g_{S^n}$, where $z$ is the distance from the center of rotation. For large $z$, the metric has the following asymptotics: the aperature $B(z)$ has leading order term $\sqrt{2(n-2)z}$, the orbital sectional curvature $K_{\text{\rm orb}}$ has leading order term $\frac{1}{2(n-2)z}$, and the radial sectional curvature $K_{\text{\rm rad}}$ has leading order term $\frac{1}{4z^2}$.  

Sometimes it is more convenient to write the metric in the form $\Phi(r)^{-1} \, dr^2+ r^2 \, g_{S^n}$, where the function $\Phi(r)$ is defined by $\Phi(B(z)) = \big ( \frac{d}{dz} B(z) \big )^2$. The function $\Phi(r)$ is known to satisfy the equation
\[\Phi(r) \Phi''(r) - \frac{1}{2} \, \Phi'(r)^2 + \frac{n-2-\Phi(r)}{r}\, \Phi'(r) + \frac{2(n-2)}{r^2} \Phi(r) (1 - \Phi(r)) = 0.\] 
The orbital and radial sectional curvatures are given by $K_{\text{\rm orb}} = \frac{1}{r^2} \, (1-\Phi(r))$ and $K_{\text{\rm rad}} = -\frac{1}{2r} \, \Phi'(r)$. It is known that $\Phi(r)$ has the following asymptotics.  Near $r= 0$, $\Phi$ is smooth and has the asymptotic expansion
\[\Phi(r) = 1 + b_0 \, r^2 + o(r^2),\]
where $b_0$ is a negative constant (since the curvature is positive). As, $r \to \infty$, $\Phi$ is smooth and has the asymptotic expansion
\[\Phi(r) = c_0 \, r^{-2} + \frac{5-n}{n-2} c_0^2 \, r^{-4} + o(r^{-4}),\]
where $c_0$ is a positive constant.

We will next find (for the convenience of the reader) the exact values of the constants $b_0$ and $c_0$ in the above asymptotics for the {\em Bryant soliton of maximal scalar curvature one}.  

Recall that the scalar curvature is given by $R = (n-1)(n-2) K_{\text{\rm orb}} + 2(n-1) K_{\text{\rm rad}}$. The maximal scalar curvature is attained at $z = 0$, at which point $K_{\text{\rm orb}} = K_{\text{\rm rad}}$. The maximal scalar curvature being equal to $1$ is equivalent to $K_{\text{\rm orb}} = K_{\text{\rm rad}} = \frac{1}{n(n-1)}$ at $z=0$. On the other hand, the asymptotic expansion of $\Phi(r)$ gives $K_{\text{\rm orb}} = \frac{1}{r^2} \, (1-\Phi(r)) = -b_0 + o(1)$ as $r \to 0$. Consequently, $b_0 = -\frac{1}{n(n-1)}$.

Bryant's asymptotics imply that for $z$ sufficiently large, the aperture satisfies $r = (1+o(1)) \, \sqrt{2(n-2)z}$, implying that $2(n-2)z = (1+o(1)) \, r^2$. The radial sectional curvature satisfies $K_{\text{\rm rad}} = (1+o(1)) \, \frac{1}{4z^2} = (1+o(1)) \, \frac{(n-2)^2}{r^4}$ for $r$ large.  On the other hand, the asymptotic expansion of $\Phi(r)$ implies $K_{\text{\rm rad}} = -\frac{1}{2r} \, \Phi'(r) = (1+o(1)) \, c_0 \, r^{-4}$ for $r$ large. Comparing the two formulae, we conclude that $c_0=(n-2)^2$.

Summarizing the above discussion we conclude the following asymptotics for the {\em Bryant soliton with maximal scalar curvature equal to one}:
\[\Phi(r) = \begin{cases} 1 - \frac{r^2}{n(n-1)} + o(r^2) &\qquad \text{\rm as $r \to 0$,} \\ (n-2)^2\, r^{-2} + (n-2)^3(5-n)r^{-4} + o(r^{-4}) &\qquad \text{\rm as $r \to \infty$.} \end{cases}\]

\end{document}